\documentclass [10pt,twoside,reqno]{amsart}
\usepackage{graphicx}
\newtheorem{thm}{Theorem}[section]

\newtheorem{lem}[thm]{Lemma}

\theoremstyle{definition}
\newtheorem{defn}[thm]{Definition}
\theoremstyle{remark}
\newtheorem{exam}[thm]{Example}
\newtheorem{rem}[thm]{Remark}
\numberwithin{equation}{section}
\newcommand{\ie}{{\it i.e.\/}\ }

\def\Hom{{\rm Hom}}

\def\Sp{{\rm Spec}}

\def\Tr{{\rm Tr}}

\newcommand{\R}{\mathbb{R}}

\def\C{{\mathbb C}}

\def\N{{\mathbb N}}

\def\Z{{\mathbb Z}}
\def\R{{\mathbb R}}

\def\cD{{\mathcal D}}

\def\cG{{\mathcal G}}

\def\cS{{\mathcal S}}

\def\Hom{{\rm Hom}}

\def\Sp{{\rm Spec}}

\def\Tr{{\rm Tr}}

\newtheorem{theorem}{Theorem}[section]

\theoremstyle{definition}

\theoremstyle{remark}
\newtheorem{rems}[theorem]{Remarks}

\usepackage[pdftex]{hyperref}
\begin{document}

\title{Short Steps in Noncommutative Geometry }%
\author{Ahmad Zainy Al-Yasry }%
%\address{}%
%\email{}%

%\thanks{}%
%\subjclass{}%
%\keywords{}%

%\date{}%
%\dedicatory{}%
%\commby{}%
% ----------------------------------------------------------------
\begin{abstract}
Noncommutative geometry (NCG) is a branch of mathematics concerned with a geometric approach to noncommutative algebras, and with the construction of spaces that are locally presented by noncommutative algebras of functions (possibly in some generalized sense). A noncommutative algebra is an associative algebra in which the multiplication is not commutative, that is, for which $xy$ does not always equal $yx$; or more generally an algebraic structure in which one of the principal binary operations is not commutative; one also allows additional structures, e.g. topology or norm, to be possibly carried by the noncommutative algebra of functions. These notes just to start understand  what we need to study Noncommutative Geometry.
\end{abstract}
\maketitle
% ----------------------------------------------------------------
\tableofcontents

\section{Noncommutative geometry}

\subsection{Motivation}
The main motivation is to extend the commutative duality between spaces and functions to the noncommutative setting. In mathematics, spaces, which are geometric in nature, can be related to numerical functions on them. In general, such functions will form a commutative ring. For instance, one may take the ring $C(X)$ of continuous complex-valued functions on a topological space X. In many cases (e.g., if $X$ is a compact Hausdorff space), we can recover $X$ from $C(X)$, and therefore it makes some sense to say that X has commutative topology.

More specifically, in topology, compact Hausdorff topological spaces can be reconstructed from the Banach algebra of functions on the space (Gel'fand-Neimark). In commutative algebraic geometry, algebraic schemes are locally prime spectra of commutative unital rings (A. Grothendieck), and schemes can be reconstructed from the categories of quasicoherent sheaves of modules on them (P. Gabriel-A. Rosenberg). For Grothendieck topologies, the cohomological properties of a site are invariant of the corresponding category of sheaves of sets viewed abstractly as a topos (A. Grothendieck). In all these cases, a space is reconstructed from the algebra of functions or its categorified version�some category of sheaves on that space.

Functions on a topological space can be multiplied and added pointwise hence they form a commutative algebra; in fact these operations are local in the topology of the base space, hence the functions form a sheaf of commutative rings over the base space.

The dream of noncommutative geometry is to generalize this duality to the duality between

\begin{itemize}
  \item noncommutative algebras, or sheaves of noncommutative algebras, or sheaf-like noncommutative algebraic or operator-algebraic structures
  \item and geometric entities of certain kind,
\end{itemize}

and interact between the algebraic and geometric description of those via this duality.

Regarding that the commutative rings correspond to usual affine schemes, and commutative $C^*$-algebras to usual topological spaces, the extension to noncommutative rings and algebras requires non-trivial generalization of topological spaces, as "non-commutative spaces". For this reason, some talk about non-commutative topology, though the term also has other meanings.

\subsection{Noncommutative $C^*$-algebras, von Neumann algebras}
(The formal duals of) non-commutative $C^*$-algebras are often now called non-commutative spaces. This is by analogy with the Gelfand representation, which shows that commutative $C^*$-algebras are dual to locally compact Hausdorff spaces. In general, one can associate to any $C^*$-algebra $S$ a topological space $S^�$; see spectrum of a $C^*$-algebra.

For the duality between $\sigma$-finite measure spaces and commutative von Neumann algebras, noncommutative von Neumann algebras are called non-commutative measure spaces.

\subsection{Noncommutative differentiable manifolds}

A smooth Riemannian manifold $M$ is a topological space with a lot of extra structure. From its algebra of continuous functions $C(M)$ we only recover $M$ topologically. The algebraic invariant that recovers the Riemannian structure is a spectral triple. It is constructed from a smooth vector bundle $E$ over $M$, e.g. the exterior algebra bundle. The Hilbert space $L^2(M,E)$ of square integrable sections of $E$ carries a representation of $C(M)$ by multiplication operators, and we consider an unbounded operator $D$ in $L^2(M,E)$ with compact resolvent (e.g. the signature operator), such that the commutators $[D,f]$ are bounded whenever $f$ is smooth. A recent deep theorem states that $M$ as a Riemannian manifold can be recovered from this data.

This suggests that one might define a noncommutative Riemannian manifold as a spectral triple $(A,H,D)$, consisting of a representation of a $C^*$-algebra A on a Hilbert space $H$, together with an unbounded operator $D$ on $H$, with compact resolvent, such that $[D,a]$ is bounded for all a in some dense subalgebra of $A$. Research in spectral triples is very active, and many examples of noncommutative manifolds have been constructed.
\subsection{Noncommutative affine and projective schemes}

In analogy to the duality between affine schemes and commutative rings, we define a category of noncommutative affine schemes as the dual of the category of associative unital rings. There are certain analogues of Zariski topology in that context so that one can glue such affine schemes to more general objects.

There are also generalizations of the Cone and of the Proj of a commutative graded ring, mimicking a Serre's theorem on Proj. Namely the category of quasicoherent sheaves of O-modules on a Proj of a commutative graded algebra is equivalent to the category of graded modules over the ring localized on Serre's subcategory of graded modules of finite length; there is also analogous theorem for coherent sheaves when the algebra is Noetherian. This theorem is extended as a definition of noncommutative projective geometry by Michael Artin and J. J. Zhang who add also some general ring-theoretic conditions (e.g. Artin-Schelter regularity).

Many properties of projective schemes extend to this context. For example, there exist an analog of the celebrated Serre duality for noncommutative projective schemes of Artin and Zhang. A. L. Rosenberg has created a rather general relative concept of noncommutative quasicompact scheme (over a base category), abstracting the Grothendieck's study of morphisms of schemes and covers in terms of categories of quasicoherent sheaves and flat localization functors. There is also another interesting approach via localization theory, due to Fred Van Oystaeyen, Luc Willaert and Alain Verschoren, where the main concept is that of a schematic algebra
\subsection{Invariants for noncommutative spaces}
Some of the motivating questions of the theory are concerned with extending known topological invariants to formal duals of noncommutative (operator) algebras and other replacements and candidates for noncommutative spaces. One of the main starting points of the Alain Connes' direction in noncommutative geometry is his discovery of a new homology theory associated to noncommutative associative algebras and noncommutative operator algebras, namely the cyclic homology and its relations to the algebraic K-theory (primarily via Connes-Chern character map).

The theory of characteristic classes of smooth manifolds has been extended to spectral triples, employing the tools of operator K-theory and cyclic cohomology. Several generalizations of now classical index theorems allow for effective extraction of numerical invariants from spectral triples. The fundamental characteristic class in cyclic cohomology, the JLO cocycle, generalizes the classical Chern character.

\section{Algebraic geometry}

Algebraic geometry is a branch of mathematics, classically studying zeros of multivariate polynomials. Modern algebraic geometry is based on the use of abstract algebraic techniques, mainly from commutative algebra, for solving geometrical problems about these sets of zeros.

The fundamental objects of study in algebraic geometry are algebraic varieties, which are geometric manifestations of solutions of systems of polynomial equations. Examples of the most studied classes of algebraic varieties are: plane algebraic curves, which include lines, circles, parabolas, ellipses, hyperbolas, cubic curves like elliptic curves and quartic curves like lemniscates, and Cassini ovals. A point of the plane belongs to an algebraic curve if its coordinates satisfy a given polynomial equation. Basic questions involve the study of the points of special interest like the singular points, the inflection points and the points at infinity. More advanced questions involve the topology of the curve and relations between the curves given by different equations.

Algebraic geometry occupies a central place in modern mathematics and has multiple conceptual connections with such diverse fields as complex analysis, topology and number theory. Initially a study of systems of polynomial equations in several variables, the subject of algebraic geometry starts where equation solving leaves off, and it becomes even more important to understand the intrinsic properties of the totality of solutions of a system of equations, than to find a specific solution; this leads into some of the deepest areas in all of mathematics, both conceptually and in terms of technique.

In the 20th century, algebraic geometry split into several subareas.
\begin{itemize}
  \item The main stream of algebraic geometry is devoted to the study of the complex points of the algebraic varieties and more generally to the points with coordinates in an algebraically closed field.
  \item The study of the points of an algebraic variety with coordinates in the field of the rational numbers or in a number field became arithmetic geometry (or more classically Diophantine geometry), a subfield of algebraic number theory.
  \item The study of the real points of an algebraic variety is the subject of real algebraic geometry.
  \item A large part of singularity theory is devoted to the singularities of algebraic varieties.
  \item With the rise of the computers, a computational algebraic geometry area has emerged, which lies at the intersection of algebraic geometry and computer algebra. It consists essentially in developing algorithms and software for studying and finding the properties of explicitly given algebraic varieties.

\end{itemize}
Much of the development of the main stream of algebraic geometry in the 20th century occurred within an abstract algebraic framework, with increasing emphasis being placed on "intrinsic" properties of algebraic varieties not dependent on any particular way of embedding the variety in an ambient coordinate space; this parallels developments in topology, differential and complex geometry. One key achievement of this abstract algebraic geometry is Grothendieck's scheme theory which allows one to use sheaf theory to study algebraic varieties in a way which is very similar to its use in the study of differential and analytic manifolds. This is obtained by extending the notion of point: In classical algebraic geometry, a point of an affine variety may be identified, through Hilbert's Nullstellensatz, with a maximal ideal of the coordinate ring, while the points of the corresponding affine scheme are all prime ideals of this ring. This means that a point of such a scheme may be either a usual point or a subvariety. This approach also enables a unification of the language and the tools of classical algebraic geometry, mainly concerned with complex points, and of algebraic number theory. Wiles's proof of the longstanding conjecture called Fermat's last theorem is an example of the power of this approach.

\subsection{Zeros of simultaneous polynomials}

In classical algebraic geometry, the main objects of interest are the vanishing sets of collections of polynomials, meaning the set of all points that simultaneously satisfy one or more polynomial equations. For instance, the two-dimensional sphere in three-dimensional Euclidean space $R^3$ could be defined as the set of all points $(x,y,z)$ with

$$x^2+y^2+z^2-1=0$$
A "slanted" circle in $R^3$ can be defined as the set of all points $(x,y,z)$ which satisfy the two polynomial equations

$$x^2+y^2+z^2-1=0$$
$$x+y+z=0$$

\subsection{Affine varieties}
First we start with a field $k$. In classical algebraic geometry, this field was always the complex numbers $C$, but many of the same results are true if we assume only that $k$ is algebraically closed. We consider the affine space of dimension $n$ over $k$, denoted $A^n(k)$ (or more simply $A^n$, when $k$ is clear from the context). When one fixes a coordinates system, one may identify $A^n(k)$ with $k^n$. The purpose of not working with $k^n$ is to emphasize that one "forgets" the vector space structure that $k^n$ carries.

A function $f : A^n \to A^1$ is said to be polynomial (or regular) if it can be written as a polynomial, that is, if there is a polynomial $p$ in $k[x_1,...,x_n]$ such that $f(M) = p(t_1,...,t_n)$ for every point $M$ with coordinates $(t_1,...,t_n)$ in $A^n$. The property of a function to be polynomial (or regular) does not depend on the choice of a coordinate system in $A^n$.

When a coordinate system is chosen, the regular functions on the affine $n$-space may be identified with the ring of polynomial functions in $n$ variables over $k$. Therefore, the set of the regular functions on $A^n$ is a ring, which is denoted $k[A^n]$.

We say that a polynomial vanishes at a point if evaluating it at that point gives zero. Let $S$ be a set of polynomials in $k[A^n]$. The vanishing set of $S$ (or vanishing locus or zero set) is the set $V(S)$ of all points in $A^n$ where every polynomial in $S$ vanishes. In other words,

$$V(S) = \{(t_1,\dots,t_n)|\forall p\in S, p(t_1,\dots,t_n) = 0\}$$
A subset of An which is $V(S)$, for some $S$, is called an algebraic set. The $V$ stands for variety (a specific type of algebraic set to be defined below).

Given a subset $U$ of $A^n$, can one recover the set of polynomials which generate it? If $U$ is any subset of $A^n$, define $I(U)$ to be the set of all polynomials whose vanishing set contains $U$. The $I$ stands for ideal: if two polynomials $f$ and $g$ both vanish on $U$, then $f+g$ vanishes on $U$, and if $h$ is any polynomial, then $hf$ vanishes on $U$, so $I(U)$ is always an ideal of the polynomial ring $k[A^n]$.

Two natural questions to ask are:
\begin{itemize}
  \item Given a subset $U$ of $A^n$, when is $U = V(I(U))?$

  \item Given a set $S$ of polynomials, when is $S = I(V(S))?$

\end{itemize}
The answer to the first question is provided by introducing the Zariski topology, a topology on An whose closed sets are the algebraic sets, and which directly reflects the algebraic structure of $k[A^n]$. Then $U = V(I(U))$ if and only if $U$ is an algebraic set or equivalently a Zariski-closed set. The answer to the second question is given by Hilbert's Nullstellensatz. In one of its forms, it says that $I(V(S))$ is the radical of the ideal generated by $S$. In more abstract language, there is a Galois connection, giving rise to two closure operators; they can be identified, and naturally play a basic role in the theory; the example is elaborated at Galois connection.
For various reasons we may not always want to work with the entire ideal corresponding to an algebraic set $U$. Hilbert's basis theorem implies that ideals in $k[A^n]$ are always finitely generated.
An algebraic set is called irreducible if it cannot be written as the union of two smaller algebraic sets. Any algebraic set is a finite union of irreducible algebraic sets and this decomposition is unique. Thus its elements are called the irreducible components of the algebraic set. An irreducible algebraic set is also called a variety. It turns out that an algebraic set is a variety if and only if it may be defined as the vanishing set of a prime ideal of the polynomial ring.
Some authors do not make a clear distinction between algebraic sets and varieties and use irreducible variety to make the distinction when needed.
\subsection{Regular functions}
Just as continuous functions are the natural maps on topological spaces and smooth functions are the natural maps on differentiable manifolds, there is a natural class of functions on an algebraic set, called regular functions or polynomial functions. A regular function on an algebraic set V contained in An is the restriction to V of a regular function on An. For an algebraic set defined on the field of the complex numbers, the regular functions are smooth and even analytic.

It may seem unnaturally restrictive to require that a regular function always extend to the ambient space, but it is very similar to the situation in a normal topological space, where the Tietze extension theorem guarantees that a continuous function on a closed subset always extends to the ambient topological space.

Just as with the regular functions on affine space, the regular functions on $V$ form a ring, which we denote by $k[V]$. This ring is called the coordinate ring of $V$.

Since regular functions on $V$ come from regular functions on $A^n$, there is a relationship between the coordinate rings. Specifically, if a regular function on $V$ is the restriction of two functions $f$ and $g$ in $k[A^n]$, then $f - g$ is a polynomial function which is null on $V$ and thus belongs to $I(V)$. Thus $k[V]$ may be identified with $k[A^n]/I(V)$.

\subsection{Morphism of affine varieties}
Using regular functions from an affine variety to $A^1$, we can define regular maps from one affine variety to another. First we will define a regular map from a variety into affine space: Let $V$ be a variety contained in $A^n$. Choose $m$ regular functions on $V$, and call them $f_1, ..., f_m$. We define a regular map $f$ from $V$ to $A^m$ by letting $f = (f_1, ..., f_m)$. In other words, each $f_i$ determines one coordinate of the range of $f$.

If $V'$ is a variety contained in $A^m$, we say that $f$ is a regular map from $V$ to $V'$ if the range of $f$ is contained in $V'$.

The definition of the regular maps apply also to algebraic sets. The regular maps are also called morphisms, as they make the collection of all affine algebraic sets into a category, where the objects are the affine algebraic sets and the morphisms are the regular maps. The affine varieties is a subcategory of the category of the algebraic sets.

Given a regular map $g$ from $V$ to $V'$ and a regular function $f$ of $k[V']$, then $f \circ g \in k[V]$. The map $f \to f \circ g$ is a ring homomorphism from $k[V']$ to $k[V]$. Conversely, every ring homomorphism from $k[V']$ to $k[V]$ defines a regular map from $V$ to $V'$. This defines an equivalence of categories between the category of algebraic sets and the opposite category of the finitely generated reduced $k$-algebras. This equivalence is one of the starting points of scheme theory.

\subsection{Rational function and birational equivalence}
Contrarily to the preceding ones, this section concerns only varieties and not algebraic sets. On the other hand, the definitions extend naturally to projective varieties (next subsection), as an affine variety and its projective completion have the same field of functions.

If $V$ is an affine variety, its coordinate ring is an integral domain and has thus a field of fractions which is denoted $k(V)$ and called the field of the rational functions on $V$ or, shortly, the function field of $V$. Its elements are the restrictions to $V$ of the rational functions over the affine space containing $V$. The domain of a rational function $f$ is not $V$ but the complement of the subvariety (a hypersurface) where the denominator of $f$ vanishes.

Like for regular maps, one may define a rational map from a variety $V$ to a variety $V'$. Like for the regular maps, the rational maps from $V$ to $V'$ may be identified to the field homomorphisms from $k(V')$ to $k(V)$.

Two affine varieties are birationally equivalent if there are two rational functions between them which are inverse one to the other in the regions where both are defined. Equivalently, they are birationally equivalent if their function fields are isomorphic.

An affine variety is a rational variety if it is birationally equivalent to an affine space. This means that the variety admits a rational parameterization. For example, the circle of equation $x^2+y^2-1=0$ is a rational curve, as it has the parameterization
$$
x=\frac{2\,t}{1+t^2}
$$
$$y=\frac{1-t^2}{1+t^2}$$
which may also be viewed as a rational map from the line to the circle.

The problem of resolution of singularities is to know if every algebraic variety is birationally equivalent to a variety whose projective completion is nonsingular (see also smooth completion). It has been positively solved in characteristic $0$ by Heisuke Hironaka in 1964 and is yet unsolved in finite characteristic.

\subsection{Projective variety}
Just as the formulas for the roots of $2^{nd}$, $3^{rd}$ and $4^{th}$ degree polynomials suggest extending real numbers to the more algebraically complete setting of the complex numbers, many properties of algebraic varieties suggest extending affine space to a more geometrically complete projective space. Whereas the complex numbers are obtained by adding the number $i$, a root of the polynomial $x^2 + 1$, projective space is obtained by adding in appropriate points "at infinity", points where parallel lines may meet.

To see how this might come about, consider the variety $V(y - x2)$. If we draw it, we get a parabola. As $x$ goes to positive infinity, the slope of the line from the origin to the point $(x, x^2)$ also goes to positive infinity. As $x$ goes to negative infinity, the slope of the same line goes to negative infinity.

Compare this to the variety $V(y - x^3)$. This is a cubic curve. As $x$ goes to positive infinity, the slope of the line from the origin to the point $(x, x^3)$ goes to positive infinity just as before. But unlike before, as $x$ goes to negative infinity, the slope of the same line goes to positive infinity as well; the exact opposite of the parabola. So the behavior "at infinity" of $V(y - x^3)$ is different from the behavior "at infinity" of $V(y - x^2)$.

The consideration of the projective completion of the two curves, which is their prolongation "at infinity" in the projective plane, allows to quantify this difference: the point at infinity of the parabola is a regular point, whose tangent is the line at infinity, while the point at infinity of the cubic curve is a cusp. Also, both curves are rational, as they are parameterized by $x$, and Riemann-Roch theorem implies that the cubic curve must have a singularity, which must be at infinity, as all its points in the affine space are regular.

Thus many of the properties of algebraic varieties, including birational equivalence and all the topological properties, depend on the behavior "at infinity" and so it is natural to study the varieties in projective space. Furthermore, the introduction of projective techniques made many theorems in algebraic geometry simpler and sharper: For example, B\'{e}zout's theorem on the number of intersection points between two varieties can be stated in its sharpest form only in projective space. For these reasons, projective space plays a fundamental role in algebraic geometry.

Nowadays, the projective space $P^n$ of dimension n is usually defined as the set of the lines passing through a point, considered as the origin, in the affine space of dimension $n+1$, or equivalently to the set of the vector lines in a vector space of dimension $n+1$. When a coordinate system has been chosen in the space of dimension $n+1$, all the points of a line have the same set of coordinates, up to the multiplication by an element of $k$. This defines the homogeneous coordinates of a point of $P^n$ as a sequence of $n+1$ elements of the base field $k$, defined up to the multiplication by a nonzero element of $k$ (the same for the whole sequence).

Given a polynomial in $n+1$ variables, it vanishes at all the point of a line passing through the origin if and only if it is homogeneous. In this case, one says that the polynomial vanishes at the corresponding point of $P^n$. This allows to define a projective algebraic set in $P^n$ as the set $V(f_1, ..., f_k)$ where a finite set of homogeneous polynomials $\{f_1, ..., f_k\}$ vanishes. Like for affine algebraic sets, there is a bijection between the projective algebraic sets and the reduced homogeneous ideals which define them. The projective varieties are the projective algebraic sets whose defining ideal is prime. In other words, a projective variety is a projective algebraic set, whose homogeneous coordinate ring is an integral domain, the projective coordinates ring being defined as the quotient of the graded ring or the polynomials in $n+1$ variables by the homogeneous (reduced) ideal defining the variety. Every projective algebraic set may be uniquely decomposed into a finite union of projective varieties.

The only regular functions which may be defined properly on a projective variety are the constant functions. Thus this notion is not used in projective situations. On the other hand, the field of the rational functions or function field is a useful notion, which, similarly as in the affine case, is defined as the set of the quotients of two homogeneous elements of the same degree in the homogeneous coordinate ring.
\subsection{Real algebraic geometry}

The real algebraic geometry is the study of the real points of the algebraic geometry.

The fact that the field of the reals number is an ordered field may not be occulted in such a study. For example, the curve of equation $$x^2+y^2-a=0$$ is a circle if  $a>0$, but does not have any real point if  $a<0$. It follows that real algebraic geometry is not only the study of the real algebraic varieties, but has been generalized to the study of the semi-algebraic sets, which are the solutions of systems of polynomial equations and polynomial inequalities. For example, a branch of the hyperbola of equation $x y-1 = 0$ is not an algebraic variety, but is a semi-algebraic set defined by $x y-1=0$ and $x>0$ or by $x y-1=0$ and $x+y>0$.

One of the challenging problems of real algebraic geometry is the unsolved Hilbert's sixteenth problem: Decide which respective positions are possible for the ovals of a nonsingular plane curve of degree 8

\subsection{Abstract modern viewpoint}
The modern approaches to algebraic geometry redefine and effectively extend the range of basic objects in various levels of generality to schemes, formal schemes, ind-schemes, algebraic spaces, algebraic stacks and so on. The need for this arises already from the useful ideas within theory of varieties, e.g. the formal functions of Zariski can be accommodated by introducing nilpotent elements in structure rings; considering spaces of loops and arcs, constructing quotients by group actions and developing formal grounds for natural intersection theory and deformation theory lead to some of the further extensions.

Most remarkably, in late 1950s, algebraic varieties were subsumed into Alexander Grothendieck's concept of a scheme. Their local objects are affine schemes or prime spectra which are locally ringed spaces which form a category which is antiequivalent to the category of commutative unital rings, extending the duality between the category of affine algebraic varieties over a field k, and the category of finitely generated reduced k-algebras. The gluing is along Zariski topology; one can glue within the category of locally ringed spaces, but also, using the Yoneda embedding, within the more abstract category of presheaves of sets over the category of affine schemes. The Zariski topology in the set theoretic sense is then replaced by a Grothendieck topology. Grothendieck introduced Grothendieck topologies having in mind more exotic but geometrically finer and more sensitive examples than the crude Zariski topology, namely the �tale topology, and the two flat Grothendieck topologies: fppf and fpqc; nowadays some other examples became prominent including Nisnevich topology. Sheaves can be furthermore generalized to stacks in the sense of Grothendieck, usually with some additional representability conditions leading to Artin stacks and, even finer, Deligne-Mumford stacks, both often called algebraic stacks.

Sometimes other algebraic sites replace the category of affine schemes. For example, Nikolai Durov has introduced commutative algebraic monads as a generalization of local objects in a generalized algebraic geometry. Versions of a tropical geometry, of an absolute geometry over a field of one element and an algebraic analogue of Arakelov's geometry were realized in this setup.

Another formal generalization is possible to Universal algebraic geometry in which every variety of algebras has its own algebraic geometry. The term variety of algebras should not be confused with algebraic variety.

The language of schemes, stacks and generalizations has proved to be a valuable way of dealing with geometric concepts and became cornerstones of modern algebraic geometry.

Algebraic stacks can be further generalized and for many practical questions like deformation theory and intersection theory, this is often the most natural approach. One can extend the Grothendieck site of affine schemes to a higher categorical site of derived affine schemes, by replacing the commutative rings with an infinity category of differential graded commutative algebras, or of simplicial commutative rings or a similar category with an appropriate variant of a Grothendieck topology. One can also replace presheaves of sets by presheaves of simplicial sets (or of infinity groupoids). Then, in presence of an appropriate homotopic machinery one can develop a notion of derived stack as such a presheaf on the infinity category of derived affine schemes, which is satisfying certain infinite categorical version of a sheaf axiom (and to be algebraic, inductively a sequence of representability conditions). Quillen model categories, Segal categories and quasicategories are some of the most often used tools to formalize this yielding the derived algebraic geometry, introduced by the school of Carlos Simpson, including Andre Hirschowitz, Bertrand To�n, Gabrielle Vezzosi, Michel Vaqui\'{e} and others; and developed further by Jacob Lurie, Bertrand To�n, and Gabrielle Vezzosi. Another (noncommutative) version of derived algebraic geometry, using A-infinity categories has been developed from early 1990-s by Maxim Kontsevich and followers.

\subsection{Analytic Geometry}

An analytic variety is defined locally as the set of common solutions of several equations involving analytic functions. It is analogous to the included concept of real or complex algebraic variety. Any complex manifold is an analytic variety. Since analytic varieties may have singular points, not all analytic varieties are manifolds.

Modern analytic geometry is essentially equivalent to real and complex algebraic geometry, as has been shown by Jean-Pierre Serre in his paper GAGA, the name of which is French for Algebraic geometry and analytic geometry. Nevertheless, the two fields remain distinct, as the methods of proof are quite different and algebraic geometry includes also geometry in finite characteristic.

\section{Vector field}
In vector calculus, a vector field is an assignment of a vector to each point in a subset of space. A vector field in the plane, for instance, can be visualized as a collection of arrows with a given magnitude and direction each attached to a point in the plane. Vector fields are often used to model, for example, the speed and direction of a moving fluid throughout space, or the strength and direction of some force, such as the magnetic or gravitational force, as it changes from point to point.

The elements of differential and integral calculus extend to vector fields in a natural way. When a vector field represents force, the line integral of a vector field represents the work done by a force moving along a path, and under this interpretation conservation of energy is exhibited as a special case of the fundamental theorem of calculus. Vector fields can usefully be thought of as representing the velocity of a moving flow in space, and this physical intuition leads to notions such as the divergence (which represents the rate of change of volume of a flow) and curl (which represents the rotation of a flow).

In coordinates, a vector field on a domain in n-dimensional Euclidean space can be represented as a vector-valued function that associates an n-tuple of real numbers to each point of the domain. This representation of a vector field depends on the coordinate system, and there is a well-defined transformation law in passing from one coordinate system to the other. Vector fields are often discussed on open subsets of Euclidean space, but also make sense on other subsets such as surfaces, where they associate an arrow tangent to the surface at each point (a tangent vector).

More generally, vector fields are defined on differentiable manifolds, which are spaces that look like Euclidean space on small scales, but may have more complicated structure on larger scales. In this setting, a vector field gives a tangent vector at each point of the manifold (that is, a section of the tangent bundle to the manifold). Vector fields are one kind of tensor field.

\subsection{Definition}

\begin{itemize}
  \item \textbf{Vector fields on subsets of Euclidean space} Given a subset $S$ in $R^n$, a vector field is represented by a vector-valued function $V: S \to R^n$ in standard Cartesian coordinates $(x_1, ..., x_n)$. If each component of $V$ is continuous, then $V$ is a continuous vector field, and more generally $V$ is a $C^k$ vector field if each component of $V$ is $k$ times continuously differentiable.

A vector field can be visualized as assigning a vector to individual points within an n-dimensional space.

Given two $C^k$-vector fields $V, W$ defined on $S$ and a real valued $C^k$-function $f$ defined on $S$, the two operations scalar multiplication and vector addition

 $$(fV)(p) := f(p)V(p)
 (V+W)(p) := V(p) + W(p)$$
define the module of $C^k$-vector fields over the ring of $C^k$-functions.
  \item \textbf{Coordinate transformation law} In physics, a vector is additionally distinguished by how its coordinates change when one measures the same vector with respect to a different background coordinate system. The transformation properties of vectors distinguish a vector as a geometrically distinct entity from a simple list of scalars, or from a covector.

Thus, suppose that $(x_1,...,x_n)$ is a choice of Cartesian coordinates, in terms of which the components of the vector $V$ are

$$V_x = (V_{1,x},\dots,V_{n,x})$$
and suppose that $(y_1,...,y_n)$ are n functions of the $x_i$ defining a different coordinate system. Then the components of the vector $V$ in the new coordinates are required to satisfy the transformation law
\begin{equation}\label{100}
    V_{i,y} = \sum_{j=1}^n \frac{\partial y_i}{\partial x_j} V_{j,x}.
\end{equation}

Such a transformation law is called contravariant. A similar transformation law characterizes vector fields in physics: specifically, a vector field is a specification of n functions in each coordinate system subject to the transformation law (\ref{100}) relating the different coordinate systems.

Vector fields are thus contrasted with scalar fields, which associate a number or scalar to every point in space, and are also contrasted with simple lists of scalar fields, which do not transform under coordinate changes.
  \item \textbf{Vector fields on manifolds} Given a differentiable manifold $M$, a vector field on $M$ is an assignment of a tangent vector to each point in $M$. More precisely, a vector field $F$ is a mapping from $M$ into the tangent bundle $TM$ so that  $p\circ F$  is the identity mapping where $p$ denotes the projection from $TM$ to $M$. In other words, a vector field is a section of the tangent bundle.

If the manifold $M$ is smooth or analytic�that is, the change of coordinates is smooth (analytic)�then one can make sense of the notion of smooth (analytic) vector fields. The collection of all smooth vector fields on a smooth manifold $M$ is often denoted by $G(TM)$ or $C^\infty(M,TM)$ (especially when thinking of vector fields as sections); the collection of all smooth vector fields is also denoted by $\scriptstyle \mathfrak{X} (M)$ (a fraktur $"X"$).
\end{itemize}
\subsection{Example}
\begin{itemize}
  \item A vector field for the movement of air on Earth will associate for every point on the surface of the Earth a vector with the wind speed and direction for that point. This can be drawn using arrows to represent the wind; the length (magnitude) of the arrow will be an indication of the wind speed. A "high" on the usual barometric pressure map would then act as a source (arrows pointing away), and a "low" would be a sink (arrows pointing towards), since air tends to move from high pressure areas to low pressure areas.
  \item Velocity field of a moving fluid. In this case, a velocity vector is associated to each point in the fluid.

  \item Streamlines, Streaklines and Pathlines are 3 types of lines that can be made from vector fields. They are :
\begin{enumerate}
   \item streaklines � as revealed in wind tunnels using smoke.
  \item streamlines (or fieldlines)� as a line depicting the instantaneous field at a given time.
  \item pathlines - showing the path that a given particle (of zero mass) would follow.

\end{enumerate}

\item Magnetic fields. The fieldlines can be revealed using small iron filings.

\item Maxwell's equations allow us to use a given set of initial conditions to deduce, for every point in Euclidean space, a magnitude and direction for the force experienced by a charged test particle at that point; the resulting vector field is the electromagnetic field.
\item A gravitational field generated by any massive object is also a vector field. For example, the gravitational field vectors for a spherically symmetric body would all point towards the sphere's center with the magnitude of the vectors reducing as radial distance from the body increases.
\end{itemize}

\subsection{Gradient field}

Vector fields can be constructed out of scalar fields using the gradient operator (denoted by the del: $\nabla$).

A vector field $V$ defined on a set $S$ is called a gradient field or a conservative field if there exists a real-valued function (a scalar field) $f$ on $S$ such that
$$
V = \nabla f = \bigg(\frac{\partial f}{\partial x_1}, \frac{\partial f}{\partial x_2}, \frac{\partial f}{\partial x_3}, \dots ,\frac{\partial f}{\partial x_n}\bigg).
$$
The associated flow is called the gradient flow, and is used in the method of gradient descent.

The path integral along any closed curve $\gamma (\gamma(0) = \gamma(1))$ in a conservative field is zero:

 $$\oint_\gamma \langle V(x), \mathrm{d}x \rangle = \oint_\gamma \langle \nabla f(x), \mathrm{d}x \rangle = f(\gamma(1)) - f(\gamma(0))$$
where the angular brackets and comma: $<,>$  denotes the inner product of two vectors (strictly speaking � the integrand $V(x)$ is a 1-form rather than a vector in the elementary sense).
\subsection{Central field}
A $C^\infty$-vector field over $R^n \backslash \{0\}$ is called a central field if

$$V(T(p)) = T(V(p)) \qquad (T \in \mathrm{O}(n, \mathbf{R}))$$

where $O(n, R)$ is the orthogonal group. We say central fields are invariant under orthogonal transformations around $0$.

The point $0$ is called the center of the field.

Since orthogonal transformations are actually rotations and reflections, the invariance conditions mean that vectors of a central field are always directed towards, or away from, $0$; this is an alternate (and simpler) definition. A central field is always a gradient field, since defining it on one semiaxis and integrating gives an antigradient.
\subsection{Operations on vector fields}

\begin{itemize}
  \item \textbf{Line integral} A common technique in physics is to integrate a vector field along a curve, i.e. to determine its line integral. Given a particle in a gravitational vector field, where each vector represents the force acting on the particle at a given point in space, the line integral is the work done on the particle when it travels along a certain path.

The line integral is constructed analogously to the Riemann integral and it exists if the curve is rectifiable (has finite length) and the vector field is continuous.

Given a vector field $V$ and a curve $\gamma$ parametrized by $[a, b]$ (where $a$ and $b$ are real) the line integral is defined as
$$\int_\gamma \langle V(x), dx>\rangle= \int_a^b \langle V(\gamma(t)), \gamma'(t)dt\rangle$$

  \item \textbf{Divergence} The divergence of a vector field on Euclidean space is a function (or scalar field). In three-dimensions, the divergence is defined by

$$\operatorname{div} \mathbf{F} = \nabla \cdot \mathbf{F} = \frac{\partial F_1}{\partial x} + \frac{\partial F_2}{\partial y}+\frac{\partial F_3}{\partial z}$$
with the obvious generalization to arbitrary dimensions. The divergence at a point represents the degree to which a small volume around the point is a source or a sink for the vector flow, a result which is made precise by the divergence theorem.

The divergence can also be defined on a Riemannian manifold, that is, a manifold with a Riemannian metric that measures the length of vectors
  \item \textbf{Curl} The curl is an operation which takes a vector field and produces another vector field. The curl is defined only in three-dimensions, but some properties of the curl can be captured in higher dimensions with the exterior derivative. In three-dimensions, it is defined by

$$\operatorname{curl}\,\mathbf{F} = \nabla \times \mathbf{F} = \left(\frac{\partial F_3}{\partial y}- \frac{\partial F_2}{\partial z}\right)\mathbf{e}_1 - \left(\frac{\partial F_3}{\partial x}- \frac{\partial F_1}{\partial z}\right)\mathbf{e}_2 + \left(\frac{\partial F_2}{\partial x}- \frac{\partial F_1}{\partial y}\right)\mathbf{e}_3$$
The curl measures the density of the angular momentum of the vector flow at a point, that is, the amount to which the flow circulates around a fixed axis. This intuitive description is made precise by Stokes' theorem.
  \item \textbf{Index of a vector field} The index of a vector field is a way of describing the behaviour of a vector field around an isolated zero (i.e. non-singular point) which can distinguish saddles from sources and sinks. Take a small sphere around the zero so that no other zeros are included. A map from this sphere to a unit sphere of dimensions $n-1$ can be constructed by dividing each vector by its length to form a unit length vector which can then be mapped to the unit sphere. The index of the vector field at the point is the degree of this map. The index of the vector field is the sum of the indices of each zero.

The index will be zero around any non singular point, it is $+1$ around sources and sinks and -1 around saddles. In two dimensions the index is equivalent to the winding number. For an ordinary sphere in three dimension space it can be shown that the index of any vector field on the sphere must be two, this leads to the hairy ball theorem which shows that every such vector field must have a zero. This theorem generalises to the Poincar\'{e}-Hopf theorem which relates the index to the Euler characteristic of the space.
\end{itemize}

\subsection{Flow curves}
Consider the flow of a fluid through a region of space. At any given time, any point of the fluid has a particular velocity associated with it; thus there is a vector field associated to any flow. The converse is also true: it is possible to associate a flow to a vector field having that vector field as its velocity.

Given a vector field $V$ defined on $S$, one defines curves $\gamma(t)$ on $S$ such that for each t in an interval $I$
$$
\gamma'(t) = V(\gamma(t))$$
By the Picard-Lindel\"{o}f theorem, if $V$ is Lipschitz continuous there is a unique $C^1$-curve $\gamma_x$ for each point $x$ in $S$ so that

$$\gamma_x(0) = x$$
$$\gamma'_x(t) = V(\gamma_x(t)) \qquad ( t \in (-\varepsilon, +\varepsilon) \subset \mathbf{R})$$
The curves $\gamma_x$ are called flow curves (or less commonly, flow lines) of the vector field $V$ and partition $S$ into equivalence classes. It is not always possible to extend the interval $(-e, +e)$ to the whole real number line. The flow may for example reach the edge of $S$ in a finite time. In two or three dimensions one can visualize the vector field as giving rise to a flow on $S$. If we drop a particle into this flow at a point $p$ it will move along the curve $\gamma_p$ in the flow depending on the initial point $p$. If $p$ is a stationary point of $V$ then the particle will remain at $p$.

Typical applications are streamline in fluid, geodesic flow, and one-parameter subgroups and the exponential map in Lie groups.

A vector field is \textbf{complete} if its flow curves exist for all time. In particular, compactly supported vector fields on a manifold are complete. If $X$ is a complete vector field on $M$, then the one-parameter group of diffeomorphisms generated by the flow along $X$ exists for all time.
\subsection{Difference between scalar and vector field}
The difference between a scalar and vector field is not that a scalar is just one number while a vector is several numbers. The difference is in how their coordinates respond to coordinate transformations. A scalar is a coordinate whereas a vector can be described by coordinates, but it is not the collection of its coordinates.

\begin{exam}
This example is about 2-dimensional Euclidean space $(R^2)$ where we examine Euclidean $(x, y)$ and polar $(r, \theta)$ coordinates (which are undefined at the origin). Thus $x = r cos\theta$  and $y = r sin\theta$  and also $r^2 = x^2 + y^2, cos\theta = x/(x^2 + y^2)^{1/2}$ and sin ? = y/(x2 + y2)1/2. Suppose we have a scalar field which is given by the constant function 1, and a vector field which attaches a vector in the r-direction with length 1 to each point. More precisely, they are given by the functions

$$s_{\mathrm{polar}}:(r, \theta) \mapsto 1, \quad v_{\mathrm{polar}}:(r, \theta) \mapsto (1, 0)$$
Let us convert these fields to Euclidean coordinates. The vector of length 1 in the $r$-direction has the $x$ coordinate $cos\theta$ and the $y$ coordinate $sin\theta$. Thus in Euclidean coordinates the same fields are described by the functions

$$s_{\mathrm{Euclidean}}:(x, y) \mapsto 1,$$
$$v_{\mathrm{Euclidean}}:(x, y) \mapsto (\cos \theta, \sin \theta) = \left(\frac{x}{\sqrt{{x^2 + y^2}}}, \frac{y}{\sqrt{x^2 + y^2}}\right)$$
We see that while the scalar field remains the same, the vector field now looks different. The same holds even in the 1-dimensional case, as illustrated by the next example.
\end{exam}
\begin{exam}
Consider the 1-dimensional Euclidean space $R$ with its standard Euclidean coordinate $x$. Suppose we have a scalar field and a vector field which are both given in the $x$ coordinate by the constant function 1,

$$s_{\mathrm{Euclidean}}:x \mapsto 1, \quad v_{\mathrm{Euclidean}}:x \mapsto 1$$
Thus, we have a scalar field which has the value 1 everywhere and a vector field which attaches a vector in the $x$-direction with magnitude 1 unit of $x$ to each point.

Now consider the coordinate $\xi := 2x$. If $x$ changes one unit then $\xi$ changes 2 units. But since we wish the integral of $v$ along a path to be independent of coordinate, this means $v^*dx=v'^*d\xi$. So from $x$ increase by 1 unit, $\xi$ increases by $1/2$ unit, so $v'$ must be 2. Thus this vector field has a magnitude of 2 in units of $\xi$. Therefore, in the $\xi$ coordinate the scalar field and the vector field are described by the functions

$$s_{\mathrm{unusual}}:\xi \mapsto 1, \quad v_{\mathrm{unusual}}:\xi \mapsto 2$$
which are different.
\end{exam}

\section{Noncommutative topology}
In mathematics, noncommutative topology is a term used for the relationship between topological and $C^*$-algebraic concepts. The term has its origins in the Gelfand-Naimark theorem, which implies the duality of the category of locally compact Hausdorff spaces and the category of commutative $C^*$-algebras. Noncommutative topology is related to analytic noncommutative geometry.

\subsection{Examples}
The premise behind noncommutative topology is that a noncommutative C*-algebra can be treated like the algebra of complex-valued continuous functions on a 'noncommutative space' which does not exist classically. Several topological properties can be formulated as properties for the C*-algebras without making reference to commutativity or the underlying space, and so have an immediate generalization. Among these are:

\begin{itemize}
  \item compactness (unital)
  \item $\sigma$-compactness ($\sigma$-unital)
  \item dimension (real or stable rank)
  \item connectedness (projectionless)
  \item Extremally disconnected spaces ($AW^*$-algebras)
\end{itemize}
In addition, certain types of continuous functions correspond with elements of C*-algebras. For example,
\begin{itemize}
  \item Real-valued continuous functions (self-adjoint elements)
  \item indicator functions of clopen sets (projections)
  \end{itemize}

There are certain examples of properties where multiple generalizations are possible and it is not clear which is preferable. For example, probability measures can correspond either to states or tracial states. Since all states are vacuously tracial states in the noncommutative case, so it is not clear whether the tracial condition is necessary to be a useful generalization.

\section{Involution}
In mathematics, an (anti-)involution, or an involutory function, is a function f that is its own inverse,
$f(f(x)) = x$ for all $x$ in the domain of $f$.

\subsection{General properties}
The identity map is a trivial example of an involution. Common examples in mathematics of nontrivial involutions include multiplication by -1 in arithmetic, the taking of reciprocals, complementation in set theory and complex conjugation. Other examples include circle inversion, rotation by a half-turn, and reciprocal ciphers such as the ROT13 transformation and the Beaufort polyalphabetic cipher.
The number of involutions, including the identity involution, on a set with $n = 0, 1, 2,$ � elements is given by a recurrence relation found by Heinrich August Rothe in 1800:

$$a_0 = a_1 = 1$$
$$a_n = a_n - 1 + (n - 1)a_n - 2, for n > 1$$
The first few terms of this sequence are 1, 1, 2, 4, 10, 26, 76, 232 (sequence A000085  in OEIS); these numbers are called the telephone numbers, and they also count the number of Young tableaux with a given number of cells. The composition $g \circ f$ of two involutions $f$ and $g$ is an involution if and only if they commute: $g \circ f=f \circ g$.

Every involution on an odd number of elements has at least one fixed point. More generally, for an involution on a finite set of elements, the number of elements and the number of fixed points have the same parity
\subsection{Involution throughout the fields of mathematics}
\begin{enumerate}
  \item \textbf{Pre-calculus}
Basic examples of involutions are the functions:

 $f(x) = -x  ,  f(x) = \frac {1}{x}$  and their combination  $f(x) = \frac {-1}{x}$
These are not the only pre-calculus involutions. Another in $ R^+$ is:

 $$f(x) = \ln\left(\frac {e^x+1}{e^x-1}\right) : x > 0 $$
  \item \textbf{Euclidean geometry}
A simple example of an involution of the three-dimensional Euclidean space is reflection against a plane. Performing a reflection twice brings a point back to its original coordinates.
Another is the so-called reflection through the origin; this is an abuse of language as it is not a reflection, though it is an involution.

These transformations are examples of affine involutions.
  \item \textbf{Projective geometry}
An involution is a projectivity of period 2, that is, a projectivity that interchanges pairs of points. Coxeter relates three theorems on involutions:

Any projectivity that interchanges two points is an involution.
\begin{itemize}
  \item The three pairs of opposite sides of a complete quadrangle meet any line (not through a vertex) in three pairs of an involution.

  \item If an involution has one fixed point, it has another, and consists of the correspondence between harmonic conjugates with respect to these two points. In this instance the involution is termed "hyperbolic", while if there are no fixed points it is "elliptic".
\end{itemize}
Another type of involution occurring in projective geometry is a polarity which is a correlation of period 2
  \item \textbf{Linear algebra} In linear algebra, an involution is a linear operator $T$ such that $T^2=I$. Except for in characteristic 2, such operators are diagonalizable with $1s$ and $-1s$ on the diagonal. If the operator is orthogonal (an orthogonal involution), it is orthonormally diagonalizable.
For example, suppose that a basis for a vector space $V$ is chosen, and that e1 and $e_2$ are basis elements. There exists a linear transformation $f$ which sends $e_1$ to $e_2$, and sends $e_2$ to $e_1$, and which is the identity on all other basis vectors. It can be checked that $f(f(x))=x$ for all $x \in V$. That is, $f$ is an involution of $V$.

This definition extends readily to modules. Given a module $M$ over a ring $R$, an $R$ endomorphism $f$ of $M$ is called an involution if $f^2$ is the identity homomorphism on $M$.
Involutions are related to idempotents; if 2 is invertible then they correspond in a one-to-one manner.
  \item \textbf{Quaternion algebra, groups, semigroups} In a quaternion algebra, an (anti-)involution is defined by the following axioms: if we consider a transformation $ x \mapsto f(x)$ then an involution is
      \begin{itemize}
        \item  $ f(f(x))=x$. An involution is its own inverse
        \item  An involution is linear:  $f(x_1+x_2)=f(x_1)+f(x_2)$  and  $f(\lambda x)=\lambda f(x)$
        \item  $f(x_1 x_2)=f(x_1) f(x_2)$ An anti-involution does not obey the last axiom but instead
        \item $f(x_1 x_2)=f(x_2) f(x_1)$
      \end{itemize}

This former law is sometimes called antidistributive. It also appears in groups as $(xy)^{-1} = y^{-1}x^{-1}$. Taken as an axiom, it leads to the notion of semigroup with involution, of which there are natural examples that are not groups, for example square matrix multiplication (i.e. the full linear monoid) with transpose as the involution.
  \item \textbf{Ring theory}
  In ring theory, the word involution is customarily taken to mean an antihomomorphism that is its own inverse function. Examples of involutions in common rings:
\begin{itemize}
  \item complex conjugation on the complex plane
  \item multiplication by j in the split-complex numbers
  \item taking the transpose in a matrix ring.
\end{itemize}
  \item \textbf{Group theory}
  In group theory, an element of a group is an involution if it has order 2; i.e. an involution is an element a such that $a \neq e$ and $a^2 = e$, where e is the identity element.
Originally, this definition agreed with the first definition above, since members of groups were always bijections from a set into itself; i.e., group was taken to mean permutation group. By the end of the 19th century, group was defined more broadly, and accordingly so was involution.
A permutation is an involution precisely if it can be written as a product of one or more non-overlapping transpositions.
The involutions of a group have a large impact on the group's structure. The study of involutions was instrumental in the classification of finite simple groups.
Coxeter groups are groups generated by involutions with the relations determined only by relations given for pairs of the generating involutions. Coxeter groups can be used, among other things, to describe the possible regular polyhedra and their generalizations to higher dimensions.
  \item \textbf{Mathematical logic} The operation of complement in Boolean algebras is an involution. Accordingly, negation in classical logic satisfies the law of double negation: $\neg\neg A$ is equivalent to A.

Generally in non-classical logics, negation that satisfies the law of double negation is called involutive. In algebraic semantics, such a negation is realized as an involution on the algebra of truth values. Examples of logics which have involutive negation are Kleene and Bochvar three-valued logics, Lukasiewicz many-valued logic, fuzzy logic IMTL, etc. Involutive negation is sometimes added as an additional connective to logics with non-involutive negation; this is usual, for example, in t-norm fuzzy logics.

The involutiveness of negation is an important characterization property for logics and the corresponding varieties of algebras. For instance, involutive negation characterizes Boolean algebras among Heyting algebras. Correspondingly, classical Boolean logic arises by adding the law of double negation to intuitionistic logic. The same relationship holds also between MV-algebras and BL-algebras (and so correspondingly between Lukasiewicz logic and fuzzy logic BL), IMTL and MTL, and other pairs of important varieties of algebras (resp. corresponding logics).
\end{enumerate}

\section{Hilbert space}
The mathematical concept of a Hilbert space, named after David Hilbert, generalizes the notion of Euclidean space. It extends the methods of vector algebra and calculus from the two-dimensional Euclidean plane and three-dimensional space to spaces with any finite or infinite number of dimensions. A Hilbert space is an abstract vector space possessing the structure of an inner product that allows length and angle to be measured. Furthermore, Hilbert spaces are complete: there are enough limits in the space to allow the techniques of calculus to be used.

Hilbert spaces arise naturally and frequently in mathematics and physics, typically as infinite-dimensional function spaces. The earliest Hilbert spaces were studied from this point of view in the first decade of the 20th century by David Hilbert, Erhard Schmidt, and Frigyes Riesz. They are indispensable tools in the theories of partial differential equations, quantum mechanics, Fourier analysis (which includes applications to signal processing and heat transfer)�and ergodic theory, which forms the mathematical underpinning of thermodynamics. John von Neumann coined the term Hilbert space for the abstract concept that underlies many of these diverse applications. The success of Hilbert space methods ushered in a very fruitful era for functional analysis. Apart from the classical Euclidean spaces, examples of Hilbert spaces include spaces of square-integrable functions, spaces of sequences, Sobolev spaces consisting of generalized functions, and Hardy spaces of holomorphic functions.

Geometric intuition plays an important role in many aspects of Hilbert space theory. Exact analogs of the Pythagorean theorem and parallelogram law hold in a Hilbert space. At a deeper level, perpendicular projection onto a subspace (the analog of "dropping the altitude" of a triangle) plays a significant role in optimization problems and other aspects of the theory. An element of a Hilbert space can be uniquely specified by its coordinates with respect to a set of coordinate axes (an orthonormal basis), in analogy with Cartesian coordinates in the plane. When that set of axes is countably infinite, this means that the Hilbert space can also usefully be thought of in terms of infinite sequences that are square-summable. Linear operators on a Hilbert space are likewise fairly concrete objects: in good cases, they are simply transformations that stretch the space by different factors in mutually perpendicular directions in a sense that is made precise by the study of their spectrum.

\subsection{Definition and illustration}

\begin{enumerate}
  \item \textbf{Motivating example: Euclidean space} One of the most familiar examples of a Hilbert space is the Euclidean space consisting of three-dimensional vectors, denoted by $R^3$, and equipped with the dot product. The dot product takes two vectors $x$ and $y$, and produces a real number $x�y$ If $x$ and $y$ are represented in Cartesian coordinates, then the dot product is defined by

$(x_1,x_2,x_3)\cdot (y_1,y_2,y_3) = x_1y_1+x_2y_2+x_3y_3$
The dot product satisfies the properties:
\begin{enumerate}
  \item It is symmetric in $x$ and $y$ : $x � y = y � x$

  \item It is linear in its first argument: $(ax_1 + bx_2) � y = ax_1 � y + bx_2 � y$ for any scalars $a, b$, and vectors $x_1, x_2,$ and $y$.

  \item It is positive definite: for all vectors $x$ then $ x � x = 0$ , with equality if and only if $x = 0$.

\end{enumerate}
An operation on pairs of vectors that, like the dot product, satisfies these three properties is known as a (real) inner product. A vector space equipped with such an inner product is known as a (real) inner product space. Every finite-dimensional inner product space is also a Hilbert space. The basic feature of the dot product that connects it with Euclidean geometry is that it is related to both the length (or norm) of a vector, denoted $||x||$, and to the angle $\theta$ between two vectors $x$ and $y$ by means of the formula

$$\mathbf{x}\cdot\mathbf{y} = \|\mathbf{x}\|\,\|\mathbf{y}\|\,\cos\theta$$

Multivariable calculus in Euclidean space relies on the ability to compute limits, and to have useful criteria for concluding that limits exist. A mathematical series
$$\sum_{n=0}^\infty \mathbf{x}_n$$
consisting of vectors in $R^3$ is absolutely convergent provided that the sum of the lengths converges as an ordinary series of real numbers:

$$\sum_{k=0}^\infty \|\mathbf{x}_k\| < \infty$$
Just as with a series of scalars, a series of vectors that converges absolutely also converges to some limit vector L in the Euclidean space, in the sense that

$$\left\|\mathbf{L}-\sum_{k=0}^N\mathbf{x}_k\right\|\to 0\quad\text{as }N\to\infty$$
This property expresses the completeness of Euclidean space: that a series that converges absolutely also converges in the ordinary sense.

Hilbert spaces are often taken over the complex numbers. The complex plane denoted by $C$ is equipped with a notion of magnitude, the complex modulus $|z|$ which is defined as the square root of the product of $z$ with its complex conjugate:

$$|z|^2 = z\overline{z}$$
If $z = x + iy$ is a decomposition of $z$ into its real and imaginary parts, then the modulus is the usual Euclidean two-dimensional length:

$$|z|=\sqrt{x^2+y^2}$$
The inner product of a pair of complex numbers $z$ and $w$ is the product of $z$ with the complex conjugate of $w$:

$$\langle z,w\rangle = z\overline{w}$$
This is complex-valued. The real part of $<z, w>$ gives the usual two-dimensional Euclidean dot product.

A second example is the space $C^2$ whose elements are pairs of complex numbers $z = (z_1, z_2)$. Then the inner product of $z$ with another such vector $w = (w_1, w_2)$ is given by

$$\langle z,w\rangle = z_1\overline{w}_1+z_2\overline{w}_2$$
The real part of $<z, w>$ is then the four-dimensional Euclidean dot product. This inner product is Hermitian symmetric, which means that the result of interchanging $z$ and $w$ is the complex conjugate:

$$\langle w,z\rangle = \overline{\langle z,w\rangle}$$
  \item \textbf{Definition}
  A Hilbert space $H$ is a real or complex inner product space that is also a complete metric space with respect to the distance function induced by the inner product.To say that $H$ is a complex inner product space means that $H$ is a complex vector space on which there is an inner product $ \langle x,y\rangle$ associating a complex number to each pair of elements $x,y$ of $H$ that satisfies the following properties:
  \begin{enumerate}
    \item The inner product of a pair of elements is equal to the complex conjugate of the inner product of the swapped elements:
$$\langle y,x\rangle = \overline{\langle x, y\rangle}$$
    \item The inner product is linear in its first argument. For all complex numbers $a$ and $b$,
$$\langle ax_1+bx_2, y\rangle = a\langle x_1, y\rangle + b\langle x_2, y\rangle$$
    \item The inner product of an element with itself is positive definite:
$$\langle x,x\rangle \ge 0$$

where the case of equality holds precisely when $x = 0$.
  \end{enumerate}

It follows from properties (a) and (b) that a complex inner product is antilinear in its second argument, meaning that

$$\langle x, ay_1+by_2\rangle = \bar{a}\langle x, y_1\rangle + \bar{b}\langle x, y_2\rangle$$
A real inner product space is defined in the same way, except that $H$ is a real vector space and the inner product takes real values. Such an inner product will be bilinear: that is, linear in each argument.

The norm is the real-valued function

$$\|x\| = \sqrt{\langle x,x \rangle}$$
and the distance d between two points $x,y \in H$ is defined in terms of the norm by

$$d(x,y)=\|x-y\| = \sqrt{\langle x-y,x-y \rangle}$$
That this function is a distance function means (a) that it is symmetric in $x$ and $y$, (b) that the distance between $x$ and itself is zero, and otherwise the distance between $x$ and $y$ must be positive, and (c) that the triangle inequality holds, meaning that the length of one leg of a triangle $xyz$ cannot exceed the sum of the lengths of the other two legs:

$$d(x,z) \le d(x,y) + d(y,z)$$
This last property is ultimately a consequence of the more fundamental Cauchy�Schwarz inequality, which asserts

$$|\langle x, y\rangle| \le \|x\|\,\|y\|$$
with equality if and only if x and y are linearly dependent.

Relative to a distance function defined in this way, any inner product space is a metric space, and sometimes is known as a pre-Hilbert space. Any pre-Hilbert space that is additionally also a complete space is a Hilbert space. Completeness is expressed using a form of the Cauchy criterion for sequences in $H$: a pre-Hilbert space $H$ is complete if every Cauchy sequence converges with respect to this norm to an element in the space. Completeness can be characterized by the following equivalent condition: if a series of vectors $\textstyle{\sum_{k=0}^\infty u_k}$ converges absolutely in the sense that

$$\sum_{k=0}^\infty\|u_k\| < \infty$$
then the series converges in $H$, in the sense that the partial sums converge to an element of $H$.

As a complete normed space, Hilbert spaces are by definition also Banach spaces. As such they are topological vector spaces, in which topological notions like the openness and closedness of subsets are well-defined. Of special importance is the notion of a closed linear subspace of a Hilbert space that, with the inner product induced by restriction, is also complete (being a closed set in a complete metric space) and therefore a Hilbert space in its own right

  \item \textbf{Second example: sequence spaces}
  The sequence space $\ell^2$ consists of all infinite sequences $z = (z_1,z_2,...)$ of complex numbers such that the series

$$\sum_{n=1}^\infty |z_n|^2$$
converges. The inner product on $\ell^2$ is defined by

$$\langle \mathbf{z},\mathbf{w}\rangle = \sum_{n=1}^\infty z_n\overline{w_n}$$
with the latter series converging as a consequence of the Cauchy�Schwarz inequality.

Completeness of the space holds provided that whenever a series of elements from $\ell^2$ converges absolutely (in norm), then it converges to an element of $\ell^2$. The proof is basic in mathematical analysis, and permits mathematical series of elements of the space to be manipulated with the same ease as series of complex numbers (or vectors in a finite-dimensional Euclidean space).
\end{enumerate}

\subsection{Examples}

\begin{enumerate}
  \item \textbf{Lebesgue spaces}
  Lebesgue spaces are function spaces associated to measure spaces $(X, M, \mu)$, where $X$ is a set, $M$ is a $\sigma$-algebra of subsets of $X$, and $\mu$ is a countably additive measure on $M$. Let $L^2(X, \mu)$ be the space of those complex-valued measurable functions on $X$ for which the Lebesgue integral of the square of the absolute value of the function is finite, i.e., for a function $f$ in $L^2(X,\mu)$,

 $\int_X |f|^2 d \mu  < \infty$
and where functions are identified if and only if they differ only on a set of measure zero.

The inner product of functions f and g in $L^2(X, \mu)$ is then defined as

$\langle f,g\rangle=\int_X f(t) \overline{g(t)} \ d \mu(t)$
For $f$ and $g$ in $L^2$, this integral exists because of the Cauchy�Schwarz inequality, and defines an inner product on the space. Equipped with this inner product, $L^2$ is in fact complete. The Lebesgue integral is essential to ensure completeness: on domains of real numbers, for instance, not enough functions are Riemann integrable.

The Lebesgue spaces appear in many natural settings. The spaces $L^2(R)$ and $L^2([0,1])$ of square-integrable functions with respect to the Lebesgue measure on the real line and unit interval, respectively, are natural domains on which to define the Fourier transform and Fourier series. In other situations, the measure may be something other than the ordinary Lebesgue measure on the real line. For instance, if $w$ is any positive measurable function, the space of all measurable functions $f$ on the interval $[0, 1]$ satisfying

$$\int_0^1 |f(t)|^2w(t)\,dt < \infty$$
is called the weighted $L^2$ space $L^2_{w([0,1])}$, and $w$ is called the weight function. The inner product is defined by

$$\langle f,g\rangle=\int_0^1 f(t) \overline{g(t)} w(t) dt$$
The weighted space $L^2_{w([0,1])}$ is identical with the Hilbert space $L^2([0,1],\mu)$ where the measure $\mu$ of a Lebesgue-measurable set $A$ is defined by

$$\mu(A) = \int_A w(t)dt$$
Weighted $L^2$ spaces like this are frequently used to study orthogonal polynomials, because different families of orthogonal polynomials are orthogonal with respect to different weighting functions.
  \item \textbf{Sobolev spaces} Sobolev spaces, denoted by $H^s$ or $W^{ s, 2}$, are Hilbert spaces. These are a special kind of function space in which differentiation may be performed, but that (unlike other Banach spaces such as the H\"{o}lder spaces) support the structure of an inner product. Because differentiation is permitted, Sobolev spaces are a convenient setting for the theory of partial differential equations. They also form the basis of the theory of direct methods in the calculus of variations.
For $s$ a non-negative integer and $\Omega \subset R^n$, the Sobolev space $H^s(\Omega)$ contains $L^2$ functions whose weak derivatives of order up to s are also $L^2$. The inner product in $H^s(\Omega)$ is

$$\langle f,g\rangle = \int_\Omega f(x)\bar{g}(x)dx + \int_\Omega D f(x)\cdot D\bar{g}(x)dx + \cdots + \int_\Omega D^s f(x)\cdot D^s \bar{g}(x) dx$$
where the dot indicates the dot product in the Euclidean space of partial derivatives of each order. Sobolev spaces can also be defined when $s$ is not an integer.

Sobolev spaces are also studied from the point of view of spectral theory, relying more specifically on the Hilbert space structure. If $\Omega$ is a suitable domain, then one can define the Sobolev space $H^s(\Omega)$ as the space of Bessel potentials roughly,

$$H^s(\Omega) = \{ (1-\Delta)^{-s/2}f | f\in L^2(\Omega)\}$$
Here $\Delta$ is the Laplacian and $(1-\Delta)^{-s/2}$ is understood in terms of the spectral mapping theorem. Apart from providing a workable definition of Sobolev spaces for non-integer $s$, this definition also has particularly desirable properties under the Fourier transform that make it ideal for the study of pseudodifferential operators. Using these methods on a compact Riemannian manifold, one can obtain for instance the Hodge decomposition, which is the basis of Hodge theory.

  \item \textbf{Spaces of holomorphic functions}
  \begin{itemize}
    \item \textbf{Hardy spaces}
    The Hardy spaces are function spaces, arising in complex analysis and harmonic analysis, whose elements are certain holomorphic functions in a complex domain. Let $U$ denote the unit disc in the complex plane. Then the Hardy space $H^2(U)$ is defined as the space of holomorphic functions $f$ on $U$ such that the means

$$M_r(f) = \frac{1}{2\pi}\int_0^{2\pi}|f(re^{i\theta})|^2\,d\theta$$
remain bounded for $r < 1$. The norm on this Hardy space is defined by

$$\|f\|_2 = \lim_{r\to 1} \sqrt{M_r(f)}$$
Hardy spaces in the disc are related to Fourier series. A function $f$ is in $H^2(U)$ if and only if

$$f(z) = \sum_{n=0}^\infty a_nz^n$$
where

$$\sum_{n=0}^\infty\,|a_n|^2 < \infty$$
Thus $H^2(U)$ consists of those functions that are $L^2$ on the circle, and whose negative frequency Fourier coefficients vanish.
    \item \textbf{Bergman spaces}
    The Bergman spaces are another family of Hilbert spaces of holomorphic functions. Let $D$ be a bounded open set in the complex plane (or a higher-dimensional complex space) and let $L^{2,h}(D)$ be the space of holomorphic functions $f$ in $D$ that are also in $L^2(D)$ in the sense that

$$\|f\|^2 = \int_D |f(z)|^2\,d\mu(z) < \infty$$
where the integral is taken with respect to the Lebesgue measure in $D$. Clearly $L^{2,h}(D)$ is a subspace of $L^2(D)$; in fact, it is a closed subspace, and so a Hilbert space in its own right. This is a consequence of the estimate, valid on compact subsets $K$ of $D$, that

$$\sup_{z\in K} |f(z)| \le C_K \|f\|_2$$
which in turn follows from Cauchy's integral formula. Thus convergence of a sequence of holomorphic functions in $L^2(D)$ implies also compact convergence, and so the limit function is also holomorphic. Another consequence of this inequality is that the linear functional that evaluates a function f at a point of $D$ is actually continuous on $L^{2,h}(D)$. The Riesz representation theorem implies that the evaluation functional can be represented as an element of $L^{2,h}(D)$. Thus, for every $z \in D$, there is a function $\eta_z \in L^{2,h}(D)$ such that

$$f(z) = \int_D f(\zeta)\overline{\eta_z(\zeta)}\,d\mu(\zeta)$$
for all $f \in L^{2,h}(D)$. The integrand

$$K(\zeta,z) = \overline{\eta_z(\zeta)}$$
is known as the Bergman kernel of $D$. This integral kernel satisfies a reproducing property

$$f(z) = \int_D f(\zeta)K(\zeta,z)\,d\mu(\zeta)$$
A Bergman space is an example of a reproducing kernel Hilbert space, which is a Hilbert space of functions along with a kernel $K(\zeta,z)$ that verifies a reproducing property analogous to this one. The Hardy space $H^2(D)$ also admits a reproducing kernel, known as the Szego kernel. Reproducing kernels are common in other areas of mathematics as well. For instance, in harmonic analysis the Poisson kernel is a reproducing kernel for the Hilbert space of square-integrable harmonic functions in the unit ball. That the latter is a Hilbert space at all is a consequence of the mean value theorem for harmonic functions.
  \end{itemize}
\end{enumerate}

\subsection{Applications}
Many of the applications of Hilbert spaces exploit the fact that Hilbert spaces support generalizations of simple geometric concepts like projection and change of basis from their usual finite dimensional setting. In particular, the spectral theory of continuous self-adjoint linear operators on a Hilbert space generalizes the usual spectral decomposition of a matrix, and this often plays a major role in applications of the theory to other areas of mathematics and physics.
\begin{enumerate}
  \item \textbf{Sturm�Liouville theory} In the theory of ordinary differential equations, spectral methods on a suitable Hilbert space are used to study the behavior of eigenvalues and eigenfunctions of differential equations. For example, the Sturm�Liouville problem arises in the study of the harmonics of waves in a violin string or a drum, and is a central problem in ordinary differential equations. The problem is a differential equation of the form

 $$-\frac{d}{dx}\left[p(x)\frac{dy}{ dx}\right]+q(x)y=\lambda w(x)y$$
for an unknown function $y$ on an interval $[a,b]$, satisfying general homogeneous Robin boundary conditions
$$
\begin{cases}
\alpha y(a)+\alpha' y'(a)=0\\
\beta y(b) + \beta' y'(b)=0.
\end{cases}
$$
The functions $p, q$, and $w$ are given in advance, and the problem is to find the function $y$ and constants $\lambda$ for which the equation has a solution. The problem only has solutions for certain values of $\lambda$, called eigenvalues of the system, and this is a consequence of the spectral theorem for compact operators applied to the integral operator defined by the Green's function for the system. Furthermore, another consequence of this general result is that the eigenvalues $\lambda$ of the system can be arranged in an increasing sequence tending to infinity
  \item \textbf{Partial differential equations} Hilbert spaces form a basic tool in the study of partial differential equations. For many classes of partial differential equations, such as linear elliptic equations, it is possible to consider a generalized solution (known as a weak solution) by enlarging the class of functions. Many weak formulations involve the class of Sobolev functions, which is a Hilbert space. A suitable weak formulation reduces to a geometrical problem the analytic problem of finding a solution or, often what is more important, showing that a solution exists and is unique for given boundary data. For linear elliptic equations, one geometrical result that ensures unique solvability for a large class of problems is the Lax�Milgram theorem. This strategy forms the rudiment of the Galerkin method (a finite element method) for numerical solution of partial differential equations A typical example is the Poisson equation $-\Delta \mu = g$ with Dirichlet boundary conditions in a bounded domain $\Omega \in R^2$. The weak formulation consists of finding a function $u$ such that, for all continuously differentiable functions $v \in \Omega$ vanishing on the boundary:

$$\int_\Omega \nabla u\cdot\nabla v = \int_\Omega gv$$
This can be recast in terms of the Hilbert space $H^1_0(\Omega)$ consisting of functions $u$ such that u, along with its weak partial derivatives, are square integrable on $\Omega$, and vanish on the boundary. The question then reduces to finding $u$ in this space such that for all $v$ in this space

$$a(u,v) = b(v)$$
where $a$ is a continuous bilinear form, and $b$ is a continuous linear functional, given respectively by

$$a(u,v) = \int_\Omega \nabla u\cdot\nabla v,\quad b(v)= \int_\Omega gv$$
Since the Poisson equation is elliptic, it follows from Poincar\'{e}'s inequality that the bilinear form a is coercive. The Lax�Milgram theorem then ensures the existence and uniqueness of solutions of this equation.

Hilbert spaces allow for many elliptic partial differential equations to be formulated in a similar way, and the Lax�Milgram theorem is then a basic tool in their analysis. With suitable modifications, similar techniques can be applied to parabolic partial differential equations and certain hyperbolic partial differential equations

  \item \textbf{Ergodic theory} The field of ergodic theory is the study of the long-term behavior of chaotic dynamical systems. The protypical case of a field that ergodic theory applies to is thermodynamics, in which�though the microscopic state of a system is extremely complicated (it is impossible to understand the ensemble of individual collisions between particles of matter)�the average behavior over sufficiently long time intervals is tractable. The laws of thermodynamics are assertions about such average behavior. In particular, one formulation of the zeroth law of thermodynamics asserts that over sufficiently long timescales, the only functionally independent measurement that one can make of a thermodynamic system in equilibrium is its total energy, in the form of temperature.

An ergodic dynamical system is one for which, apart from the energy�measured by the Hamiltonian�there are no other functionally independent conserved quantities on the phase space. More explicitly, suppose that the energy $E$ is fixed, and let $\Omega_E$ be the subset of the phase space consisting of all states of energy $E$ (an energy surface), and let $T_t$ denote the evolution operator on the phase space. The dynamical system is ergodic if there are no continuous non-constant functions on $\Omega_E$ such that

$$f(T_tw) = f(w)$$
for all w on $\Omega_E$ and all time $t$. Liouville's theorem implies that there exists a measure $\mu$ on the energy surface that is invariant under the time translation. As a result, time translation is a unitary transformation of the Hilbert space $L^2(\Omega_E,\mu)$ consisting of square-integrable functions on the energy surface $\Omega_E$ with respect to the inner product

$$\langle f,g\rangle_{L^2(\Omega_E,\mu)} = \int_E f\bar{g}\,d\mu$$
The von Neumann mean ergodic theorem states the following:

If $U_t$ is a (strongly continuous) one-parameter semigroup of unitary operators on a Hilbert space $H$, and $P$ is the orthogonal projection onto the space of common fixed points of $U_t$, $\{x \in H | U_tx = x  \forall t > 0\}$, then
$$Px = \lim_{T\to\infty}\frac{1}{T}\int_0^TU_tx dt$$
For an ergodic system, the fixed set of the time evolution consists only of the constant functions, so the ergodic theorem implies the following for any function $f \in L^2(\Omega_E,\mu)$

$$\underset{T\to\infty}{L^2\!-\!\lim} \frac{1}{T}\int_0^T f(T_tw)\,dt = \int_{\Omega_E} f(y)\,d\mu(y)$$
That is, the long time average of an observable $f$ is equal to its expectation value over an energy surface
  \item \textbf{Fourier analysis} One of the basic goals of Fourier analysis is to decompose a function into a (possibly infinite) linear combination of given basis functions: the associated Fourier series. The classical Fourier series associated to a function $f$ defined on the interval $[0, 1]$ is a series of the form

$$\sum_{n=-\infty}^\infty a_n e^{2\pi in\theta}$$ where $$a_n = \int_0^1f(\theta)e^{-2\pi in\theta} d\theta$$
The example of adding up the first few terms in a Fourier series for a sawtooth function is shown in the figure. The basis functions are sine waves with wavelengths $\lambda/n$ ($n$=integer) shorter than the wavelength $\lambda$ of the sawtooth itself (except for $n=1$, the fundamental wave). All basis functions have nodes at the nodes of the sawtooth, but all but the fundamental have additional nodes. The oscillation of the summed terms about the sawtooth is called the Gibbs phenomenon.

A significant problem in classical Fourier series asks in what sense the Fourier series converges, if at all, to the function $f$. Hilbert space methods provide one possible answer to this question. The functions $e_n(\theta) = e^{2\pi in\theta}$ form an orthogonal basis of the Hilbert space $L^2([0,1])$. Consequently, any square-integrable function can be expressed as a series

$$f(\theta) = \sum_n a_n e_n(\theta),\quad a_n = \langle f,e_n\rangle$$
and, moreover, this series converges in the Hilbert space sense (that is, in the $L^2$ mean).

The problem can also be studied from the abstract point of view: every Hilbert space has an orthonormal basis, and every element of the Hilbert space can be written in a unique way as a sum of multiples of these basis elements. The coefficients appearing on these basis elements are sometimes known abstractly as the Fourier coefficients of the element of the space. The abstraction is especially useful when it is more natural to use different basis functions for a space such as $L^2([0,1])$. In many circumstances, it is desirable not to decompose a function into trigonometric functions, but rather into orthogonal polynomials or wavelets for instance, and in higher dimensions into spherical harmonics.
For instance, if en are any orthonormal basis functions of $L^2[0,1]$, then a given function in $L^2[0,1]$ can be approximated as a finite linear combination

$$f(x) \approx f_n (x) = a_1 e_1 (x) + a_2 e_2(x) + \cdots + a_n e_n (x)$$
The coefficients $\{a_j\}$ are selected to make the magnitude of the difference $||f - f_n||^2$ as small as possible. Geometrically, the best approximation is the orthogonal projection of $f$ onto the subspace consisting of all linear combinations of the $\{e_j\}$, and can be calculated by

$$a_j = \int_0^1 \overline{e_j(x)}f (x) dx$$
That this formula minimizes the difference $||f - f_n||^2$ is a consequence of Bessel's inequality and Parseval's formula.

In various applications to physical problems, a function can be decomposed into physically meaningful eigenfunctions of a differential operator (typically the Laplace operator): this forms the foundation for the spectral study of functions, in reference to the spectrum of the differential operator. A concrete physical application involves the problem of hearing the shape of a drum: given the fundamental modes of vibration that a drumhead is capable of producing, can one infer the shape of the drum itself? The mathematical formulation of this question involves the Dirichlet eigenvalues of the Laplace equation in the plane, that represent the fundamental modes of vibration in direct analogy with the integers that represent the fundamental modes of vibration of the violin string.

Spectral theory also underlies certain aspects of the Fourier transform of a function. Whereas Fourier analysis decomposes a function defined on a compact set into the discrete spectrum of the Laplacian (which corresponds to the vibrations of a violin string or drum), the Fourier transform of a function is the decomposition of a function defined on all of Euclidean space into its components in the continuous spectrum of the Laplacian. The Fourier transformation is also geometrical, in a sense made precise by the Plancherel theorem, that asserts that it is an isometry of one Hilbert space (the "time domain") with another (the "frequency domain"). This isometry property of the Fourier transformation is a recurring theme in abstract harmonic analysis, as evidenced for instance by the Plancherel theorem for spherical functions occurring in noncommutative harmonic analysis.

  \item \textbf{Quantum mechanics} In the mathematically rigorous formulation of quantum mechanics, developed by John von Neumann, the possible states (more precisely, the pure states) of a quantum mechanical system are represented by unit vectors (called state vectors) residing in a complex separable Hilbert space, known as the state space, well defined up to a complex number of norm 1 (the phase factor). In other words, the possible states are points in the projectivization of a Hilbert space, usually called the complex projective space. The exact nature of this Hilbert space is dependent on the system; for example, the position and momentum states for a single non-relativistic spin zero particle is the space of all square-integrable functions, while the states for the spin of a single proton are unit elements of the two-dimensional complex Hilbert space of spinors. Each observable is represented by a self-adjoint linear operator acting on the state space. Each eigenstate of an observable corresponds to an eigenvector of the operator, and the associated eigenvalue corresponds to the value of the observable in that eigenstate.

The inner product between two state vectors is a complex number known as a probability amplitude. During an ideal measurement of a quantum mechanical system, the probability that a system collapses from a given initial state to a particular eigenstate is given by the square of the absolute value of the probability amplitudes between the initial and final states. The possible results of a measurement are the eigenvalues of the operator�which explains the choice of self-adjoint operators, for all the eigenvalues must be real. The probability distribution of an observable in a given state can be found by computing the spectral decomposition of the corresponding operator.

For a general system, states are typically not pure, but instead are represented as statistical mixtures of pure states, or mixed states, given by density matrices: self-adjoint operators of trace one on a Hilbert space. Moreover, for general quantum mechanical systems, the effects of a single measurement can influence other parts of a system in a manner that is described instead by a positive operator valued measure. Thus the structure both of the states and observables in the general theory is considerably more complicated than the idealization for pure states.
\end{enumerate}

\subsection{Properties}
\begin{enumerate}
  \item \textbf{Pythagorean identity} Two vectors $u$ and $v$ in a Hilbert space $H$ are orthogonal when $\langle u, v\rangle = 0$. The notation for this is $u \bot v$. More generally, when $S$ is a subset in $H$, the notation $u \bot S$ means that $u$ is orthogonal to every element from $S$.
When $u$ and $v$ are orthogonal, one has

$$\|u + v\|^2 = \langle u + v, u + v \rangle = \langle u, u \rangle + 2 \, \mathrm{Re} \langle u, v \rangle + \langle v, v \rangle= \|u\|^2 + \|v\|^2$$
By induction on $n$, this is extended to any family $u_1,...,u_n$ of n orthogonal vectors,

$$\|u_1 + \cdots + u_n\|^2 = \|u_1\|^2 + \cdots + \|u_n\|^2$$
Whereas the Pythagorean identity as stated is valid in any inner product space, completeness is required for the extension of the Pythagorean identity to series. A series $\sum u_k$ of orthogonal vectors converges in $H$ if and only if the series of squares of norms converges, and

$$\bigl\|\sum_{k=0}^\infty u_k \bigr\|^2 = \sum_{k=0}^\infty \|u_k\|^2$$
Furthermore, the sum of a series of orthogonal vectors is independent of the order in which it is taken.
  \item \textbf{Parallelogram identity and polarization} By definition, every Hilbert space is also a Banach space. Furthermore, in every Hilbert space the following parallelogram identity holds:

$$\|u+v\|^2+\|u-v\|^2=2(\|u\|^2+\|v\|^2)$$
Conversely, every Banach space in which the parallelogram identity holds is a Hilbert space, and the inner product is uniquely determined by the norm by the polarization identity. For real Hilbert spaces, the polarization identity is

$$\langle u,v\rangle = \frac{1}{4}\left(\|u+v\|^2-\|u-v\|^2\right)$$
For complex Hilbert spaces, it is

$$\langle u,v\rangle = \frac{1}{4}\left(\|u+v\|^2-\|u-v\|^2+i\|u+iv\|^2-i\|u-iv\|^2\right)$$
The parallelogram law implies that any Hilbert space is a uniformly convex Banach space
  \item \textbf{Best approximation} This subsection employs the Hilbert projection theorem. If $C$ is a non-empty closed convex subset of a Hilbert space $H$ and $x$ a point in $H$, there exists a unique point $y \in C$ that minimizes the distance between $x$ and points in $C$,

$$ y \in C, \ \ \ \|x - y\| = \mathrm{dist}(x, C) = \min \{ \|x - z\| : z \in C \}$$
This is equivalent to saying that there is a point with minimal norm in the translated convex set $D = C - x$. The proof consists in showing that every minimizing sequence $(d_n) \subset D$ is Cauchy (using the parallelogram identity) hence converges (using completeness) to a point in $D$ that has minimal norm. More generally, this holds in any uniformly convex Banach space.

When this result is applied to a closed subspace $F$ of $H$, it can be shown that the point $y \in F$ closest to $x$ is characterized by

 $$y \in F, \ \ x - y \perp F$$
This point $y$ is the orthogonal projection of $x$ onto $F$, and the mapping $P_F : x \rightarrow y$ is linear (see Orthogonal complements and projections). This result is especially significant in applied mathematics, especially numerical analysis, where it forms the basis of least squares methods.

In particular, when $F$ is not equal to $H$, one can find a non-zero vector $v$ orthogonal to $F$ (select $x$ not in $F$ and $v = x - y$). A very useful criterion is obtained by applying this observation to the closed subspace $F$ generated by a subset $S$ of $H$.

A subset $S$ of $H$ spans a dense vector subspace if (and only if) the vector 0 is the sole vector $v \in H$ orthogonal to $S$.
  \item \textbf{Duality} The dual space H* is the space of all continuous linear functions from the space $H$ into the base field. It carries a natural norm, defined by

$$\|\varphi\| = \sup_{\|x\|=1, x\in H} |\varphi(x)|$$
This norm satisfies the parallelogram law, and so the dual space is also an inner product space. The dual space is also complete, and so it is a Hilbert space in its own right.

The Riesz representation theorem affords a convenient description of the dual. To every element $u$ of $H$, there is a unique element $\varphi_u$ of $H^*$, defined by

$$\varphi_u(x) = \langle x,u\rangle$$
The mapping $u\mapsto \varphi_u$ is an antilinear mapping from $H$ to $H^*$. The Riesz representation theorem states that this mapping is an antilinear isomorphism. Thus to every element $\varphi$ of the dual $H^*$ there exists one and only one $u_\varphi$ in H such that

$$\langle x, u_\varphi\rangle = \varphi(x)$$
for all $x \in H$. The inner product on the dual space $H^*$ satisfies

 $$\langle \varphi, \psi \rangle = \langle u_\psi, u_\varphi \rangle$$
The reversal of order on the right-hand side restores linearity in $\varphi$ from the antilinearity of $u_\varphi$. In the real case, the antilinear isomorphism from $H$ to its dual is actually an isomorphism, and so real Hilbert spaces are naturally isomorphic to their own duals.

The representing vector $u_\varphi$ is obtained in the following way. When $\varphi \neq 0$, the $kernel F = Ker(\varphi)$ is a closed vector subspace of $H$, not equal to $H$, hence there exists a non-zero vector $v$ orthogonal to $F$. The vector $u$ is a suitable scalar multiple $\lambda v$ of $v$. The requirement that $\varphi(v) =\langle v, u \rangle$ yields

 $u = \langle v, v \rangle^{-1} \, \overline{\varphi (v)}  v$
This correspondence $\varphi \leftrightarrow u$ is exploited by the bra�ket notation popular in physics. It is common in physics to assume that the inner product, denoted by $\langle x|y \rangle$, is linear on the right,

$$\langle x| y \rangle = \langle y, x \rangle$$
The result $\langle x|y \rangle$ can be seen as the action of the linear functional $\langle x| (the bra) on the vector  |y \rangle$ (the ket).

The Riesz representation theorem relies fundamentally not just on the presence of an inner product, but also on the completeness of the space. In fact, the theorem implies that the topological dual of any inner product space can be identified with its completion. An immediate consequence of the Riesz representation theorem is also that a Hilbert space $H$ is reflexive, meaning that the natural map from $H$ into its double dual space is an isomorphism.
  \item \textbf{Weakly convergent sequences} In a Hilbert space $H$, a sequence $\{x_n\}$ is weakly convergent to a vector $x \in H$ when

$$\lim_n \langle x_n, v \rangle = \langle x, v \rangle$$
for every $v \in H$

For example, any orthonormal sequence $\{f_n\}$ converges weakly to 0, as a consequence of Bessel's inequality. Every weakly convergent sequence $\{x_n\}$ is bounded, by the uniform boundedness principle.

Conversely, every bounded sequence in a Hilbert space admits weakly convergent subsequences (Alaoglu's theorem). This fact may be used to prove minimization results for continuous convex functionals, in the same way that the Bolzano�Weierstrass theorem is used for continuous functions on Rd. Among several variants, one simple statement is as follows:
If $f: H \rightarrow R$ is a convex continuous function such that $f(x)$ tends to $+\infty$ when $||x||$ tends to $\infty$, then $f$ admits a minimum at some point $x_0 \in H$.
This fact (and its various generalizations) are fundamental for direct methods in the calculus of variations. Minimization results for convex functionals are also a direct consequence of the slightly more abstract fact that closed bounded convex subsets in a Hilbert space $H$ are weakly compact, since $H$ is reflexive. The existence of weakly convergent subsequences is a special case of the Eberlein��mulian theorem.
\end{enumerate}

\subsection{Operators on Hilbert spaces}
\begin{enumerate}
  \item \textbf{Bounded operators} The continuous linear operators $A : H_1 \rightarrow H_2$ from a Hilbert space $H_1$ to a second Hilbert space $H_2$ are bounded in the sense that they map bounded sets to bounded sets. Conversely, if an operator is bounded, then it is continuous. The space of such bounded linear operators has a norm, the operator norm given by

$$\lVert A \rVert = \sup \left\{\,\lVert Ax \rVert : \lVert x \rVert \leq 1\,\right\}$$
The sum and the composite of two bounded linear operators is again bounded and linear. For $y$ in $H_2$, the map that sends $x \in H_1$ to $\langle Ax, y\rangle$ is linear and continuous, and according to the Riesz representation theorem can therefore be represented in the form

$$\langle x, A^* y \rangle = \langle Ax, y \rangle$$
for some vector $A^*y$ in $H_1$. This defines another bounded linear operator $A^*: H_2 \rightarrow H_1$, the adjoint of $A$. One can see that $A^{**} = A$.

The set $B(H)$ of all bounded linear operators on $H$, together with the addition and composition operations, the norm and the adjoint operation, is a $C^*$-algebra, which is a type of operator algebra.

An element $A$ of $B(H)$ is called self-adjoint or Hermitian if $A^*= A$. If $A$ is Hermitian and $\langle Ax, x \rangle = 0$ for every $x$, then $A$ is called non-negative, written $A = 0$; if equality holds only when $x = 0$, then $A$ is called positive. The set of self adjoint operators admits a partial order, in which $A \geq B$ if $A - B \geq 0$. If $A$ has the form $B^*B$ for some $B$, then $A$ is non-negative; if $B$ is invertible, then $A$ is positive. A converse is also true in the sense that, for a non-negative operator $A$, there exists a unique non-negative square root $B$ such that

$$A = B^2=B^*B$$,
In a sense made precise by the spectral theorem, self-adjoint operators can usefully be thought of as operators that are "real". An element $A$ of $B(H)$ is called normal if $A^*A = AA^*$. Normal operators decompose into the sum of a self-adjoint operators and an imaginary multiple of a self adjoint operator

$$A = \frac{A+A^*}{2} + i\frac{A-A^*}{2i}$$
that commute with each other. Normal operators can also usefully be thought of in terms of their real and imaginary parts.

An element $U$ of $B(H)$ is called unitary if $U$ is invertible and its inverse is given by $U^*$. This can also be expressed by requiring that $U$ be onto and $\langle Ux, Uy \rangle = \langle x, y \rangle$ for all $x$ and $y$ in H. The unitary operators form a group under composition, which is the isometry group of $H$.

An element of $B(H)$ is compact if it sends bounded sets to relatively compact sets. Equivalently, a bounded operator $T$ is compact if, for any bounded sequence $\{x_k\}$, the sequence $\{Tx_k\}$ has a convergent subsequence. Many integral operators are compact, and in fact define a special class of operators known as Hilbert-Schmidt operators that are especially important in the study of integral equations. Fredholm operators differ from a compact operator by a multiple of the identity, and are equivalently characterized as operators with a finite dimensional kernel and cokernel. The index of a Fredholm operator $T$ is defined by

$$\operatorname{index}\, T = \dim\ker T - \dim\operatorname{coker} T$$
The index is homotopy invariant, and plays a deep role in differential geometry via the Atiyah-Singer index theorem.
  \item \textbf{Unbounded operators} Unbounded operators are also tractable in Hilbert spaces, and have important applications to quantum mechanics. An unbounded operator $T$ on a Hilbert space $H$ is defined as a linear operator whose domain $D(T)$ is a linear subspace of $H$. Often the domain $D(T)$ is a dense subspace of $H$, in which case $T$ is known as a densely defined operator.

The adjoint of a densely defined unbounded operator is defined in essentially the same manner as for bounded operators. Self-adjoint unbounded operators play the role of the observables in the mathematical formulation of quantum mechanics. Examples of self-adjoint unbounded operators on the Hilbert space $L^2(R)$ are:
\begin{itemize}
  \item A suitable extension of the differential operator
 $$(A f)(x) = -i \frac{d}{dx} f(x),$$
where $i$ is the imaginary unit and $f$ is a differentiable function of compact support.
  \item The multiplication-by-$x$ operator:
$$ (B f) (x) = x f(x)$$
\end{itemize}
These correspond to the momentum and position observables, respectively. Note that neither $A$ nor $B$ is defined on all of $H$, since in the case of $A$ the derivative need not exist, and in the case of $B$ the product function need not be square integrable. In both cases, the set of possible arguments form dense subspaces of $L^2(R)$.

\end{enumerate}

\subsection{Constructions}
\begin{enumerate}
  \item \textbf{Direct sums} Two Hilbert spaces $H_1$ and $H_2$ can be combined into another Hilbert space, called the (orthogonal) direct sum, and denoted

$$H_1\oplus H_2$$
consisting of the set of all ordered pairs $(x_1, x_2)$ where $x_i \in H_i, i = 1,2,$ and inner product defined by

$$\langle (x_1,x_2), (y_1,y_2)\rangle_{H_1\oplus H_2} = \langle x_1,y_1\rangle_{H_1} + \langle x_2,y_2\rangle_{H_2}$$
More generally, if $H_i$ is a family of Hilbert spaces indexed by $i \in I$, then the direct sum of the $H_i$, denoted

$$\bigoplus_{i\in I}H_i$$
consists of the set of all indexed families

$$x=(x_i\in H_i|i\in I) \in \prod_{i\in I}H_i$$
in the Cartesian product of the $H_i$ such that

$$\sum_{i\in I} \|x_i\|^2 < \infty$$
The inner product is defined by

$$\langle x, y\rangle = \sum_{i\in I} \langle x_i, y_i\rangle_{H_i}$$
Each of the $H_i$ is included as a closed subspace in the direct sum of all of the $H_i$. Moreover, the $H_i$ are pairwise orthogonal. Conversely, if there is a system of closed subspaces, $V_i, i \in I$, in a Hilbert space $H$, that are pairwise orthogonal and whose union is dense in $H$, then $H$ is canonically isomorphic to the direct sum of $V_i$. In this case, $H$ is called the internal direct sum of the $V_i$. A direct sum (internal or external) is also equipped with a family of orthogonal projections $E_i$ onto the $i^{th}$ direct summand $H_i$. These projections are bounded, self-adjoint, idempotent operators that satisfy the orthogonality condition

$$E_iE_j = 0,\quad i\not= j$$
The spectral theorem for compact self-adjoint operators on a Hilbert space $H$ states that $H$ splits into an orthogonal direct sum of the eigenspaces of an operator, and also gives an explicit decomposition of the operator as a sum of projections onto the eigenspaces. The direct sum of Hilbert spaces also appears in quantum mechanics as the Fock space of a system containing a variable number of particles, where each Hilbert space in the direct sum corresponds to an additional degree of freedom for the quantum mechanical system. In representation theory, the Peter-Weyl theorem guarantees that any unitary representation of a compact group on a Hilbert space splits as the direct sum of finite-dimensional representations.
  \item \textbf{Tensor products} If $H_1$ and $H_2$, then one defines an inner product on the (ordinary) tensor product as follows. On simple tensors, let

 $$\langle x_1 \otimes x_2, \, y_1 \otimes y_2 \rangle = \langle x_1, y_1 \rangle \, \langle x_2, y_2 \rangle$$
This formula then extends by sesquilinearity to an inner product on $H_1 \otimes H_2$. The Hilbertian tensor product of $H_1$ and $H_2$, sometimes denoted by $H_1\widehat{\otimes}H_2$, is the Hilbert space obtained by completing $H_1 \otimes H_2$ for the metric associated to this inner product.
An example is provided by the Hilbert space $L^2([0, 1])$. The Hilbertian tensor product of two copies of $L^2([0, 1])$ is isometrically and linearly isomorphic to the space $L^2([0, 1]^2)$ of square-integrable functions on the square $[0, 1]^2$. This isomorphism sends a simple tensor $f_1 \otimes f_2$ to the function

 $$(s, t) \mapsto f_1(s) \, f_2(t)$$
on the square.

This example is typical in the following sense. Associated to every simple tensor product $x_1 \otimes x_2$ is the rank one operator from $H^*_1$ to $H_2$ that maps a given $x^*\in H^*_1$ as

$$ x^* \mapsto x^*(x_1) x_2$$
This mapping defined on simple tensors extends to a linear identification between $H_1 \otimes H_2$ and the space of finite rank operators from $H^*_1$ to $H_2$. This extends to a linear isometry of the Hilbertian tensor product $H_1\widehat{\otimes}H_2$ with the Hilbert space $HS(H^*_1, H_2)$ of Hilbert-Schmidt operators from $H^*_1$ to $H_2$.
\end{enumerate}

\section{Banach algebra}

In mathematics, especially functional analysis, a Banach algebra, named after Stefan Banach, is an associative algebra $A$ over the real or complex numbers (or over a non-archimedean complete normed field) that at the same time is also a Banach space, i.e. normed and complete. The algebra multiplication and the Banach space norm are required to be related by the following inequality:

$$ \forall x,y \in A: \|xy\| \leq \|x\|\|y\|$$
(i.e., the norm of the product is less than or equal to the product of the norms). This ensures that the multiplication operation is continuous. This property is found in the real and complex numbers; for instance, $|-6�5| \leq |-6|�|5|$
If in the above we relax Banach space to normed space the analogous structure is called a normed algebra.

A Banach algebra is called "unital" if it has an identity element for the multiplication whose norm is 1, and "commutative" if its multiplication is commutative. Any Banach algebra $A$ (whether it has an identity element or not) can be embedded isometrically into a unital Banach algebra $A_e$ so as to form a closed ideal of $A_e$. Often one assumes a priori that the algebra under consideration is unital: for one can develop much of the theory by considering $A_e$ and then applying the outcome in the original algebra. However, this is not the case all the time. For example, one cannot define all the trigonometric functions in a Banach algebra without identity.

The theory of real Banach algebras can be very different from the theory of complex Banach algebras. For example, the spectrum of an element of a nontrivial complex Banach algebra can never be empty, whereas in a real Banach algebra it could be empty for some elements.

Banach algebras can also be defined over fields of $p$-adic numbers. This is part of $p$-adic analysis.
\subsection{Examples}

The prototypical example of a Banach algebra is $C_0(X)$, the space of (complex-valued) continuous functions on a locally compact (Hausdorff) space that vanish at infinity. $C_0(X)$ is unital if and only if $X$ is compact. The complex conjugation being an involution, $C_0(X)$ is in fact a $C^*$-algebra. More generally, every $C^*$-algebra is a Banach algebra.

\begin{itemize}
  \item The set of real (or complex) numbers is a Banach algebra with norm given by the absolute value.
  \item The set of all real or complex $n \times n$ matrices becomes a unital Banach algebra if we equip it with a sub-multiplicative matrix norm.
  \item Take the Banach space $R^n$ (or $C^n$) with norm $||x|| = max |x_i|$ and define multiplication componentwise: $$(x_1,...,x_n)(y_1,...,y_n) = (x_1y_1,...,x_ny_n)$$.
  \item The quaternions form a 4-dimensional real Banach algebra, with the norm being given by the absolute value of quaternions.
  \item The algebra of all bounded real- or complex-valued functions defined on some set (with pointwise multiplication and the supremum norm) is a unital Banach algebra.
  \item The algebra of all bounded continuous real- or complex-valued functions on some locally compact space (again with pointwise operations and supremum norm) is a Banach algebra.
  \item The algebra of all continuous linear operators on a Banach space $E$ (with functional composition as multiplication and the operator norm as norm) is a unital Banach algebra. The set of all compact operators on $E$ is a closed ideal in this algebra.
  \item If $G$ is a locally compact Hausdorff topological group and $\mu$ its Haar measure, then the Banach space $L^1(G)$ of all $\mu$-integrable functions on $G$ becomes a Banach algebra under the convolution $xy(g) = \int x(h) y(h^{-1}g) d\mu(h) for x, y in L^1(G)$.
  \item Uniform algebra: A Banach algebra that is a subalgebra of the complex algebra $C(X)$ with the supremum norm and that contains the constants and separates the points of $X$ (which must be a compact Hausdorff space).
  \item Natural Banach function algebra: A uniform algebra whose all characters are evaluations at points of $X$.
  \item $C^*$-algebra: A Banach algebra that is a closed $*$-subalgebra of the algebra of bounded operators on some Hilbert space.
  \item Measure algebra: A Banach algebra consisting of all Radon measures on some locally compact group, where the product of two measures is given by convolution.

\end{itemize}

\subsection{Properties}

Several elementary functions which are defined via power series may be defined in any unital Banach algebra; examples include the exponential function and the trigonometric functions, and more generally any entire function. (In particular, the exponential map can be used to define abstract index groups.) The formula for the geometric series remains valid in general unital Banach algebras. The binomial theorem also holds for two commuting elements of a Banach algebra. The set of invertible elements in any unital Banach algebra is an open set, and the inversion operation on this set is continuous, (and hence homeomorphism) so that it forms a topological group under multiplication. If a Banach algebra has unit 1, then 1 cannot be a commutator; i.e.,$ xy - yx \ne 1 $ for any $x, y \in A$.The various algebras of functions given in the examples above have very different properties from standard examples of algebras such as the reals. For example:
\begin{itemize}
  \item Every real Banach algebra which is a division algebra is isomorphic to the reals, the complexes, or the quaternions. Hence, the only complex Banach algebra which is a division algebra is the complexes. (This is known as the Gelfand-Mazur theorem.)
  \item Every unital real Banach algebra with no zero divisors, and in which every principal ideal is closed, is isomorphic to the reals, the complexes, or the quaternions.
  \item Every commutative real unital Noetherian Banach algebra with no zero divisors is isomorphic to the real or complex numbers.
  \item Every commutative real unital Noetherian Banach algebra (possibly having zero divisors) is finite-dimensional.
  \item Permanently singular elements in Banach algebras are topological divisors of zero, i.e., considering extensions $B$ of Banach algebras $A$ some elements that are singular in the given algebra $A$ have a multiplicative inverse element in a Banach algebra extension $B$. Topological divisors of zero in $A$ are permanently singular in all Banach extension $B$ of $A$.

\end{itemize}

\subsection{Spectral theory}
Unital Banach algebras over the complex field provide a general setting to develop spectral theory. The spectrum of an element $x \in A$, denoted by $\sigma(x)$, consists of all those complex scalars $\lambda$ such that $x - \lambda 1$ is not invertible in $A$. The spectrum of any element $x$ is a closed subset of the closed disc in $C$ with radius $||x||$ and center 0, and thus is compact. Moreover, the spectrum $\sigma(x)$ of an element $x$ is non-empty and satisfies the spectral radius formula:
$$sup\{|\lambda| : \lambda \in \sigma(x)\} = \lim_{n \rightarrow \infty}\|x^n\|^{1/n}$$
Given $x \in A$, the holomorphic functional calculus allows to define $f(x) \in A$ for any function $f$ holomorphic in a neighborhood of $\sigma(x)$. Furthermore, the spectral mapping theorem holds
$$\sigma(f(x)) = f(\sigma(x))$$

When the Banach algebra $A$ is the algebra $L(X)$ of bounded linear operators on a complex Banach space $X$ (e.g., the algebra of square matrices), the notion of the spectrum in A coincides with the usual one in the operator theory. For $f \in C(X)$ (with a compact Hausdorff space $X$), one sees that:

$$\sigma(f) = \{ f(t) : t \in X \}$$.
The norm of a normal element $x$ of a $C^*$-algebra coincides with its spectral radius. This generalizes an analogous fact for normal operators.

Let $A$ be a complex unital Banach algebra in which every non-zero element x is invertible (a division algebra). For every $a \in A$, there is $\lambda \in C$ such that $a - \lambda 1$ is not invertible (because the spectrum of a is not empty) hence $a = \lambda 1$ : this algebra $A$ is naturally isomorphic to $C$ (the complex case of the Gelfand-Mazur theorem).

\subsection{characters}

A character $\chi$ is a linear functional on $A$ which is at the same time multiplicative, $\chi(ab) = \chi(a) \chi(b)$, and satisfies $\chi(1) = 1$. Every character is automatically continuous from $A$ to $C$, since the kernel of a character is a maximal ideal, which is closed. Moreover, the norm (i.e., operator norm) of a character is one. Equipped with the topology of pointwise convergence on $A$ (i.e., the topology induced by the weak-* topology of $A^*$), the character space, $\delta(A)$, is a Hausdorff compact space.

For any $x \in A$,

$$\sigma(x) = \sigma(\hat x)$$
where $\hat x$ is the Gelfand representation of $x$ defined as follows: $\hat x$ is the continuous function from $\delta(A)$ to $C$ given by $\hat x(\chi) = \chi(x)$.? The spectrum of $\hat x$, in the formula above, is the spectrum as element of the algebra $C(\lambda(A))$ of complex continuous functions on the compact space $\lambda(A)$. Explicitly,

$$\sigma(\hat x) = \{ \chi(x) : \chi \in \Delta(A) \}$$
As an algebra, a unital commutative Banach algebra is semisimple (i.e., its Jacobson radical is zero) if and only if its Gelfand representation has trivial kernel. An important example of such an algebra is a commutative $C^*$-algebra. In fact, when $A$ is a commutative unital $C^*$-algebra, the Gelfand representation is then an isometric *-isomorphism between $A$ and $C(\delta(A))$

\section{$C^*$-algebra}

$C^*$-algebras (pronounced "C-star") are an important area of research in functional analysis, a branch of mathematics. A $C^*$-algebra is a complex algebra $A$ of continuous linear operators on a complex Hilbert space with two additional properties:

\begin{itemize}
  \item A is a topologically closed set in the norm topology of operators.
  \item A is closed under the operation of taking adjoints of operators.
\end{itemize}

$C^*$-algebras were first considered primarily for their use in quantum mechanics to model algebras of physical observables. This line of research began with Werner Heisenberg's matrix mechanics and in a more mathematically developed form with Pascual Jordan around 1933. Subsequently John von Neumann attempted to establish a general framework for these algebras which culminated in a series of papers on rings of operators. These papers considered a special class of $C^*$-algebras which are now known as von Neumann algebras.

Around 1943, the work of Israel Gelfand and Mark Naimark yielded an abstract characterisation of $C^*$-algebras making no reference to operators on a Hilbert space.

$C^*$-algebras are now an important tool in the theory of unitary representations of locally compact groups, and are also used in algebraic formulations of quantum mechanics. Another active area of research is the program to obtain classification, or to determine the extent of which classification is possible, for separable simple nuclear $C^*$-algebras.
\subsection{Abstract characterization}
We begin with the abstract characterization of $C^*$-algebras given in the 1943 paper by Gelfand and Naimark.
A $C^*$-algebra, $A$ is a Banach algebra over the field of complex numbers, together with a map $* : A \rightarrow A$. One writes $x^*$ for the image of an element $x$ of $A$. The map * has the following properties
\begin{itemize}
  \item It is an involution, for every $x \in A$
$$ x^{**} = (x^*)^* =  x$$
  \item For all $x, y \in A$:
 $$(x + y)^* = x^* + y^*$$
 $$(x y)^* = y^* x^*$$
  \item For every complex number $\lambda$ in $C$ and every $x \in A$:
 $$(\lambda x)^* = \overline{\lambda} x^*$$
  \item For all $x \in A$:
  $$\|x^* x \| = \|x\|\|x^*\|$$
\end{itemize}

\begin{rem}The first three identities say that $A$ is a *-algebra. The last identity is called the $C^*$ identity and is equivalent to:

$$\|xx^*\| = \|x\|^2$$
\end{rem}

which is sometimes called the $B^*$-identity. For history behind the names $C^*$- and $B^*$-algebras, see the history section below.

The $C^*$-identity is a very strong requirement. For instance, together with the spectral radius formula, it implies that the $C^*$-norm is uniquely determined by the algebraic structure:

 $$\|x\|^2 = \|x^* x\| = \sup\{|\lambda| : x^* x - \lambda \,1 \text{ is not invertible} \}$$
A bounded linear map, $\pi : A \rightarrow B$, between $C^*$-algebras $A$ and $B$ is called a *-homomorphism if
\begin{itemize}
\item For $x, y \in A$
 $$\pi(x y) = \pi(x) \pi(y) $$
\item For $x \in A$
 $$\pi(x^*) = \pi(x)^* $$
 \end{itemize}
In the case of $C^*$-algebras, any *-homomorphism $\pi$ between $C^*$-algebras is non-expansive, i.e. bounded with norm = 1. Furthermore, an injective *-homomorphism between $C^*$-algebras is isometric. These are consequences of the $C^*$-identity.

A bijective *-homomorphism $\pi$ is called a $C^*$-isomorphism, in which case $A$ and $B$ are said to be isomorphic.

\subsection{Some history: $B^*$-algebras and $C^*$-algebras}
The term $B^*$-algebra was introduced by C. E. Rickart in 1946 to describe Banach *-algebras that satisfy the condition:

$\lVert x x^* \rVert = \lVert x \rVert ^2$ for all $x$ in the given $B^*$-algebra. ($B^*$-condition)
This condition automatically implies that the *-involution is isometric, that is, $||x|| = ||x*||$. Hence $||xx*|| = ||x|| ||x*||$, and therefore, a $B^*$-algebra is also a $C^*$-algebra. Conversely, the $C^*$-condition implies the $B^*$-condition. This is nontrivial, and can be proved without using the condition $||x|| = ||x*||$. For these reasons, the term B*-algebra is rarely used in current terminology, and has been replaced by the term 'C*-algebra'.

The term C*-algebra was introduced by I. E. Segal in 1947 to describe norm-closed subalgebras of $B(H)$, namely, the space of bounded operators on some Hilbert space H. 'C' stood for 'closed'. In his paper Segal defines a C*-algebra as a "uniformly closed, self-adjoint algebra of bounded operators on a Hilbert space".
\begin{enumerate}
  \item \textbf{Structure of C*-algebras}
C*-algebras have a large number of properties that are technically convenient. Some of these properties can be established by using the continuous functional calculus or by reduction to commutative C*-algebras. In the latter case, we can use the fact that the structure of these is completely determined by the Gelfand isomorphism.
  \item \textbf{Self-adjoint elements}
Self-adjoint elements are those of the form $x=x^*$. The set of elements of a C*-algebra A of the form $x^*x$ forms a closed convex cone. This cone is identical to the elements of the form $xx^*$. Elements of this cone are called non-negative (or sometimes positive, even though this terminology conflicts with its use for elements of $R$.)

The set of self-adjoint elements of a C*-algebra A naturally has the structure of a partially ordered vector space; the ordering is usually denoted $\geq$. In this ordering, a self-adjoint element $x$ of $A$ satisfies $x = 0$ if and only if the spectrum of $x$ is non-negative, if and only if $x = s*s$ for some $s$. Two self-adjoint elements $x$ and $y$ of $A$ satisfy $x = y$ if $x-y = 0$.

This partially ordered subspace allows the definition of a positive linear functional on a C*-algebra, which in turn is used to define the states of a C*-algebra, which in turn can be used to construct the spectrum of a C*-algebra using the GNS construction
  \item \textbf{Quotients and approximate identities}
Any C*-algebra A has an approximate identity. In fact, there is a directed family $\{e_\lambda\}_{\lambda \in I}$ of self-adjoint elements of $A$ such that

 $$x e_\lambda \rightarrow x$$
 $$0 \leq e_\lambda \leq e_\mu \leq 1\quad \mbox{ whenever } \lambda \leq \mu$$
In case $A$ is separable, $A$ has a sequential approximate identity. More generally, $A$ will have a sequential approximate identity if and only if $A$ contains a strictly positive element, i.e. a positive element $h$ such that $hAh$ is dense in $A$.
Using approximate identities, one can show that the algebraic quotient of a C*-algebra by a closed proper two-sided ideal, with the natural norm, is a C*-algebra.

Similarly, a closed two-sided ideal of a C*-algebra is itself a C*-algebra.

\end{enumerate}

\subsection{Examples}

\begin{enumerate}
  \item \textbf{Finite-dimensional C*-algebras} The algebra $M(n, C)$ of $n \times n$ matrices over $C$ becomes a C*-algebra if we consider matrices as operators on the Euclidean space, $C^n$, and use the operator norm $||.||$ on matrices. The involution is given by the conjugate transpose. More generally, one can consider finite direct sums of matrix algebras. In fact, all C*-algebras that are finite dimensional as vector spaces are of this form, up to isomorphism. The self-adjoint requirement means finite-dimensional C*-algebras are semisimple, from which fact one can deduce the following theorem of Artin�Wedderburn type:

\begin{thm} A finite-dimensional C*-algebra, $A$, is canonically isomorphic to a finite direct sum

 $$A = \bigoplus_{e \in \min A } A_e$$
where min $A$ is the set of minimal nonzero self-adjoint central projections of $A$.
\end{thm}

Each C*-algebra, $A_e$, is isomorphic (in a noncanonical way) to the full matrix algebra $M(dim(e), C)$. The finite family indexed on min $A$ given by $\{dim(e)\}_e$ is called the dimension vector of $A$. This vector uniquely determines the isomorphism class of a finite-dimensional C*-algebra. In the language of $K$-theory, this vector is the positive cone of the $K_0$ group of $A$.

An immediate generalization of finite dimensional C*-algebras are the approximately finite dimensional C*-algebras.

  \item \textbf{C*-algebras of operators}
  The prototypical example of a C*-algebra is the algebra $B(H)$ of bounded (equivalently continuous) linear operators defined on a complex Hilbert space $H$; here $x^*$ denotes the adjoint operator of the operator $x : H \rightarrow H$. In fact, every C*-algebra, $A$, is *-isomorphic to a norm-closed adjoint closed subalgebra of $B(H)$ for a suitable Hilbert space, $H$; this is the content of the Gelfand�Naimark theorem.

  \item \textbf{C*-algebras of compact operators}
  Let $H$ be a separable infinite-dimensional Hilbert space. The algebra $K(H)$ of compact operators on $H$ is a norm closed subalgebra of $B(H)$. It is also closed under involution; hence it is a C*-algebra.

Concrete C*-algebras of compact operators admit a characterization similar to Wedderburn's theorem for finite dimensional C*-algebras:

\begin{thm} If $A$ is a C*-subalgebra of $K(H)$, then there exists Hilbert spaces $\{H_i\}_{i\in I}$ such that

$$ A \cong \bigoplus_{i \in I } K(H_i)$$
where the (C*-)direct sum consists of elements $(T_i)$ of the Cartesian product $\prod K(H_i)$ with $||T_i|| \rightarrow 0$.
\end{thm}
Though $K(H)$ does not have an identity element, a sequential approximate identity for $K(H)$ can be developed. To be specific, $H$ is isomorphic to the space of square summable sequences $l^2$; we may assume that $H =l^2$. For each natural number $n$ let $H_n$ be the subspace of sequences of $l^2$ which vanish for indices $k \leq n$ and let $e_n$ be the orthogonal projection onto $H_n$. The sequence $\{e_n\}_n$ is an approximate identity for $K(H)$.

$K(H)$ is a two-sided closed ideal of $B(H)$. For separable Hilbert spaces, it is the unique ideal. The quotient of $B(H)$ by $K(H)$ is the Calkin algebra.

  \item \textbf{Commutative C*-algebras}
  Let $X$ be a locally compact Hausdorff space. The space $C_0(X)$ of complex-valued continuous functions on $X$ that vanish at infinity (defined in the article on local compactness) form a commutative C*-algebra $C_0(X)$ under pointwise multiplication and addition. The involution is pointwise conjugation. $C_0(X)$ has a multiplicative unit element if and only if $X$ is compact. As does any C*-algebra, $C_0(X)$ has an approximate identity. In the case of $C_0(X)$ this is immediate: consider the directed set of compact subsets of $X$, and for each compact $K$ let $f_{K'}$ be a function of compact support which is identically 1 on $K$. Such functions exist by the Tietze extension theorem which applies to locally compact Hausdorff spaces. Function sequence $\{f_K\}$ is an approximate identity.

The Gelfand representation states that every commutative C*-algebra is *-isomorphic to the algebra $C_0(X)$, where $X$ is the space of characters equipped with the weak* topology. Furthermore if $C_0(X)$ is isomorphic to $C_0(Y)$ as C*-algebras, it follows that $X$ and $Y$ are homeomorphic. This characterization is one of the motivations for the noncommutative topology and noncommutative geometry programs.

  \item \textbf{C*-enveloping algebra}
  Given a Banach *-algebra A with an approximate identity, there is a unique (up to C*-isomorphism) C*-algebra $E(A)$ and *-morphism $\pi$ from $A$ into $E(A)$ which is universal, that is, every other continuous *-morphism $\pi': A \rightarrow B$ factors uniquely through $\pi$. The algebra $E(A)$ is called the C*-enveloping algebra of the Banach *-algebra $A$.

Of particular importance is the C*-algebra of a locally compact group $G$. This is defined as the enveloping C*-algebra of the group algebra of $G$. The C*-algebra of $G$ provides context for general harmonic analysis of $G$ in the case $G$ is non-abelian. In particular, the dual of a locally compact group is defined to be the primitive ideal space of the group C*-algebra. See spectrum of a C*-algebra.

  \item \textbf{Von Neumann algebras}
  Von Neumann algebras, known as $W^*$ algebras before the 1960s, are a special kind of C*-algebra. They are required to be closed in the weak operator topology, which is weaker than the norm topology.

The Sherman�Takeda theorem implies that any C*-algebra has a universal enveloping $W^*$-algebra, such that any homomorphism to a $W^*$-algebra factors through it.
\end{enumerate}

\subsection{Type for C*-algebras}
A C*-algebra A is of type I if and only if for all non-degenerate representations $\pi$ of A the von Neumann algebra $\pi(A)''$ (that is, the bicommutant of $\pi(A))$ is a type I von Neumann algebra. In fact it is sufficient to consider only factor representations, i.e. representations $\pi$ for which $\pi(A)''$ is a factor.

A locally compact group is said to be of type I if and only if its group C*-algebra is type I.

However, if a C*-algebra has non-type I representations, then by results of James Glimm it also has representations of type II and type III. Thus for C*-algebras and locally compact groups, it is only meaningful to speak of type I and non type I properties.

\subsection{C*-algebras and quantum field theory}
In quantum mechanics, one typically describes a physical system with a C*-algebra A with unit element; the self-adjoint elements of $A$ (elements $x$ with $x^* = x$) are thought of as the observables, the measurable quantities, of the system. A state of the system is defined as a positive functional on $A$ (a C-linear map $\varphi : A \rightarrow C$ with $\varphi(u^*u) = 0$ for $all u \in A$) such that $\varphi(1) = 1$. The expected value of the observable $x$, if the system is in state $\varphi$, is then $\varphi(x)$.

This C*-algebra approach is used in the Haag-Kastler axiomatization of local quantum field theory, where every open set of Minkowski spacetime is associated with a C*-algebra.

\section{Spectrum of a $C^*$-algebra}

In mathematics, the spectrum of a $C^*$-algebra or dual of a $C^*$-algebra $A$, denoted $\check{A}$, is the set of unitary equivalence classes of irreducible *-representations of $A$. A *-representation $p$ of $A$ on a Hilbert space $H$ is irreducible if, and only if, there is no closed subspace $K$ different from $H$ and $\{0\}$ which is invariant under all operators $p(x)$ with $x\in A$. We implicitly assume that irreducible representation means non-null irreducible representation, thus excluding trivial (i.e. identically 0) representations on one-dimensional spaces. As explained below, the spectrum $\check{A}$ is also naturally a topological space; this generalizes the notion of the spectrum of a ring.

One of the most important applications of this concept is to provide a notion of dual object for any locally compact group. This dual object is suitable for formulating a Fourier transform and a Plancherel theorem for unimodular separable locally compact groups of type I and a decomposition theorem for arbitrary representations of separable locally compact groups of type I. The resulting duality theory for locally compact groups is however much weaker than the Tannaka�Krein duality theory for compact topological groups or Pontryagin duality for locally compact abelian groups, both of which are complete invariants. That the dual is not a complete invariant is easily seen as the dual of any finite-dimensional full matrix algebra $M_n(C)$ consists of a single point.

\subsection{Primitive spectrum}

The topology of $\check{A}$ can be defined in several equivalent ways. We first define it in terms of the primitive spectrum. The primitive spectrum of A is the set of primitive ideals $Prim(A)$ of $A$, where a primitive ideal is the kernel of an irreducible *-representation. The set of primitive ideals is a topological space with the hull-kernel topology (or Jacobson topology). This is defined as follows: If $X$ is a set of primitive ideals, its hull-kernel closure is

 $$\overline{X} = \left \{\rho \in \operatorname{Prim}(A): \rho \supseteq \bigcap_{\pi \in X} \pi \right \}$$

Hull-kernel closure is easily shown to be an idempotent operation, that is

$$ \overline{\overline{X}} = \overline{X}$$,
and it can be shown to satisfy the Kuratowski closure axioms. As a consequence, it can be shown that there is a unique topology $T$ on $Prim(A)$ such that the closure of a set $X$ with respect to t is identical to the hull-kernel closure of $X$.

Since unitarily equivalent representations have the same kernel, the map $\pi \rightarrow ker(\pi)$ factors through a surjective map

 $$\operatorname{k}: \hat{A} \to \operatorname{Prim}(A)$$
We use the map $k$ to define the topology on $\check{A}$ as follows:
The open sets of $\check{A}$ are inverse images $k^{-1}(U)$ of open subsets $U$ of $Prim(A)$. This is indeed a topology.
The hull-kernel topology is an analogue for non-commutative rings of the Zariski topology for commutative rings.
The topology on $\check{A}$ induced from the hull-kernel topology has other characterizations in terms of states of $A$.

\subsection{Examples}
\begin{enumerate}
  \item \textbf{Commutative $C^*$-algebras}:The spectrum of a commutative $C^*$-algebra $A$ coincides with the Gelfand dual of $A$ (not to be confused with the dual $A'$ of the Banach space $A$). In particular, suppose $X$ is a compact Hausdorff space. Then there is a natural homeomorphism

 $$\operatorname{I}: X \cong \operatorname{Prim}( \operatorname{C}(X))$$
This mapping is defined by

  $$\operatorname{I}(x) = \{f \in \operatorname{C}(X): f(x) = 0 \}$$
$I(x)$ is a closed maximal ideal in $C(X)$ so is in fact primitive. For a commutative $C^*$-algebra,

 $$\hat{A} \cong \operatorname{Prim}(A)$$
  \item \textbf{The $C^*$-algebra of bounded operators}: Let $H$ be a separable Hilbert space. $L(H)$ has two norm-closed *-ideals: $I_0 = \{0\}$ and the ideal $K = K(H)$ of compact operators. Thus as a set, $Prim(L(H)) = \{I_0, K\}$. Now

\begin{itemize}
  \item $\{K\}$ is a closed subset of $Prim(L(H))$.

  \item The closure of $\{I_0\}$ is $Prim(L(H))$.
\end{itemize}

Thus $Prim(L(H))$ is a non-Hausdorff space.

The spectrum of $L(H)$ on the other hand is much larger. There are many inequivalent irreducible representations with kernel $K(H)$ or with kernel $\{0\}$.
  \item \textbf{Finite-dimensional $C^*$-algebras}

Suppose $A$ is a finite-dimensional $C^*$-algebra. It is known $A$ is isomorphic to a finite direct sum of full matrix algebras:

$$ A \cong \bigoplus_{e \in \operatorname{min}(A)} A_e$$,
where $min(A)$ are the minimal central projections of $A$. The spectrum of $A$ is canonically isomorphic to $min(A)$ with the discrete topology. For finite-dimensional $C^*$-algebras, we also have the isomorphism

 $$\hat{A} \cong \operatorname{Prim}(A)$$

\end{enumerate}

\subsection{Other characterizations of the spectrum}

The hull-kernel topology is easy to describe abstractly, but in practice for $C^*$-algebras associated to locally compact topological groups, other characterizations of the topology on the spectrum in terms of positive definite functions are desirable.

In fact, the topology on $ \check{A}$ is intimately connected with the concept of weak containment of representations as is shown by the following:

\begin{thm} Let $S$ be a subset of $ \check{A}$. Then the following are equivalent for an irreducible representation $\pi$;
The equivalence class of p in $ \check{A}$ is in the closure of $S$ Every state associated to $\pi$, that is one of the form
 $$f_\xi(x) = \langle \xi  \mid \pi(x) \xi \rangle$$
with $||\xi|| = 1$, is the weak limit of states associated to representations in $S$.
The second condition means exactly that p is weakly contained in $S$.
\end{thm}
The GNS construction is a recipe for associating states of a $C^*$-algebra A to representations of $A$. By one of the basic theorems associated to the GNS construction, a state $f$ is pure if and only if the associated representation $\pi_f$ is irreducible. Moreover, the mapping $$\kappa:PureState(A) \mapsto \check{A}$$ defined by $f \mapsto \pi_f$ is a surjective map.

From the previous theorem one can easily prove the following;

\begin{thm}The mapping
 $\kappa: \operatorname{PureState}(A) \to \hat{A}$
given by the GNS construction is continuous and open.
\end{thm}

\subsection{The space $Irr_n(A)$}

There is yet another characterization of the topology on $\check{A}$ which arises by considering the space of representations as a topological space with an appropriate pointwise convergence topology. More precisely, let n be a cardinal number and let $H_n$ be the canonical Hilbert space of dimension $n$.

$Irr_n(A)$ is the space of irreducible *-representations of $A$ on $H_n$ with the point-weak topology. In terms of convergence of nets, this topology is defined by $\pi_i \rightarrow \pi$; if and only if

$$\langle \pi_i(x) \xi \mid \eta \rangle \to \langle \pi(x) \xi \mid \eta \rangle \quad \forall \xi, \eta \in H_n \ x \in A$$
It turns out that this topology on $Irr_n(A)$ is the same as the point-strong topology, i.e. $\pi_i \rightarrow \pi$ if and only if

 $$\pi_i(x) \xi \to \pi(x) \xi \quad \mbox{ normwise } \forall \xi \in H_n \ x \in A$$

 \begin{thm} Let $\check{A}_n$ be the subset of $\check{A}$ consisting of equivalence classes of representations whose underlying Hilbert space has dimension $n$. The canonical map $Irr_n(A) \rightarrow  \check{A}_n$ is continuous and open. In particular, $\check{A}_n$ can be regarded as the quotient topological space of $Irr_n(A)$ under unitary equivalence.
 \end{thm}
\begin{rem} The piecing together of the various $\check{A}_n$ can be quite complicated.
\end{rem}

\subsection{Algebraic primitive spectra}
Since a $C^*$-algebra $A$ is a ring, we can also consider the set of primitive ideals of $A$, where $A$ is regarded algebraically. For a ring an ideal is primitive if and only if it is the annihilator of a simple module. It turns out that for a $C^*$-algebra $A$, an ideal is algebraically primitive if and only if it is primitive in the sense defined above.

 \begin{thm}Let A be a $C^*$-algebra. Any algebraically irreducible representation of $A$ on a complex vector space is algebraically equivalent to a topologically irreducible *-representation on a Hilbert space. Topologically irreducible *-representations on a Hilbert space are algebraically isomorphic if and only if they are unitarily equivalent.
 \end{thm}

If $G$ is a locally compact group, the topology on dual space of the group $C^*$-algebra $C^*(G)$ of G is called the Fell topology.

\section{Von Neumann algebra}
In mathematics, a von Neumann algebra or $W^*$-algebra is a $*$-algebra of bounded operators on a Hilbert space that is closed in the weak operator topology and contains the identity operator. They were originally introduced by John von Neumann, motivated by his study of single operators, group representations, ergodic theory and quantum mechanics. His double commutant theorem shows that the analytic definition is equivalent to a purely algebraic definition as an algebra of symmetries.

Two basic examples of von Neumann algebras are as follows. The ring $L^\infty(R)$ of essentially bounded measurable functions on the real line is a commutative von Neumann algebra, which acts by pointwise multiplication on the Hilbert space $L^2(R)$ of square integrable functions. The algebra $B(H)$ of all bounded operators on a Hilbert space $H$ is a von Neumann algebra, non-commutative if the Hilbert space has dimension at least 2.

Von Neumann algebras were first studied by von Neumann (1930) in 1929; he and Francis Murray developed the basic theory, under the original name of rings of operators, in a series of papers written in the 1930s and 1940s (F.J. Murray, J. von Neumann 1936, 1937, 1943; J. von Neumann 1938, 1940, 1943, 1949), reprinted in the collected works of von Neumann (1961).

Introductory accounts of von Neumann algebras are given in the online notes of Jones (2003) and Wassermann (1991) and the books by Dixmier (1981), Schwartz (1967), Blackadar (2005) and Sakai (1971). The three volume work by Takesaki (1979) gives an encyclopedic account of the theory. The book by Connes (1994) discusses more advanced topics.
\subsection{Definitions}

There are three common ways to define von Neumann algebras.

The first and most common way is to define them as weakly closed $*$-algebras of bounded operators (on a Hilbert space) containing the identity. In this definition the weak (operator) topology can be replaced by many other common topologies including the strong, ultrastrong or ultraweak operator topologies. The $*$-algebras of bounded operators that are closed in the norm topology are $C^*$-algebras, so in particular any von Neumann algebra is a $C^*$-algebra.

The second definition is that a von Neumann algebra is a subset of the bounded operators closed under $*$ and equal to its double commutant, or equivalently the commutant of some subset closed under $*$. The von Neumann double commutant theorem (von Neumann 1930) says that the first two definitions are equivalent.

The first two definitions describe a von Neumann algebra concretely as a set of operators acting on some given Hilbert space. Sakai (1971) showed that von Neumann algebras can also be defined abstractly as $C^*$-algebras that have a predual; in other words the von Neumann algebra, considered as a Banach space, is the dual of some other Banach space called the predual. The predual of a von Neumann algebra is in fact unique up to isomorphism. Some authors use "von Neumann algebra" for the algebras together with a Hilbert space action, and "$W^*$-algebra" for the abstract concept, so a von Neumann algebra is a $W^*$-algebra together with a Hilbert space and a suitable faithful unital action on the Hilbert space. The concrete and abstract definitions of a von Neumann algebra are similar to the concrete and abstract definitions of a $C^*$-algebra, which can be defined either as norm-closed $*$-algebras of operators on a Hilbert space, or as Banach $*$-algebras such that $||aa^*||=||a|| ||a^*||$.
\subsection{Terminology}
Some of the terminology in von Neumann algebra theory can be confusing, and the terms often have different meanings outside the subject.
\begin{itemize}
  \item A factor is a von Neumann algebra with trivial center, i.e. a center consisting only of scalar operators.
  \item A finite von Neumann algebra is one which is the direct integral of finite factors (meaning the von Neumann algebra has a faithful normal tracial state $\tau: M \to C$.
  \item A von Neumann algebra that acts on a separable Hilbert space is called separable. Note that such algebras are rarely separable in the norm topology.
  \item The von Neumann algebra generated by a set of bounded operators on a Hilbert space is the smallest von Neumann algebra containing all those operators.
  \item The tensor product of two von Neumann algebras acting on two Hilbert spaces is defined to be the von Neumann algebra generated by their algebraic tensor product, considered as operators on the Hilbert space tensor product of the Hilbert spaces.
\end{itemize}
By forgetting about the topology on a von Neumann algebra, we can consider it a (unital) $*$-algebra, or just a ring. Von Neumann algebras are semihereditary: every finitely generated submodule of a projective module is itself projective. There have been several attempts to axiomatize the underlying rings of von Neumann algebras, including Baer $*$-rings and $AW^*$ algebras. The $*$-algebra of affiliated operators of a finite von Neumann algebra is a von Neumann regular ring. (The von Neumann algebra itself is in general not von Neumann regular.)

\subsection{Commutative von Neumann algebras}
The relationship between commutative von Neumann algebras and measure spaces is analogous to that between commutative $C^*$-algebras and locally compact Hausdorff spaces. Every commutative von Neumann algebra is isomorphic to $L^\infty(X)$ for some measure space $(X, \mu)$ and conversely, for every $\sigma$-finite measure space $X$, the $*$-algebra $L^\infty(X)$ is a von Neumann algebra.

Due to this analogy, the theory of von Neumann algebras has been called noncommutative measure theory, while the theory of $C^*$-algebras is sometimes called noncommutative topology (Connes 1994).
\subsection{Projections}
Operators $E$ in a von Neumann algebra for which $E = EE = E^*$ are called projections; they are exactly the operators which give an orthogonal projection of $H$ onto some closed subspace. A subspace of the Hilbert space $H$ is said to belong to the von Neumann algebra $M$ if it is the image of some projection in $M$. This establishes a $1:1$ correspondence between projections of $M$ and subspaces that belong to $M$. Informally these are the closed subspaces that can be described using elements of $M$, or that $M$ "knows" about.

It can be shown that the closure of the image of any operator in $M$ and the kernel of any operator in $M$ belongs to $M$. Also, the closure of the image under an operator of $M$ of any subspace belonging to $M$ also belongs to $M$. (These results are a consequence of the polar decomposition).

\subsection{Comparison Theory of Projections}
The basic theory of projections was worked out by Murray and von Neumann (1936). Two subspaces belonging to $M$ are called (Murray�von Neumann) equivalent if there is a partial isometry mapping the first isomorphically onto the other that is an element of the von Neumann algebra (informally, if $M$ "knows" that the subspaces are isomorphic). This induces a natural equivalence relation on projections by defining $E$ to be equivalent to $F$ if the corresponding subspaces are equivalent, or in other words if there is a partial isometry of $H$ that maps the image of $E$ isometrically to the image of $F$ and is an element of the von Neumann algebra. Another way of stating this is that $E$ is equivalent to $F$ if $E=uu^*$ and $F=u^*u$ for some partial isometry $u$ in $M$.

The equivalence relation $\sim$ thus defined is additive in the following sense: Suppose $E_1 \sim F_1$ and $E_2 \sim F_2$. If $E_1 \perp E_2$ and $F_1 \perp F_2$, then $E_1 + E_2 \sim F_1 + F_2$. Additivity would not generally hold if one were to require unitary equivalence in the definition of $\sim$, i.e. if we say $E$ is equivalent to $F$ if $u^*Eu = F$ for some unitary $u$.

The subspaces belonging to $M$ are partially ordered by inclusion, and this induces a partial order $=$ of projections. There is also a natural partial order on the set of equivalence classes of projections, induced by the partial order $=$ of projections. If $M$ is a factor, $=$ is a total order on equivalence classes of projections, described in the section on traces below.

A projection (or subspace belonging to $M$) $E$ is said to be a finite projection if there is no projection $F < E$ (meaning $F \leq E$ and $F \neq E$) that is equivalent to $E$. For example, all finite-dimensional projections (or subspaces) are finite (since isometries between Hilbert spaces leave the dimension fixed), but the identity operator on an infinite-dimensional Hilbert space is not finite in the von Neumann algebra of all bounded operators on it, since it is isometrically isomorphic to a proper subset of itself. However it is possible for infinite dimensional subspaces to be finite.

Orthogonal projections are noncommutative analogues of indicator functions in $L^\infty(R)$. $L^\infty(R)$ is the $||.||_\infty$-closure of the subspace generated by the indicator functions. Similarly, a von Neumann algebra is generated by its projections; this is a consequence of the spectral theorem for self-adjoint operators.

The projections of a finite factor form a continuous geometry.

\subsection{Factors}
A von Neumann algebra $N$ whose center consists only of multiples of the identity operator is called a factor. von Neumann (1949) showed that every von Neumann algebra on a separable Hilbert space is isomorphic to a direct integral of factors. This decomposition is essentially unique. Thus, the problem of classifying isomorphism classes of von Neumann algebras on separable Hilbert spaces can be reduced to that of classifying isomorphism classes of factors.

Murray and von Neumann (1936) showed that every factor has one of 3 types as described below. The type classification can be extended to von Neumann algebras that are not factors, and a von Neumann algebra is of type $X$ if it can be decomposed as a direct integral of type $X$ factors; for example, every commutative von Neumann algebra has type $I_1$. Every von Neumann algebra can be written uniquely as a sum of von Neumann algebras of types $I$, $II$, and $III$.

There are several other ways to divide factors into classes that are sometimes used:
\begin{itemize}
  \item A factor is called discrete (or occasionally tame) if it has type $I$, and continuous (or occasionally wild) if it has type $II$ or $III$.
  \item A factor is called semifinite if it has type $I$ or $II$, and purely infinite if it has type $III$.
  \item A factor is called finite if the projection 1 is finite and properly infinite otherwise. Factors of types $I$ and $II$ may be either finite or properly infinite, but factors of type $III$ are always properly infinite.
\end{itemize}
\begin{enumerate}
  \item \textbf{Type $I$ factors}
  A factor is said to be of type $I$ if there is a minimal projection $E \neq 0$, i.e. a projection $E$ such that there is no other projection $F$ with $0 < F < E$. Any factor of type $I$ is isomorphic to the von Neumann algebra of all bounded operators on some Hilbert space; since there is one Hilbert space for every cardinal number, isomorphism classes of factors of type $I$ correspond exactly to the cardinal numbers. Since many authors consider von Neumann algebras only on separable Hilbert spaces, it is customary to call the bounded operators on a Hilbert space of finite dimension $n$ a factor of type $I_n$, and the bounded operators on a separable infinite-dimensional Hilbert space, a factor of type $I_\infty$.
  \item \textbf{Type $II$ factors} A factor is said to be of type $II$ if there are no minimal projections but there are non-zero finite projections. This implies that every projection $E$ can be halved in the sense that there are equivalent projections $F$ and $G$ such that $E = F + G$. If the identity operator in a type $II$ factor is finite, the factor is said to be of type $II_1$; otherwise, it is said to be of type $II_\infty$. The best understood factors of type $II$ are the hyperfinite type $II_1$ factor and the hyperfinite type $II_\infty$ factor, found by Murray and von Neumann (1936). These are the unique hyperfinite factors of types $II_1$ and $II_\infty$; there are an uncountable number of other factors of these types that are the subject of intensive study. Murray and von Neumann (1937) proved the fundamental result that a factor of type $II_1$ has a unique finite tracial state, and the set of traces of projections is $[0,1]$.

A factor of type II8 has a semifinite trace, unique up to rescaling, and the set of traces of projections is $[0,\infty]$. The set of real numbers $\lambda$ such that there is an automorphism rescaling the trace by a factor of $\lambda$ is called the fundamental group of the type $II_\lambda$ factor.

The tensor product of a factor of type $II_1$ and an infinite type $I$ factor has type $II_\infty$, and conversely any factor of type $II_\infty$ can be constructed like this. The fundamental group of a type $II_1$ factor is defined to be the fundamental group of its tensor product with the infinite (separable) factor of type $I$. For many years it was an open problem to find a type $II$ factor whose fundamental group was not the group of positive reals, but Connes then showed that the von Neumann group algebra of a countable discrete group with Kazhdan's property $T$ (the trivial representation is isolated in the dual space), such as $SL(3,Z)$, has a countable fundamental group. Subsequently Sorin Popa showed that the fundamental group can be trivial for certain groups, including the semidirect product of $Z^2$ by $SL(2,Z)$.

An example of a type $II_1$ factor is the von Neumann group algebra of a countable infinite discrete group such that every non-trivial conjugacy class is infinite. McDuff (1969) found an uncountable family of such groups with non-isomorphic von Neumann group algebras, thus showing the existence of uncountably many different separable type $II_1$ factors.
  \item \textbf{Type III factors} Lastly, type $III$ factors are factors that do not contain any nonzero finite projections at all. In their first paper Murray and von Neumann (1936) were unable to decide whether or not they existed; the first examples were later found by von Neumann (1940). Since the identity operator is always infinite in those factors, they were sometimes called type $III_\infty$ in the past, but recently that notation has been superseded by the notation $III_\lambda$, where $\lambda$ is a real number in the interval $[0,1]$. More precisely, if the Connes spectrum (of its modular group) is 1 then the factor is of type $III_0$, if the Connes spectrum is all integral powers of $\lambda$ for $0 < \lambda < 1$, then the type is $III_\lambda$, and if the Connes spectrum is all positive reals then the type is $III_1$. (The Connes spectrum is a closed subgroup of the positive reals, so these are the only possibilities.) The only trace on type $III$ factors takes value $\infty$ on all non-zero positive elements, and any two non-zero projections are equivalent. At one time type $III$ factors were considered to be intractable objects, but Tomita�Takesaki theory has led to a good structure theory. In particular, any type $III$ factor can be written in a canonical way as the crossed product of a type $II_\infty$ factor and the real numbers.
\end{enumerate}

\subsection{The predual}
Any von Neumann algebra $M$ has a predual $M^*$, which is the Banach space of all ultraweakly continuous linear functionals on $M$. As the name suggests, $M$ is (as a Banach space) the dual of its predual. The predual is unique in the sense that any other Banach space whose dual is $M$ is canonically isomorphic to $M^*$. Sakai (1971) showed that the existence of a predual characterizes von Neumann algebras among $C^*$ algebras.

The definition of the predual given above seems to depend on the choice of Hilbert space that $M$ acts on, as this determines the ultraweak topology. However the predual can also be defined without using the Hilbert space that $M$ acts on, by defining it to be the space generated by all positive normal linear functionals on $M$. (Here "normal" means that it preserves suprema when applied to increasing nets of self adjoint operators; or equivalently to increasing sequences of projections.)

The predual $M^*$ is a closed subspace of the dual $M^*$ (which consists of all norm-continuous linear functionals on $M$) but is generally smaller. The proof that $M^*$ is (usually) not the same as $M^*$ is nonconstructive and uses the axiom of choice in an essential way; it is very hard to exhibit explicit elements of $M^*$ that are not in $M^*$. For example, exotic positive linear forms on the von Neumann algebra $l^\infty(Z)$ are given by free ultrafilters; they correspond to exotic $*$-homomorphisms into $C$ and describe the Stone�Cech compactification of $Z$.

\subsection{Examples}
\begin{enumerate}
  \item The predual of the von Neumann algebra $L^\infty(R)$ of essentially bounded functions on $R$ is the Banach space $L_1(R)$ of integrable functions. The dual of $L^\infty(R)$ is strictly larger than $L^1(R)$ For example, a functional on $L^\infty(R)$ that extends the Dirac measure $\delta_0$ on the closed subspace of bounded continuous functions $C^0_b(R)$ cannot be represented as a function in $L^1(R)$.
  \item The predual of the von Neumann algebra $B(H)$ of bounded operators on a Hilbert space $H$ is the Banach space of all trace class operators with the trace norm $||A||= Tr(|A|)$. The Banach space of trace class operators is itself the dual of the $C^*$-algebra of compact operators (which is not a von Neumann algebra).

\end{enumerate}

\subsection{Weights, states, and traces}
Weights and their special cases states and traces are discussed in detail in (Takesaki 1979).

\begin{itemize}
  \item A weight $\omega$ on a von Neumann algebra is a linear map from the set of positive elements (those of the form $a^*a$) to $[0,\infty]$.
  \item A positive linear functional is a weight with $\omega(1)$ finite (or rather the extension of $\omega$ to the whole algebra by linearity).
  \item A state is a weight with $\omega(1) = 1$.
  \item A tracial state is a trace with $\omega(1) = 1$.
  \item A trace is a weight with $\omega(aa^*) = \omega(a^*a)$ for all $a$.

\end{itemize}

Any factor has a trace such that the trace of a non-zero projection is non-zero and the trace of a projection is infinite if and only if the projection is infinite. Such a trace is unique up to rescaling. For factors that are separable or finite, two projections are equivalent if and only if they have the same trace. The type of a factor can be read off from the possible values of this trace as follows:
\begin{itemize}
\item Type $I_n: 0, x, 2x, ....,nx$ for some positive $x$ (usually normalized to be $1/n$ or 1).
\item Type $I_\infty: 0, x, 2x, ....,\infty$ for some positive $x$ (usually normalized to be 1).
\item Type $II_1: [0,x]$ for some positive $x$ (usually normalized to be 1).
\item Type $II_\infty: [0,\infty]$.
\item Type $III: 0,\infty$.
\end{itemize}
If a von Neumann algebra acts on a Hilbert space containing a norm 1 vector $v$, then the functional $a \to (av,v)$ is a normal state. This construction can be reversed to give an action on a Hilbert space from a normal state: this is the GNS construction for normal states.

\subsection{Modules over a factor}
Given an abstract separable factor, one can ask for a classification of its modules, meaning the separable Hilbert spaces that it acts on. The answer is given as follows: every such module H can be given an $M$-dimension $dim M(H)$ (not its dimension as a complex vector space) such that modules are isomorphic if and only if they have the same $M$-dimension. The $M$-dimension is additive, and a module is isomorphic to a subspace of another module if and only if it has smaller or equal $M$-dimension.

A module is called standard if it has a cyclic separating vector. Each factor has a standard representation, which is unique up to isomorphism. The standard representation has an antilinear involution $J$ such that $JMJ = M'$. For finite factors the standard module is given by the GNS construction applied to the unique normal tracial state and the $M$-dimension is normalized so that the standard module has $M$-dimension 1, while for infinite factors the standard module is the module with $M$-dimension equal to $\infty$.

The possible $M$-dimensions of modules are given as follows:
\begin{itemize}
  \item Type In (n finite): The M-dimension can be any of $0/n, 1/n, 2/n, 3/n, ..., \infty$. The standard module has $M$-dimension 1 (and complex dimension $n^2$.)
  \item Type $I\infty$ The $M$-dimension can be any of $0, 1, 2, 3, ..., \infty$ The standard representation of $B(H)$ is $H\otimes H$; its $M$-dimension is $\infty$.
  \item Type $II_1$: The $M$-dimension can be anything in $[0, \infty]$. It is normalized so that the standard module has $M$-dimension 1. The $M$-dimension is also called the coupling constant of the module $H$.
  \item Type $II_\infty$: The $M$-dimension can be anything in $[0, 8]$. There is in general no canonical way to normalize it; the factor may have outer automorphisms multiplying the $M$-dimension by constants. The standard representation is the one with $M$-dimension $\infty$.
  \item Type $III$: The $M$-dimension can be 0 or $\infty$. Any two non-zero modules are isomorphic, and all non-zero modules are standard.
\end{itemize}
\subsection{Amenable von Neumann algebras}
Connes (1976) and others proved that the following conditions on a von Neumann algebra $M$ on a separable Hilbert space $H$ are all equivalent:

$M$ is hyperfinite or AFD or approximately finite dimensional or approximately finite: this means the algebra contains an ascending sequence of finite dimensional subalgebras with dense union. (Warning: some authors use "hyperfinite" to mean "AFD and finite".)
\begin{itemize}
  \item $M$ is amenable: this means that the derivations of $M$ with values in a normal dual Banach bimodule are all inner.
  \item $M$ has Schwartz's property $P$: for any bounded operator $T$ on $H$ the weak operator closed convex hull of the elements $uTu^*$ contains an element commuting with $M$.
  \item $M$ is semidiscrete: this means the identity map from $M$ to $M$ is a weak pointwise limit of completely positive maps of finite rank.
  \item $M$ has property $E$ or the Hakeda-Tomiyama extension property: this means that there is a projection of norm 1 from bounded operators on $H$ to $M'$.
  \item $M$ is injective: any completely positive linear map from any self adjoint closed subspace containing 1 of any unital $C^*$-algebra $A$ to $M$ can be extended to a completely positive map from $A$ to $M$.
  \item There is no generally accepted term for the class of algebras above; Connes has suggested that amenable should be the standard term.
\end{itemize}

The amenable factors have been classified: there is a unique one of each of the types $I_n, I_\infty, II_1, II_\infty, III_\lambda$, for $0 < \lambda \leq 1$, and the ones of type $III_0$ correspond to certain ergodic flows. (For type $III_0$ calling this a classification is a little misleading, as it is known that there is no easy way to classify the corresponding ergodic flows.) The ones of type $I$ and $II_1$ were classified by Murray and von Neumann (1943), and the remaining ones were classified by Connes (1976), except for the type $III_1$ case which was completed by Haagerup.

All amenable factors can be constructed using the group-measure space construction of Murray and von Neumann for a single ergodic transformation. In fact they are precisely the factors arising as crossed products by free ergodic actions of $Z$ or $Z/nZ$ on abelian von Neumann algebras $L^\infty(X)$. Type $I$ factors occur when the measure space $X$ is atomic and the action transitive. When $X$ is diffuse or non-atomic, it is equivalent to $[0,1]$ as a measure space. Type $II$ factors occur when $X$ admits an equivalent finite $(II_1)$ or infinite $(II_\infty)$ measure, invariant under an action of $Z$. Type $III$ factors occur in the remaining cases where there is no invariant measure, but only an invariant measure class: these factors are called Krieger factors.
\subsection{Tensor products of von Neumann algebras}

The Hilbert space tensor product of two Hilbert spaces is the completion of their algebraic tensor product. One can define a tensor product of von Neumann algebras (a completion of the algebraic tensor product of the algebras considered as rings), which is again a von Neumann algebra, and act on the tensor product of the corresponding Hilbert spaces. The tensor product of two finite algebras is finite, and the tensor product of an infinite algebra and a non-zero algebra is infinite. The type of the tensor product of two von Neumann algebras $(I, II, or III)$ is the maximum of their types. The commutation theorem for tensor products states that

$$(M\otimes N)^\prime = M^\prime\otimes N^\prime$$
where $M'$ denotes the commutant of $M$.

The tensor product of an infinite number of von Neumann algebras, if done naively, is usually a ridiculously large non-separable algebra. Instead von Neumann (1938) showed that one should choose a state on each of the von Neumann algebras, use this to define a state on the algebraic tensor product, which can be used to product a Hilbert space and a (reasonably small) von Neumann algebra. Araki and Woods (1968) studied the case where all the factors are finite matrix algebras; these factors are called Araki-Woods factors or ITPFI factors (ITPFI stands for "infinite tensor product of finite type $I$ factors"). The type of the infinite tensor product can vary dramatically as the states are changed; for example, the infinite tensor product of an infinite number of type $I_2$ factors can have any type depending on the choice of states. In particular Powers (1967) found an uncountable family of non-isomorphic hyperfinite type $III_\lambda$ factors for $0 < \lambda < 1$, called Powers factors, by taking an infinite tensor product of type $I_2$ factors, each with the state given by:

$$x\mapsto {\rm Tr}\begin{pmatrix}{1\over \lambda+1}&0\\ 0&{\lambda\over \lambda+1}\\ \end{pmatrix} x$$
All hyperfinite von Neumann algebras not of type III0 are isomorphic to Araki-Woods factors, but there are uncountably many of type $III_0$ that are not.
\subsection{Bimodules and subfactors}

A bimodule (or correspondence) is a Hilbert space $H$ with module actions of two commuting von Neumann algebras. Bimodules have a much richer structure than that of modules. Any bimodule over two factors always gives a subfactor since one of the factors is always contained in the commutant of the other. There is also a subtle relative tensor product operation due to Connes on bimodules. The theory of subfactors, initiated by Vaughan Jones, reconciles these two seemingly different points of view.

Bimodules are also important for the von Neumann group algebra M of a discrete group $G$. Indeed if $V$ is any unitary representation of $G$, then, regarding $G$ as the diagonal subgroup of $G \times G$, the corresponding induced representation on $l^2 (G, V)$ is naturally a bimodule for two commuting copies of $M$. Important representation theoretic properties of $G$ can be formulated entirely in terms of bimodules and therefore make sense for the von Neumann algebra itself. For example Connes and Jones gave a definition of an analogue of Kazhdan's Property $T$ for von Neumann algebras in this way.
\subsection{Non-amenable factors}

Von Neumann algebras of type I are always amenable, but for the other types there are an uncountable number of different non-amenable factors, which seem very hard to classify, or even distinguish from each other. Nevertheless Voiculescu has shown that the class of non-amenable factors coming from the group-measure space construction is disjoint from the class coming from group von Neumann algebras of free groups. Later Narutaka Ozawa proved that group von Neumann algebras of hyperbolic groups yield prime type II1 factors, i.e. ones that cannot be factored as tensor products of type II1 factors, a result first proved by Leeming Ge for free group factors using Voiculescu's free entropy. Popa's work on fundamental groups of non-amenable factors represents another significant advance. The theory of factors "beyond the hyperfinite" is rapidly expanding at present, with many new and surprising results; it has close links with rigidity phenomena in geometric group theory and ergodic theory.

\subsection{Examples}
\begin{itemize}
  \item The essentially bounded functions on a s$\sigma$-finite measure space form a commutative $(type I_1)$ von Neumann algebra acting on the $L^2$ functions. For certain non-$\sigma$-finite measure spaces, usually considered pathological, $L^\infty(X)$ is not a von Neumann algebra; for example, the $\sigma$-algebra of measurable sets might be the countable-cocountable algebra on an uncountable set.
  \item The bounded operators on any Hilbert space form a von Neumann algebra, indeed a factor, of type $I$.
  \item If we have any unitary representation of a group $G$ on a Hilbert space $H$ then the bounded operators commuting with $G$ form a von Neumann algebra $G'$, whose projections correspond exactly to the closed subspaces of $H$ invariant under $G$. Equivalent subrepresentations correspond to equivalent projections in $G'$. The double commutant $G''$ of $G$ is also a von Neumann algebra.
  \item The von Neumann group algebra of a discrete group $G$ is the algebra of all bounded operators on $H = l^2(G)$ commuting with the action of $G$ on $H$ through right multiplication. One can show that this is the von Neumann algebra generated by the operators corresponding to multiplication from the left with an element $g \in G$. It is a factor (of type $II_1$) if every non-trivial conjugacy class of $G$ is infinite (for example, a non-abelian free group), and is the hyperfinite factor of type $II_1$ if in addition $G$ is a union of finite subgroups (for example, the group of all permutations of the integers fixing all but a finite number of elements).
  \item The tensor product of two von Neumann algebras, or of a countable number with states, is a von Neumann algebra as described in the section above.
  \item The crossed product of a von Neumann algebra by a discrete (or more generally locally compact) group can be defined, and is a von Neumann algebra. Special cases are the group-measure space construction of Murray and von Neumann and Krieger factors.
  \item The von Neumann algebras of a measurable equivalence relation and a measurable groupoid can be defined. These examples generalise von Neumann group algebras and the group-measure space construction
\end{itemize}
\subsection{Applications}
Von Neumann algebras have found applications in diverse areas of mathematics like knot theory, statistical mechanics, Quantum field theory, Local quantum physics, Free probability, Noncommutative geometry, representation theory, geometry, and probability.

For instance, $C^*$-algebra provides an alternative axiomatization to probability theory. In this case the method goes by the name of Gelfand�Naimark�Segal construction. This is analogous to the two approaches to measure and integration, where one has the choice to construct measures of sets first and define integrals later, or construct integrals first and define set measures as integrals of characteristic functions.

\section{Gelfand representation}

In mathematics, the Gelfand representation in functional analysis (named after I. M. Gelfand) has two related meanings:
\begin{itemize}
  \item a way of representing commutative Banach algebras as algebras of continuous functions;
  \item the fact that for commutative C*-algebras, this representation is an isometric isomorphism.
\end{itemize}

In the former case, one may regard the Gelfand representation as a far-reaching generalization of the Fourier transform of an integrable function. In the latter case, the Gelfand-Naimark representation theorem is one avenue in the development of spectral theory for normal operators, and generalizes the notion of diagonalizing a normal matrix.
\subsection{Historical remarks}

One of Gelfand's original applications (and one which historically motivated much of the study of Banach algebras) was to give a much shorter and more conceptual proof of a celebrated lemma of Norbert Wiener, characterizing the elements of the group algebras $L^1(R)$ and $\ell^1({\mathbf Z})$ whose translates span dense subspaces in the respective algebras.

\subsection{The model algebra}
For any locally compact Hausdorff topological space $X$, the space $C_0(X)$ of continuous complex-valued functions on X which vanish at infinity is in a natural way a commutative C*-algebra:
\begin{itemize}
  \item The structure of algebra over the complex numbers is obtained by considering the pointwise operations of addition and multiplication.
  \item The involution is pointwise complex conjugation.
  \item The norm is the uniform norm on functions.

\end{itemize}
Note that $A$ is unital if and only if $X$ is compact, in which case $C_0(X)$ is equal to $C(X)$, the algebra of all continuous complex-valued functions on $X$

\subsection{Gelfand representation of a commutative Banach algebra}
Let $A$ be a commutative Banach algebra, defined over the field $C$ of complex numbers. A non-zero algebra homomorphism $\varphi: A \rightarrow C$ is called a character of $A$; the set of all characters of $A$ is denoted by $\Phi_A $.

It can be shown that every character on $A$ is automatically continuous, and hence $\Phi_A $ is a subset of the space $A^*$ of continuous linear functionals on $A$; moreover, when equipped with the relative weak-* topology, $\Phi_A $ turns out to be locally compact and Hausdorff. (This follows from the Banach�Alaoglu theorem.) The space $\Phi_A $ is compact (in the topology just defined) if and only if the algebra $A$ has an identity element.

Given $a \in A$, one defines the function $\widehat{a}:\Phi_A\to{\mathbb C}$ by $\widehat{a}(\phi)=\phi(a)$. The definition of $\Phi_A $ and the topology on it ensure that $\widehat{a}$ is continuous and vanishes at infinity, and that the map $a \mapsto \widehat{a}$ defines a norm-decreasing, unit-preserving algebra homomorphism from $A$ to $C_0(\Phi_A)$. This homomorphism is the Gelfand representation of $A$, and $\widehat{a}$ is the Gelfand transform of the element $a$. In general, the representation is neither injective nor surjective.

In the case where $A$ has an identity element, there is a bijection between $\Phi_A $ and the set of maximal ideals in $A$ (this relies on the Gelfand�Mazur theorem). As a consequence, the kernel of the Gelfand representation $A \rightarrow C_0(\Phi_A)$ may be identified with the Jacobson radical of $A$. Thus the Gelfand representation is injective if and only if $A$ is (Jacobson) semisimple

\subsection{Examples}
\begin{enumerate}
  \item In the case where $A = L^1(R)$, the group algebra of $R$, then $\Phi_A $ is homeomorphic to $R$ and the Gelfand transform of $f \in L^1(R)$ is the Fourier transform $\tilde{f}$.

  \item In the case where $A = L^1(R_+)$, the $L_1$-convolution algebra of the real half-line, then $\Phi_A $ is homeomorphic to $\{z \in C: Re(z) = 0\}$, and the Gelfand transform of an element $f \in L^1(R_+)$ is the Laplace transform ${\mathcal L}f$.

\end{enumerate}
\subsection{The C*-algebra case}

As motivation, consider the special case $A = C_0(X)$. Given $x \in X$, let $\varphi_x \in A^*$ be pointwise evaluation at $x$, i.e. $\varphi_x(f) = f(x)$. Then $\varphi_x$ is a character on $A$, and it can be shown that all characters of $A$ are of this form; a more precise analysis shows that we may identify $\Phi_A $ with $X$, not just as sets but as topological spaces. The Gelfand representation is then an isomorphism

$$C_0(X)\to C_0(\Phi_A)$$
\begin{enumerate}
  \item \textbf{The spectrum of a commutative C*-algebra}
The spectrum or Gelfand space of a commutative C*-algebra $A$, denoted $\check{A}$, consists of the set of non-zero *-homomorphisms from $A$ to the complex numbers. Elements of the spectrum are called characters on $A$. (It can be shown that every algebra homomorphism from $A$ to the complex numbers is automatically a *-homomorphism, so that this definition of the term 'character' agrees with the one above.)
In particular, the spectrum of a commutative C*-algebra is a locally compact Hausdorff space: In the unital case, i.e. where the C*-algebra has a multiplicative unit element 1, all characters f must be unital, i.e. $f(1)$ is the complex number one. This excludes the zero homomorphism. So $\check{A}$ is closed under weak-* convergence and the spectrum is actually compact. In the non-unital case, the weak-* closure of $\check{A}$ is $\check{A} \bigcup \{0\}$, where 0 is the zero homomorphism, and the removal of a single point from a compact Hausdorff space yields a locally compact Hausdorff space.
Note that spectrum is an overloaded word. It also refers to the spectrum $\sigma(x)$ of an element $x$ of an algebra with unit 1, that is the set of complex numbers $r$ for which $x - r 1$ is not invertible in $A$. For unital C*-algebras, the two notions are connected in the following way: $\sigma(x)$ is the set of complex numbers $f(x)$ where $f$ ranges over Gelfand space of $A$. Together with the spectral radius formula, this shows that $\check{A}$ is a subset of the unit ball of $A^*$ and as such can be given the relative weak-* topology. This is the topology of pointwise convergence. $A$ net $\{f_k\}_k$ of elements of the spectrum of $A$ converges to $f$ if and only if for each $x \in A$, the net of complex numbers $\{f_k(x)\}_k$ converges to $f(x)$.
If $A$ is a separable C*-algebra, the weak-* topology is metrizable on bounded subsets. Thus the spectrum of a separable commutative C*-algebra $A$ can be regarded as a metric space. So the topology can be characterized via convergence of sequences.
Equivalently, $\sigma(x)$ is the range of $\gamma(x)$, where $\gamma$ is the Gelfand representation.
  \item \textbf{Statement of the commutative Gelfand-Naimark theorem}
Let $A$ be a commutative C*-algebra and let $X$ be the spectrum of $A$. Let
$$\gamma:A \to C_0(X)$$
be the Gelfand representation defined above.
\begin{thm} The Gelfand map $\gamma$ is an isometric *-isomorphism from $A$ onto $C_0(X)$.\end{thm}
The spectrum of a commutative C*-algebra can also be viewed as the set of all maximal ideals $M$ of $A$, with the hull-kernel topology. (See the earlier remarks for the general, commutative Banach algebra case.) For any such $M$ the quotient algebra $A/M$ is one-dimensional (by the Gelfand-Mazur theorem), and therefore any $a \in A$ gives rise to a complex-valued function on $Y$.

In the case of C*-algebras with unit, the spectrum map gives rise to a contravariant functor from the category of C*-algebras with unit and unit-preserving continuous *-homomorphisms, to the category of compact Hausdorff spaces and continuous maps. This functor is one half of a contravariant equivalence between these two categories (its adjoint being the functor that assigns to each compact Hausdorff space $X$ the C*-algebra $C_0(X)$). In particular, given compact Hausdorff spaces $X$ and $Y$, then $C(X)$ is isomorphic to $C(Y)$ (as a C*-algebra) if and only if $X$ is homeomorphic to $Y$.

The 'full' Gelfand�Naimark theorem is a result for arbitrary (abstract) noncommutative C*-algebras A, which though not quite analogous to the Gelfand representation, does provide a concrete representation of $A$ as an algebra of operators.

\end{enumerate}

\subsection{Applications}
One of the most significant applications is the existence of a continuous functional calculus for normal elements in C*-algebra A: An element $x$ is normal if and only if $x$ commutes with its adjoint $x^*$, or equivalently if and only if it generates a commutative C*-algebra $C^*(x)$. By the Gelfand isomorphism applied to $C^*(x)$ this is *-isomorphic to an algebra of continuous functions on a locally compact space. This observation leads almost immediately to:

 \begin{thm} Let $A$ be a C*-algebra with identity and $x$ an element of $A$. Then there is a *-morphism $f \rightarrow f(x)$ from the algebra of continuous functions on the spectrum $\sigma(x)$ into $A$ such that
\begin{itemize}
  \item It maps 1 to the multiplicative identity of A;
  \item It maps the identity function on the spectrum to x.
\end{itemize}
\end{thm}
This allows us to apply continuous functions to bounded normal operators on Hilbert space.
\section{Noncommutative quantum field theory}
In mathematical physics, noncommutative quantum field theory (or quantum field theory on noncommutative spacetime) is an application of noncommutative mathematics to the spacetime of quantum field theory that is an outgrowth of noncommutative geometry and index theory in which the coordinate functions are noncommutative. One commonly studied version of such theories has the "canonical" commutation relation:

$$[x^{\mu}, x^{\nu}]=i \theta^{\mu \nu}$$
which means that (with any given set of axes), it is impossible to accurately measure the position of a particle with respect to more than one axis. In fact, this leads to an uncertainty relation for the coordinates analogous to the Heisenberg uncertainty principle.

Various lower limits have been claimed for the noncommutative scale, (i.e. how accurately positions can be measured) but there is currently no experimental evidence in favour of such theory or grounds for ruling them out.

One of the novel features of noncommutative field theories is the UV/IR mixing phenomenon in which the physics at high energies affects the physics at low energies which does not occur in quantum field theories in which the coordinates commute.\textit{In theoretical physics, it is usually possible to organize physical phenomena according to the energy scale or distance scale. The theory of renormalization group is based on this paradigm. The short-distance, ultraviolet (UV) physics does not directly affect qualitative features of the long-distance, infrared (IR) physics, and vice versa.
This separation of scales holds in quantum field theory. However, in its generalizations such as noncommutative field theory and quantum gravity - string theory in particular - it is expected that interrelations between UV and IR physics start to emerge. In many cases, these interrelations, UV/IR mixing, may be demonstrated explicitly.}

Other features include violation of Lorentz invariance due to the preferred direction of noncommutativity. Relativistic invariance can however be retained in the sense of twisted Poincar\'{e} invariance of the theory. The causality condition is modified from that of the commutative theories.
\subsection{History and motivation}
Heisenberg was the first to suggest extending noncommutativity to the coordinates as a possible way of removing the infinite quantities appearing in field theories before the renormalization procedure was developed and had gained acceptance. The first paper on the subject was published in 1947 by Hartland Snyder. The success of the renormalization method resulted in little attention being paid to the subject for some time. In the 1980s, mathematicians, most notably Alain Connes, developed noncommutative geometry. Among other things, this work generalized the notion of differential structure to a noncommutative setting. This led to an operator algebraic description of noncommutative space-times, with the problem that it classically corresponds to a manifold with positively defined metric tensor, so that there is no description of (noncommutative) causality in this approach. However it also led to the development of a Yang�Mills theory on a noncommutative torus.

The particle physics community became interested in the noncommutative approach because of a paper by Nathan Seiberg and Edward Witten. They argued in the context of string theory that the coordinate functions of the endpoints of open strings constrained to a D-brane in the presence of a constant Neveu-Schwarz B-field�equivalent to a constant magnetic field on the brane�would satisfy the noncommutative algebra set out above. The implication is that a quantum field theory on noncommutative spacetime can be interpreted as a low energy limit of the theory of open strings.

Two papers, one by Sergio Doplicher, Klaus Fredenhagen and John Roberts and the other by D. V. Ahluwalia, set out another motivation for the possible noncommutativity of space-time. The arguments go as follows: According to general relativity, when the energy density grows sufficiently large, a black hole is formed. On the other hand according to the Heisenberg uncertainty principle, a measurement of a space-time separation causes an uncertainty in momentum inversely proportional to the extent of the separation. Thus energy whose scale corresponds to the uncertainty in momentum is localized in the system within a region corresponding to the uncertainty in position. When the separation is small enough, the Schwarzschild radius of the system is reached and a black hole is formed, which prevents any information from escaping the system. Thus there is a lower bound for the measurement of length. A sufficient condition for preventing gravitational collapse can be expressed as an uncertainty relation for the coordinates. This relation can in turn be derived from a commutation relation for the coordinates.

It is worth stressing that, differently from other approaches, in particular those relying upon Connes' ideas, here the noncommutative spacetime is a proper spacetime, i.e. it extends the idea of a four-dimensional pseudo-Riemannian manifold. On the other hand, differently from Connes' noncommutative geometry, the proposed model turns out to be coordinates dependent from scratch. In Doplicher Fredenhagen Roberts' paper noncommutativity of coordinates concerns all four spacetime coordinates and not only spatial ones.

\section{$K$-Theory}
In mathematics, K-theory is, roughly speaking, the study of certain kinds of invariants of large matrices. It originated as the study of a ring generated by vector bundles over a topological space or scheme. In algebraic topology, it is an extraordinary cohomology theory known as topological K-theory. In algebra and algebraic geometry, it is referred to as algebraic K-theory. It is also a fundamental tool in the field of operator algebras.

K-theory involves the construction of families of K-functors that map from topological spaces or schemes to associated rings; these rings reflect some aspects of the structure of the original spaces or schemes. As with functors to groups in algebraic topology, the reason for this functorial mapping is that it is easier to compute some topological properties from the mapped rings than from the original spaces or schemes. Examples of results gleaned from the K-theory approach include Bott periodicity, the Atiyah-Singer index theorem and the Adams operations.

In high energy physics, K-theory and in particular twisted K-theory have appeared in Type II string theory where it has been conjectured that they classify D-branes, Ramond�Ramond field strengths and also certain spinors on generalized complex manifolds. In condensed matter physics K-theory has been used to classify topological insulators, superconductors and stable Fermi surfaces. For more details, see K-theory (physics)
\subsection{Early history}

The subject can be said to begin with Alexander Grothendieck (1957), who used it to formulate his Grothendieck�Riemann�Roch theorem. It takes its name from the German Klasse, meaning "class". Grothendieck needed to work with coherent sheaves on an algebraic variety $X$. Rather than working directly with the sheaves, he defined a group using isomorphism classes of sheaves as generators of the group, subject to a relation that identifies any extension of two sheaves with their sum. The resulting group is called $K(X)$ when only locally free sheaves are used, or $G(X)$ when all are coherent sheaves. Either of these two constructions is referred to as the Grothendieck group; $K(X)$ has cohomological behavior and $G(X)$ has homological behavior.

If $X$ is a smooth variety, the two groups are the same. If it is a smooth affine variety, then all extensions of locally free sheaves split, so the group has an alternative definition.
In topology, by applying the same construction to vector bundles, Michael Atiyah and Friedrich Hirzebruch defined $K(X)$ for a topological space $X$ in 1959, and using the Bott periodicity theorem they made it the basis of an extraordinary cohomology theory. It played a major role in the second proof of the Index Theorem (circa 1962). Furthermore this approach led to a noncommutative $K$-theory for $C*$-algebras.
Already in 1955, Jean-Pierre Serre had used the analogy of vector bundles with projective modules to formulate Serre's conjecture, which states that every finitely generated projective module over a polynomial ring is free; this assertion is correct, but was not settled until 20 years later. (Swan's theorem is another aspect of this analogy.)
The other historical origin of algebraic $K$-theory was the work of Whitehead and others on what later became known as Whitehead torsion.
There followed a period in which there were various partial definitions of higher K-theory functors. Finally, two useful and equivalent definitions were given by Daniel Quillen using homotopy theory in 1969 and 1972. A variant was also given by Friedhelm Waldhausen in order to study the algebraic $K$-theory of spaces, which is related to the study of pseudo-isotopies. Much modern research on higher K-theory is related to algebraic geometry and the study of motivic cohomology.
The corresponding constructions involving an auxiliary quadratic form received the general name $L$-theory. It is a major tool of surgery theory.
In string theory the $K$-theory classification of Ramond�Ramond field strengths and the charges of stable $D$-branes was first proposed in 1997.
\subsection{Chern characters}

Chern classes can be used to construct a homomorphism of rings from the topological $K$-theory of a space to (the completion of) its rational cohomology. For a line bundle $L$, the Chern character $ch$ is defined by

$$\operatorname{ch}(L) = \exp(c_{1}(L)) := \sum_{m=0}^\infty \frac{c_1(L)^m}{m!}$$
More generally, if $V = L_1 \oplus ... \oplus L_n$ is a direct sum of line bundles, with first Chern classes $ x_i = c_1(L_i)$, the Chern character is defined additively

 $$\operatorname{ch}(V)  = e^{x_1} + \dots + e^{x_n} :=\sum_{m=0}^\infty \frac{1}{m!}(x_1^m + ... + x_n^m)$$
The Chern character is useful in part because it facilitates the computation of the Chern class of a tensor product. The Chern character is used in the Hirzebruch-Riemann-Roch theorem.

\subsection{Equivariant K-theory}
The equivariant algebraic $K$-theory is an algebraic $K$-theory associated to the category $\operatorname{Coh}^G(X)$ of equivariant coherent sheaves on an algebraic scheme $X$ with action of a linear algebraic group $G$, via Quillen's Q-construction; thus, by definition,

$$K_i^G(X) = \pi_i(B^+ \operatorname{Coh}^G(X))$$
In particular, $K_0^G(C)$ is the Grothendieck group of $\operatorname{Coh}^G(X)$. The theory was developed by R. W. Thomason in 1980s. Specifically, he proved equivariant analogs of fundamental theorems such as the localization theorem.

One of the major examples of this idea is the generalization of topological $K$-theory to noncommutative $C^*$-algebras in the form of operator

$K$-theory. A further development in this is a bivariant version of $K$-theory called $KK$-theory, which has a composition product
$$KK(A,B)\times KK(B,C) \rightarrow KK(A,C) $$

of which the ring structure in ordinary $K$-theory is a special case. The product gives the structure of a category to $KK$. It has been related to correspondences of algebraic varieties.

\section{Spectral triple}
In noncommutative geometry and related branches of mathematics and mathematical physics, a spectral triple is a set of data which encodes a geometric phenomenon in an analytic way. The definition typically involves a Hilbert space, an algebra of operators on it and an unbounded self-adjoint operator, endowed with supplemental structures. It was conceived by Alain Connes who was motivated by the Atiyah-Singer index theorem and sought its extension to 'noncommutative' spaces. Some authors refer to this notion as unbounded K-cycles or as unbounded Fredholm modules
\subsection{Motivation}
A motivating example of spectral triple is given by the algebra of functions on a compact spin manifold, acting on the Hilbert space of $L^2$-spinors, accompanied by the Dirac operator associated to the spin structure. From the knowledge of these objects one is able to recover the original manifold as a metric space: the manifold as a topological space is recovered as the spectrum of the algebra, while the (absolute value of) Dirac operator retains the metric. On the other hand, the phase part of the Dirac operator, in conjunction with the algebra of functions, gives a K-cycle which encodes index-theoretic information. The local index formula expresses the pairing of the K-group of the manifold with this K-cycle in two ways: the 'analytic/global' side involves the usual trace on the Hilbert space and commutators of functions with the phase operator (which corresponds to the 'index' part of the index theorem), while the 'geometric/local' side involves the Dixmier trace and commutators with the Dirac operator (which corresponds to the 'characteristic class integration' part of the index theorem).

Extensions of the index theorem can be considered in cases, typically when one has an action of a group on the manifold, or when the manifold is endowed with a foliation structure, among others. In those cases the algebraic system of the 'functions' which expresses the underlying geometric object is no longer commutative, but one may able to find the space of square integrable spinors (or, sections of a Clifford module) on which the algebra acts, and the corresponding 'Dirac' operator on it satisfying certain boundedness of commutators implied by the pseudo-differential calculus.
Given a spectral triple $(A, H, D)$, one can apply several important operations to it. The most fundamental one is the polar decomposition $D = F|D|$ of $D$ into a self adjoint unitary operator $F$ (the 'phase' of $D$) and a densely defined positive operator $|D|$ (the 'metric' part).
The positive $|D|$ operator defines a metric on the set of pure states on the norm closure of $A$.

\subsection{Definition}
An odd spectral triple is a triple $(A, H, D)$ consisting of a Hilbert space $H$, an algebra $A$ of operators on $H$ (usually closed under taking adjoints) and a densely defined self adjoint operator $D$ satisfying $\|[a, D]\| < \infty$ for any $a \in A$. An even spectral triple is an odd spectral triple with a $Z/2Z$-grading on $H$, such that the elements in $A$ are even while $D$ is odd with respect to this grading. One could also say that an even spectral triple is given by a quartet $(A, H, D, \gamma)$ such that $\gamma$ is a self adjoint unitary on $H$ satisfying $a\gamma = \gamma a$ for any $a \in A$ and $D \gamma = - \gamma D$.

A \textbf{finitely summable} spectral triple is a spectral triple $(A, H, D)$ such that $a.D$ for any $a \in A$ has a compact resolvent which belongs to the class of $L^{p+}$-operators for a fixed $p$ (when $A$ contains the identity operator on $H$, it is enough to require $D-1$ in $L^{p+}(H))$. When this condition is satisfied, the triple $(A, H, D)$ is said to be $p$-summable. A spectral triple is said to be $\theta$-summable when $e^{-tD^2}$ is of trace class for any $t > 0$.

Let $\delta(T)$ denote the commutator of $|D|$ with an operator $T$ on $H$. A spectral triple is said to be regular when the elements in $A$ and the operators of the form $[a, D]$ for $a \in A$ are in the domain of the iterates $\delta_n$ of $\delta$.

When a spectral triple $(A, H, D)$ is $p$-summable, one may define its zeta function $\zeta_D(s) = Tr(|D|^{-s})$; more generally there are zeta functions $\zeta_b(s) = Tr(b|D|^{-s})$ for each element $b$ in the algebra $B$ generated by $\delta_n(A)$ and $\delta_n([a, D])$ for positive integers $n$. They are related to the heat kernel $exp(-t|D|)$ by a Mellin transform. The collection of the poles of the analytic continuation of $\delta_b$ for $b \in B$ is called the dimension spectrum of $(A, H, D)$.

A real spectral triple is a spectral triple $(A, H, D)$ accompanied with an anti-linear involution $J$ on $H$, satisfying $[a, JbJ] = 0$ for $a, b \in A$. In the even case it is usually assumed that J is even with respect to the grading on H.
\subsection{Pairing with K-theory}

The self adjoint unitary $F$ gives a map of the K-theory of $A$ into integers by taking Fredholm index as follows. In the even case, each projection e in $A$ decomposes as $e_0 \oplus e_1$ under the grading and $e_1Fe_0$ becomes a Fredholm operator from $e_0H$ to $e_1H$. Thus $e \rightarrow Ind e_1Fe_0$ defines an additive mapping of $K_0(A)$ to $Z$. In the odd case the eigenspace decomposition of $F$ gives a grading on $H$, and each invertible element in $A$ gives a Fredholm operator $(F + 1) u (F - 1)/4$ from $(F - 1)H$ to $(F + 1)H$. Thus $u \rightarrow Ind (F + 1) u (F - 1)/4$ gives an additive mapping from $K_1(A)$ to $Z$.

When the spectral triple is finitely summable, one may write the above indexes using the (super) trace, and a product of $F, e (resp. u)$ and commutator of $F$ with $e (resp. u)$. This can be encoded as a $(p + 1)$-functional on $A$ satisfying some algebraic conditions and give Hochschild / cyclic cohomology cocycles, which describe the above maps from K-theory to the integers.

\section{Signature operator}
In mathematics, the signature operator is an elliptic differential operator defined on a certain subspace of the space of differential forms on an even-dimensional compact Riemannian manifold, whose analytic index is the same as the topological signature of the manifold if the dimension of the manifold is a multiple of four. It is an instance of a Dirac-type operator.
\subsection{Definition in the even-dimensional case}
Let $M$ be a compact Riemannian manifold of even dimension $2l$. Let

 $$d : \Omega^p(M)\rightarrow \Omega^{p+1}(M)$$
be the exterior derivative on $i^{th}$ order differential forms on $M$. The Riemannian metric on $M$ allows us to define the Hodge star operator $\star$ and with it the inner product

$$\langle\omega,\eta\rangle=\int_M\omega\wedge\star\eta$$
on forms. Denote by

$$ d^*: \Omega^{p+1}(M)\rightarrow \Omega^p(M)$$
the adjoint operator of the exterior differential $d$. This operator can be expressed purely in terms of the Hodge star operator as follows:

$$d^*=  (-1)^{2l(p+1) + 2l + 1} \star d  \star=  - \star d  \star$$
Now consider $d + d^*$ acting on the space of all forms $\Omega(M)=\bigoplus_{p=0}^{2l}\Omega^{p}(M)$. One way to consider this as a graded operator is the following: Let $\tau$ be an involution on the space of all forms defined by:

 $$\tau(\omega)=i^{p(p-1)+l}\star \omega\quad,\quad\omega \in \Omega^p(M)$$
It is verified that $d + d^*$ anti-commutes with $\tau$ and, consequently, switches the $(\pm 1)$-eigenspaces $\Omega_{\pm}(M)$ of $\tau$

Consequently,

$$ d + d^* = \begin{pmatrix} 0 & D \\ D^* & 0 \end{pmatrix}$$
Definition: The operator  $d + d^*$ with the above grading respectively the above operator $D: \Omega_+(M) \rightarrow \Omega_-(M)$  is called the signature operator of $M$. In the odd-dimensional case one defines the signature operator to be $i(d+d^*)\tau$ acting on the even-dimensional forms of $M$.

\subsection{Hirzebruch Signature Theorem}
If  $l = 2k$ , so that the dimension of $M$ is a multiple of four, then Hodge theory implies that:

$$\mathrm{index}(D) = \mathrm{sign}(M)$$
where the right hand side is the topological signature (i.e. the signature of a quadratic form on $H^{2k}(M)$  defined by the cup product).

The Heat Equation approach to the Atiyah-Singer index theorem can then be used to show that:

$$\mathrm{sign}(M) = \int_M L(p_1,\ldots,p_l)$$
where $L$ is the Hirzebruch $L$-Polynomial, and the $p_i$  the Pontrjagin forms on $M$.

\section{Characteristic class}
In mathematics, a characteristic class is a way of associating to each principal bundle on a topological space $X$ a cohomology class of $X$. The cohomology class measures the extent to which the bundle is "twisted"-particularly, whether it possesses sections or not. In other words, characteristic classes are global invariants which measure the deviation of a local product structure from a global product structure. They are one of the unifying geometric concepts in algebraic topology, differential geometry and algebraic geometry.
The notion of characteristic class arose in 1935 in the work of Stiefel and Whitney about vector fields on manifolds.
\subsection{Definition}
Let $G$ be a topological group, and for a topological space $X$, write $b_G(X)$ for the set of isomorphism classes of principal $G$-bundles. This $b_G$ is a contravariant functor from Top (the category of topological spaces and continuous functions) to Set (the category of sets and functions), sending a map f to the pullback operation $f^*$.
A \textbf{characteristic class} $c$ of principal $G$-bundles is then a natural transformation from $b_G$ to a cohomology functor $H^*$, regarded also as a functor to Set. In other words, a characteristic class associates to any principal $G$-bundle $P \rightarrow X$ in $b_G(X)$ an element $c(P)$ in $H^*(X)$ such that, if $f : Y \rightarrow X$ is a continuous map, then $c(f^*P) = f^*c(P)$. On the left is the class of the pullback of $P$ to $Y$; on the right is the image of the class of $P$ under the induced map in cohomology

\subsection{Characteristic numbers}
Characteristic classes are elements of cohomology groups; one can obtain integers from characteristic classes, called characteristic numbers. Respectively: Stiefel-Whitney numbers, Chern numbers, Pontryagin numbers, and the Euler characteristic.
Given an oriented manifold $M$ of dimension $n$ with fundamental class $[M] \in H_n(M)$, and a $G$-bundle with characteristic classes $c_1,\dots,c_k$, one can pair a product of characteristic classes of total degree $n$ with the fundamental class. The number of distinct characteristic numbers is the number of monomials of degree $n$ in the characteristic classes, or equivalently the partitions of $n$ into $\mbox{deg}\,c_i$.
Formally, given $i_1,\dots,i_l$ such that $\sum \mbox{deg}\,c_{i_j} = n$, the corresponding characteristic number is:
$$c_{i_1}\smile c_{i_2}\smile \dots \smile c_{i_m}([M])$$
where $\smile$ denotes the cup product of cohomology classes. These are notated various as either the product of characteristic classes, such as $c_1^2$ or by some alternative notation, such as $P_{1,1}$ for the Pontryagin number corresponding to $p_1^2$, or $\chi$ for the Euler characteristic.
From the point of view of de Rham cohomology, one can take differential forms representing the characteristic classes, take a wedge product so that one obtains a top dimensional form, then integrates over the manifold; this is analogous to taking the product in cohomology and pairing with the fundamental class.
This also works for non-orientable manifolds, which have a $\mathbf{Z}/2\mathbf{Z}$-orientation, in which case one obtains $\mathbf{Z}/2\mathbf{Z}$-valued characteristic numbers, such as the Stiefel-Whitney numbers.
Characteristic numbers solve the oriented and unoriented bordism questions: two manifolds are (respectively oriented or unoriented) cobordant if and only if their characteristic numbers are equal.
\subsection{Motivation}
Characteristic classes are in an essential way phenomena of cohomology theory � they are contravariant constructions, in the way that a section is a kind of function on a space, and to lead to a contradiction from the existence of a section we do need that variance. In fact cohomology theory grew up after homology and homotopy theory, which are both covariant theories based on mapping into a space; and characteristic class theory in its infancy in the 1930s (as part of obstruction theory) was one major reason why a 'dual' theory to homology was sought. The characteristic class approach to curvature invariants was a particular reason to make a theory, to prove a general Gauss-Bonnet theorem.

When the theory was put on an organised basis around 1950 (with the definitions reduced to homotopy theory) it became clear that the most fundamental characteristic classes known at that time (the Stiefel-Whitney class, the Chern class, and the Pontryagin classes) were reflections of the classical linear groups and their maximal torus structure. What is more, the Chern class itself was not so new, having been reflected in the Schubert calculus on Grassmannians, and the work of the Italian school of algebraic geometry. On the other hand there was now a framework which produced families of classes, whenever there was a vector bundle involved.

The prime mechanism then appeared to be this: Given a space X carrying a vector bundle, that implied in the homotopy category a mapping from $X$ to a classifying space $BG$, for the relevant linear group $G$. For the homotopy theory the relevant information is carried by compact subgroups such as the orthogonal groups and unitary groups of $G$. Once the cohomology $H^*(BG)$ was calculated, once and for all, the contravariance property of cohomology meant that characteristic classes for the bundle would be defined in $H^*(X)$ in the same dimensions. For example the Chern class is really one class with graded components in each even dimension.

This is still the classic explanation, though in a given geometric theory it is profitable to take extra structure into account. When cohomology became 'extraordinary' with the arrival of K-theory and cobordism theory from 1955 onwards, it was really only necessary to change the letter $H$ everywhere to say what the characteristic classes were.

Characteristic classes were later found for foliations of manifolds; they have (in a modified sense, for foliations with some allowed singularities) a classifying space theory in homotopy theory.

In later work after the rapprochement of mathematics and physics, new characteristic classes were found by Simon Donaldson and Dieter Kotschick in the instanton theory. The work and point of view of Chern have also proved important: see Chern-Simons theory.
\subsection{Stability}
In the language of stable homotopy theory, the Chern class, Stiefel-Whitney class, and Pontryagin class are stable, while the Euler class is unstable.

Concretely, a stable class is one that does not change when one adds a trivial bundle: $c(V \oplus 1) = c(V)$. More abstractly, it means that the cohomology class in the classifying space for $BG(n)$ pulls back from the cohomology class in $BG(n+1)$ under the inclusion $BG(n) \to BG(n+1)$ (which corresponds to the inclusion $\mathbf{R}^n \to \mathbf{R}^{n+1}$ and similar). Equivalently, all finite characteristic classes pull back from a stable class in $BG$.

This is not the case for the Euler class, as detailed there, not least because the Euler class of a $k$-dimensional bundle lives in $H^k(X)$ (hence pulls back from $H^k(BO(k))$, so it can�t pull back from a class in $H^{k+1}$, as the dimensions differ

\section{Noncommutative algebraic geometry}
Noncommutative algebraic geometry is a branch of mathematics, and more specifically a direction in noncommutative geometry that studies the geometric properties of formal duals of non-commutative algebraic objects such as rings as well as geometric objects derived from them (e.g. by gluing along localizations, taking noncommutative stack quotients etc.). For example, noncommutative algebraic geometry is supposed to extend a notion of an algebraic scheme by suitable gluing of spectra of noncommutative rings; depending on how literally and how generally this aim (and a notion of spectrum) is understood in noncommutative setting, this has been achieved in various level of success. The noncommutative ring generalizes here a commutative ring of regular functions on a commutative scheme. Functions on usual spaces in the traditional (commutative) algebraic geometry multiply by points; as the values of these functions commute, the functions also commute: a times b equals b times a. It is remarkable that viewing noncommutative associative algebras as algebras of functions on "noncommutative" would-be space is a far reaching geometric intuition, though it formally looks like a fallacy.

Much of motivations for noncommutative geometry, and in particular for the noncommutative algebraic geometry is from physics; especially from quantum physics, where the algebras of observables are indeed viewed as noncommutative analogues of functions, hence why not looking at their geometric aspects.

One of the values of the field is that it also provides new techniques to study objects in commutative algebraic geometry such as Brauer groups.

The methods of noncommutative algebraic geometry are analogs of the methods of commutative algebraic geometry, but frequently the foundations are different. Local behavior in commutative algebraic geometry is captured by commutative algebra and especially the study of local rings. These do not have a ring-theoretic analogues in the noncommutative setting; though in a categorical setup one can talk about stacks of local categories of quasi-coherent sheaves over noncommutative spectra. Global properties such as those arising from homological algebra and K-theory more frequently carry over to the noncommutative setting
\subsection{Modern viewpoint via categories of sheaves}

In modern times, one accepts a paradigm implicit in Pierre Gabriel's thesis and partly justified by Gabriel-Rosenberg reconstruction theorem (after Pierre Gabriel and Alexander Rosenberg) that a commutative scheme can be reconstructed, up to isomorphism of schemes, solely from the abelian category of quasicoherent sheaves on the scheme. Alexander Grothendieck taught us that to do a geometry one does not need a space, it is enough to have a category of sheaves on that would be space; this idea has been transmitted to noncommutative algebra via Yuri Manin. There are, a bit weaker, reconstruction theorems from the derived categories of (quasi)coherent sheaves motivating the derived noncommutative algebraic geometry.

\subsection{Non-commutative deformations of commutative rings}
As a motivating example, consider the one-dimensional Weyl algebra over the complex numbers $C$. This is the quotient of the free ring $C\{x, y\}$ by the relation

$$xy - yx = 1$$
This ring represents the polynomial differential operators in a single variable $x; y$ stands in for the differential operator $\partial_x$. This ring fits into a one-parameter family given by the relations $xy - yx = a$. When a is not zero, then this relation determines a ring isomorphic to the Weyl algebra. When a is zero, however, the relation is the commutativity relation for $x$ and $y$, and the resulting quotient ring is the polynomial ring in two variables, $C[x, y]$. Geometrically, the polynomial ring in two variables represents the two-dimensional affine space $A^2$, so the existence of this one-parameter family says that affine space admits non-commutative deformations to the space determined by the Weyl algebra. In fact, this deformation is related to the symbol of a differential operator and the fact that $A^2$ is the cotangent bundle of the affine line.

Studying the Weyl algebra can lead to information about affine space: The Dixmier conjecture about the Weyl algebra is equivalent to the Jacobian conjecture about affine space.

\subsection{Non-commutative localization}
Commutative algebraic geometry begins by constructing the spectrum of a ring. The points of the algebraic variety (or more generally, scheme) are the prime ideals of the ring, and the functions on the algebraic variety are the elements of the ring. A noncommutative ring, however, may not have any proper non-zero two-sided prime ideals. For instance, this is true of the Weyl algebra of polynomial differential operators on affine space: The Weyl algebra is a simple ring. Furthermore, the theory of non-commutative localization and descent theory is much more difficult in the non-commutative setting than in the commutative setting. While it works sometimes, there are rings that cannot be localized in the required fashion. Nevertheless, it is possible to prove some theorems in this setting.
\subsection{Proj of a noncommutative ring}
One of the basic constructions in commutative algebraic geometry is the $Proj$ of a graded commutative ring. This construction builds a projective algebraic variety together with a very ample line bundle whose homogeneous coordinate ring is the original ring. Building the underlying topological space of the variety requires localizing the ring, but building sheaves on that space does not. By a theorem of Jean-Pierre Serre, quasi-coherent sheaves on $Proj$ of a graded ring are the same as graded modules over the ring up to finite dimensional factors. The philosophy of topos theory promoted by Alexander Grothendieck says that the category of sheaves on a space can serve as the space itself. Consequently, in non-commutative algebraic geometry one often defines $Proj$ in the following fashion: Let $R$ be a graded $C$-algebra, and let $Mod-R$ denote the category of graded right $R$-modules. Let $F$ denote the subcategory of $Mod-R$ consisting of all modules of finite length. $Proj$ $R$ is defined to be the quotient of the abelian category $Mod-R$ by $F$. Equivalently, it is a localization of $Mod-R$ in which two modules become isomorphic if, after taking their direct sums with appropriately chosen objects of $F$, they are isomorphic in $Mod-R$.
This approach leads to a theory of non-commutative projective geometry. A non-commutative smooth projective curve turns out to be a smooth commutative curve, but for singular curves or smooth higher-dimensional spaces, the non-commutative setting allows new objects.

\section{Sheaf}
In mathematics, a sheaf is a tool for systematically tracking locally defined data attached to the open sets of a topological space. The data can be restricted to smaller open sets, and the data assigned to an open set is equivalent to all collections of compatible data assigned to collections of smaller open sets covering the original one. For example, such data can consist of the rings of continuous or smooth real-valued functions defined on each open set. Sheaves are by design quite general and abstract objects, and their correct definition is rather technical. They exist in several varieties such as sheaves of sets or sheaves of rings, depending on the type of data assigned to open sets.

There are also maps (or morphisms) from one sheaf to another; sheaves (of a specific type, such as sheaves of abelian groups) with their morphisms on a fixed topological space form a category. On the other hand, to each continuous map there is associated both a direct image functor, taking sheaves and their morphisms on the domain to sheaves and morphisms on the codomain, and an inverse image functor operating in the opposite direction. These functors, and certain variants of them, are essential parts of sheaf theory.

Due to their general nature and versatility, sheaves have several applications in topology and especially in algebraic and differential geometry. First, geometric structures such as that of a differentiable manifold or a scheme can be expressed in terms of a sheaf of rings on the space. In such contexts several geometric constructions such as vector bundles or divisors are naturally specified in terms of sheaves. Second, sheaves provide the framework for a very general cohomology theory, which encompasses also the "usual" topological cohomology theories such as singular cohomology. Especially in algebraic geometry and the theory of complex manifolds, sheaf cohomology provides a powerful link between topological and geometric properties of spaces. Sheaves also provide the basis for the theory of D-modules, which provide applications to the theory of differential equations. In addition, generalisations of sheaves to more general settings than topological spaces, such as Grothendieck topology, have provided applications to mathematical logic and number theory.

\subsection{Introduction}
In topology, differential geometry and algebraic geometry, several structures defined on a topological space (e.g., a differentiable manifold) can be naturally localised or restricted to open subsets of the space: typical examples include continuous real or complex-valued functions, n times differentiable (real or complex-valued) functions, bounded real-valued functions, vector fields, and sections of any vector bundle on the space.

Presheaves formalise the situation common to the examples above: a presheaf (of sets) on a topological space is a structure that associates to each open set $U$ of the space a set $F(U)$ of sections on $U$, and to each open set $V$ included in $U$ a map $F(U) \to F(V)$ giving restrictions of sections over $U$ to $V$. Each of the examples above defines a presheaf with by taking the restriction maps to be the usual restriction of functions, vector fields and sections of a vector bundle. Moreover, in each of these examples the sets of sections have additional algebraic structure: pointwise operations make them abelian groups, and in the examples of real and complex-valued functions the sets of sections even have a ring structure. In addition, in each example the restriction maps are homomorphisms of the corresponding algebraic structure. This observation leads to the natural definition of presheaves with additional algebraic structure such as presheaves of groups, of abelian groups, of rings: sets of sections are required to have the specified algebraic structure, and the restrictions are required to be homomorphisms. Thus for example continuous real-valued functions on a topological space form a presheaf of rings on the space.

Given a presheaf, a natural question to ask is to what extent its sections over an open set $U$ are specified by their restrictions to smaller open sets $V_i$ of an open cover of $U$. A presheaf is separated if its sections are "locally determined": whenever two sections over $U$ coincide when restricted to each of $V_i$, the two sections are identical. All examples of presheaves discussed above are separated, since in each case the sections are specified by their values at the points of the underlying space. Finally, a separated presheaf is a sheaf if compatible sections can be glued together, i.e., whenever there is a section of the presheaf over each of the covering sets $V_i$, chosen so that they match on the overlaps of the covering sets, these sections correspond to a (unique) section on $U$, of which they are restrictions. It is easy to verify that all examples above except the presheaf of bounded functions are in fact sheaves: in all cases the criterion of being a section of the presheaf is local in a sense that it is enough to verify it in an arbitrary neighbourhood of each point.

On the other hand, it is clear that a function can be bounded on each set of an (infinite) open cover of a space without being bounded on all of the space; thus bounded functions provide an example of a presheaf that in general fails to be a sheaf. Another example of a presheaf that fails to be a sheaf is the constant presheaf that associates the same fixed set (or abelian group, or ring,...) to each open set: it follows from the gluing property of sheaves that the set of sections on a disjoint union of two open sets is the Cartesian product of the sets of sections over the two open sets. The correct way to define the constant sheaf $F_A$ (associated to for instance a set $A$) on a topological space is to require sections on an open set $U$ to be continuous maps from $U$ to $A$ equipped with the discrete topology; then in particular $F_A(U) = A$ for connected $U$.

Maps between sheaves or presheaves (called morphisms), consist of maps between the sets of sections over each open set of the underlying space, compatible with restrictions of sections. If the presheaves or sheaves considered are provided with additional algebraic structure, these maps are assumed to be homomorphisms. Sheaves endowed with nontrivial endomorphisms, such as the action of an algebraic torus or a Galois group, are of particular interest.

Presheaves and sheaves are typically denoted by capital letters, F being particularly common, presumably for the French word for sheaves, faisceau. Use of script letters such as $\mathcal{F}$ is also common.
\subsection{Formal definitions}
The first step in defining a sheaf is to define a presheaf, which captures the idea of associating data and restriction maps to the open sets of a topological space. The second step is to require the normalisation and gluing axioms. A presheaf that satisfies these axioms is a sheaf.
\begin{enumerate}
  \item \textbf{Presheaves} Let $X$ be a topological space, and let $C$ be a category. Usually $C$ is the category of sets, the category of groups, the category of abelian groups, or the category of commutative rings. A presheaf $F$ on $X$ is a functor with values in $C$ given by the following data:
      \begin{itemize}
        \item For each open set $U$ of $X$, there corresponds an object $F(U)$ in $C$
        \item For each inclusion of open sets $V \subseteq U$, there corresponds a morphism $res_{V,U} : F(U) \to F(V)$ in the category $C$.
      \end{itemize}
      The morphisms $res_{V,U}$ are called restriction morphisms. If $s \in F(U)$, then its restriction $res_{V,U}(s)$ is often denoted $s|_V$ by analogy with restriction of functions. The restriction morphisms are required to satisfy two properties:
      \begin{itemize}
        \item For every open set $U$ of $X$, the restriction morphism $res_{U,U} : F(U) \to F(U)$ is the identity morphism on $F(U)$
        \item If we have three open sets $W \subseteq V \subseteq U$, then the composite $res_{W,V} \circ res_{V,U} = res_{W,U}$.
      \end{itemize}
      Informally, the second axiom says it doesn't matter whether we restrict to $W$ in one step or restrict first to $V$, then to $W$.

There is a compact way to express the notion of a presheaf in terms of category theory. First we define the category of open sets on $X$ to be the category $O(X)$ whose objects are the open sets of $X$ and whose morphisms are inclusions. Then a $C$-valued presheaf on $X$ is the same as a contravariant functor from $O(X)$ to $C$. This definition can be generalized to the case when the source category is not of the form $O(X)$ for any $X$; see presheaf (category theory).

If $F$ is a $C$-valued presheaf on $X$, and $U$ is an open subset of $X$, then $F(U)$ is called the sections of $F$ over $U$. If $C$ is a concrete category, then each element of $F(U)$ is called a section. A section over $X$ is called a global section. A common notation (used also below) for the restriction $res_{V,U}(s)$ of a section is $s|_V$. This terminology and notation is by analogy with sections of fiber bundles or sections of the \'{e}tal\'{e} space of a sheaf; see below. $F(U)$ is also often denoted $G(U,F)$, especially in contexts such as sheaf cohomology where $U$ tends to be fixed and $F$ tends to be variable.

  \item \textbf{Sheaves} For simplicity, consider first the case where the sheaf takes values in the category of sets. In fact, this definition applies more generally to the situation where the category is a concrete category whose underlying set functor is conservative, meaning that if the underlying map of sets is a bijection, then the original morphism is an isomorphism.

A sheaf is a presheaf with values in the category of sets that satisfies the following two axioms:

\begin{enumerate}
  \item (Locality) If $(U_i)$ is an open covering of an open set $U$, and if $s,t \in F(U)$ are such that $s|_{U_i} = t|_{U_i}$ for each set $U_i$ of the covering, then $s = t$; and
  \item (Gluing) If $(U_i)$ is an open covering of an open set $U$, and if for each $i$ a section $s_i \in F(U_i)$ is given such that for each pair $U_i,U_j$ of the covering sets the restrictions of $s_i$ and $s_j$ agree on the overlaps: $s_i|_{U_i \cap U_j} = s_j|_{U_i \cap U_j}$, then there is a section $s \in F(U)$ such that $s|_{U_i} = s_i$ for each $i$.

\end{enumerate}
The section $s$ whose existence is guaranteed by axiom 2 is called the gluing, concatenation, or collation of the sections $s_i$. By axiom 1 it is unique. Sections $s_i$ satisfying the condition of axiom 2 are often called compatible; thus axioms 1 and 2 together state that compatible sections can be uniquely glued together. A separated presheaf, or monopresheaf, is a presheaf satisfying axiom 1.

If $C$ has products, the sheaf axioms are equivalent to the requirement that, for any open covering $U_i$, the first arrow in the following diagram is an equalizer:

$$F(U) \rightarrow \prod_{i} F(U_i) {{{} \atop \longrightarrow}\atop{\longrightarrow \atop {}}} \prod_{i, j} F(U_i \cap U_j)$$
Here the first map is the product of the restriction maps

$$\operatorname{res}_{U_i, U} \colon F(U) \rightarrow F(U_i)$$
and the pair of arrows the products of the two sets of restrictions

$$\operatorname{res}_{U_i \cap U_j, U_i} \colon F(U_i) \rightarrow F(U_i \cap U_j)$$
and

$$\operatorname{res}_{U_i \cap U_j, U_j} \colon F(U_j) \rightarrow F(U_i \cap U_j)$$
For a separated presheaf, the first arrow need only be injective.

In general, for an open set $U$ and open covering $(U_i)$, construct a category $J$ whose objects are the sets $U_i$ and the intersections $U_i \cap U_j$ and whose morphisms are the inclusions of $U_i \cap U_j$ in $U_i$ and $U_j$. The sheaf axioms for $U$ and $(U_i)$ are that the limit of the functor $F$ restricted to the category $J$ must be isomorphic to $F(U)$.

Notice that the empty subset of a topological space is covered by the empty family of sets. The product of an empty family or the limit of an empty family is a terminal object, and consequently the value of a sheaf on the empty set must be a terminal object. If sheaf values are in the category of sets, applying the local identity axiom to the empty family shows that over the empty set, there is at most one section, and applying the gluing axiom to the empty family shows that there is at least one section. This property is called the normalisation axiom.

It can be shown that to specify a sheaf, it is enough to specify its restriction to the open sets of a basis for the topology of the underlying space. Moreover, it can also be shown that it is enough to verify the sheaf axioms above relative to the open sets of a covering. Thus a sheaf can often be defined by giving its values on the open sets of a basis, and verifying the sheaf axioms relative to the basis. (see gluing axiom-Sheaves on a basis of open sets.)
  \item \textbf{Morphisms} Heuristically speaking, a morphism of sheaves is analogous to a function between them. However, because sheaves contain data relative to every open set of a topological space, a morphism of sheaves is defined as a collection of functions, one for each open set, that satisfy a compatibility condition.

Let $F$ and $G$ be two sheaves on $X$ with values in the category $C$. A morphism $\varphi : G \to F$ consists of a morphism $\varphi(U) : G(U) \to F(U)$ for each open set $U$ of $X$, subject to the condition that this morphism is compatible with restrictions. In other words, for every open subset $V$ of an open set $U$ the following diagram

\begin{figure}
\begin{center}
\includegraphics[width=4cm]{SheafMorphism-01a.jpg}
\caption{}
\end{center}
\end{figure}

is commutative.

Recall that we could also express a sheaf as a special kind of functor. In this language, a morphism of sheaves is a natural transformation of the corresponding functors. With this notion of morphism, there is a category of $C$-valued sheaves on $X$ for any $C$. The objects are the $C$-valued sheaves, and the morphisms are morphisms of sheaves. An isomorphism of sheaves is an isomorphism in this category.

It can be proved that an isomorphism of sheaves is an isomorphism on each open set $U$. In other words, $\varphi$ is an isomorphism if and only if for each $U$, $\varphi(U)$ is an isomorphism. The same is true of monomorphisms, but not of epimorphisms. See sheaf cohomology.

Notice that we did not use the gluing axiom in defining a morphism of sheaves. Consequently, the above definition makes sense for presheaves as well. The category of $C$-valued presheaves is then a functor category, the category of contravariant functors from $O(X)$ to $C$.
\end{enumerate}

\subsection{Examples}

Because sheaves encode exactly the data needed to pass between local and global situations, there are many examples of sheaves occurring throughout mathematics. Here are some additional examples of sheaves:
\begin{itemize}
  \item Any continuous map of topological spaces determines a sheaf of sets. Let $f : Y \to X$ be a continuous map. We define a sheaf $\Gamma(Y/X)$ on $X$ by setting $\Gamma(Y/X)(U)$ equal to the sections $U \to Y$, that is, $\Gamma(Y/X)(U)$ is the set of all continuous functions $s:U \to Y$ such that $f \circ s = id_U$. Restriction is given by restriction of functions. This sheaf is called the sheaf of sections of $f$, and it is especially important when $f$ is the projection of a fiber bundle onto its base space. Notice that if the image of f does not contain $U$, then $\Gamma(Y/X)(U)$ is empty. For a concrete example, take $X = C \setminus \{0\}, Y = C$, and $f(z) = exp(z). \Gamma(Y/X)(U)$ is the set of branches of the logarithm on $U$.
  \item Fix a point $x \in X$ and an object $S$ in a category $C$. The skyscraper sheaf over $x$ with stalk $S$ is the sheaf $S_x$ defined as follows: If $U$ is an open set containing $x$, then $S_x(U) = S$. If $U$ does not contain $x$, then $S_x(U)$ is the terminal object of $C$. The restriction maps are either the identity on $S$, if both open sets contain $x$, or the unique map from $S$ to the terminal object of $C$
\end{itemize}
\subsection{Turning a presheaf into a sheaf}

It is frequently useful to take the data contained in a presheaf and to express it as a sheaf. It turns out that there is a best possible way to do this. It takes a presheaf $F$ and produces a new sheaf $aF$ called the sheaving, sheafification or sheaf associated to the presheaf $F$. a is called the sheaving functor, sheafification functor, or associated sheaf functor. There is a natural morphism of presheaves $i : F \to aF$ that has the universal property that for any sheaf $G$ and any morphism of presheaves $f : F \to G$, there is a unique morphism of sheaves $\tilde{f} : aF \rightarrow G$ such that $f = \tilde{f} i$. In fact a is the left adjoint functor to the inclusion functor (or forgetful functor) from the category of sheaves to the category of presheaves, and $i$ is the unit of the adjunction. In this way, the category of sheaves turns into a Giraud subcategory of presheaves.

One concrete way of constructing the sheaf $aF$ is to identify it with the sheaf of sections of an appropriate topological space. This space is analogous to the \'{e}tal\'{e} space of a sheaf. Briefly, the underlying set of the topological space is the disjoint union of the stalks of $F$, denoted $Spe  F$. There is a natural map $\varphi : Spe F \to X$ that sends each germ to the point of $x$ over which it lies. For each open set $U$ and each section s of $F$ over $U$, we define a section $\bar{s}$ of $\varphi$ that sends $x$ to the germ $s_x$. Then $Spe  F$ is given the finest topology for which all sections $\bar{s}$ are continuous, and $aF$ is the sheaf of continuous sections of $\varphi$ for this topology.

There are other constructions of the sheaf $aF$. In particular, Grothendieck and Verdier (SGA 4 II 3.0.5) define a functor $L$ from presheaves to presheaves which, when applied to a presheaf, yields a separated presheaf and, when applied to a separated presheaf, yields a sheaf. Applying the functor $L$ twice therefore turns a presheaf into a sheaf, and in fact $LLF$ is the associated sheaf $aF$.

\subsection{Operations}
If $K$ is a subsheaf of a sheaf $F$ of abelian groups, then the quotient sheaf $Q$ is the sheaf associated to the presheaf $U \mapsto F(U)/K(U)$; in other words, the quotient sheaf fits into an exact sequence of sheaves of abelian groups;

$$0 \to K \to F \to Q \to 0$$
(this is also called a sheaf extension.)

Let $F, G$ be sheaves of abelian groups. The set of morphisms of sheaves from $F$ to $G$ forms an abelian group (by the abelian group structure of $G$). The sheaf hom of $F$ and $G$, denoted by,

$$\mathcal{H}om(F, G)$$
is the sheaf of abelian groups $U \mapsto \operatorname{Hom}(F|_U, G|_U)$ where $F|_U$ is the sheaf on $U$ given by $(F|_U)(V) = F(V)$ (Note sheafification is not needed here.) The tensor product of $F$ and $G$ is the sheaf associated to the presheaf $U \mapsto F(U) \otimes G(U)$.

All of these operations extend to sheaves of modules over a sheaf of rings $A$; the above is the special case when $A$ is the constant sheaf $\underline{\mathbf{Z}}$.

\subsection{Images of sheaves}
The definition of a morphism on sheaves makes sense only for sheaves on the same space $X$. This is because the data contained in a sheaf is indexed by the open sets of the space. If we have two sheaves on different spaces, then their data is indexed differently. There is no way to go directly from one set of data to the other.

However, it is possible to move a sheaf from one space to another using a continuous function. Let $f : X \to Y$ be a continuous function from a topological space $X$ to a topological space $Y$. If we have a sheaf on $X$, we can move it to $Y$, and vice versa. There are four ways in which sheaves can be moved.
\begin{itemize}
  \item A sheaf $\mathcal{F}$ on $X$ can be moved to $Y$ using the direct image functor $f_*$ or the direct image with proper support functor $f_!$.
  \item A sheaf $\mathcal{G}$ on $Y$ can be moved to $X$ using the inverse image functor $f^{-1}$ or the twisted inverse image functor $f^!$.

\end{itemize}
The twisted inverse image functor $f^!$ is, in general, only defined as a functor between derived categories. These functors come in adjoint pairs: $f^{-1}$ and $f_*$ are left and right adjoints of each other, and $Rf_!$ and $f^!$ are left and right adjoints of each other. The functors are intertwined with each other by Grothendieck duality and Verdier duality.

There is a different inverse image functor for sheaves of modules over sheaves of rings. This functor is usually denoted $f^*$ and it is distinct from $f^{-1}$. See inverse image functor.

\subsection{Stalks of a sheaf}
The stalk $\mathcal{F}_x$ of a sheaf $\mathcal{F}$ captures the properties of a sheaf "around" a point $x \in X$. Here, "around" means that, conceptually speaking, one looks at smaller and smaller neighborhoods of the point. Of course, no single neighborhood will be small enough, so we will have to take a limit of some sort.

The stalk is defined by

$$\mathcal{F}_x = \varinjlim_{U\ni x} \mathcal{F}(U)$$
the direct limit being over all open subsets of $X$ containing the given point $x$. In other words, an element of the stalk is given by a section over some open neighborhood of $x$, and two such sections are considered equivalent if their restrictions agree on a smaller neighborhood.

The natural morphism $F(U) \to F_x$ takes a section $s$ in $F(U)$ to its germ. This generalises the usual definition of a germ.

A different way of defining the stalk is

$$\mathcal{F}_x := i^{-1}\mathcal{F}(\{x\})$$
where $i$ is the inclusion of the one-point space $\{x\}$ into $X$. The equivalence follows from the definition of the inverse image.

In many situations, knowing the stalks of a sheaf is enough to control the sheaf itself. For example, whether or not a morphism of sheaves is a monomorphism, epimorphism, or isomorphism can be tested on the stalks. They also find use in constructions such as Godement resolutions.

\subsection{Ringed spaces and locally ringed spaces}
A pair $(X, \mathcal{O}_X)$ consisting of a topological space $X$ and a sheaf of rings on $X$ is called a ringed space. Many types of spaces can be defined as certain types of ringed spaces. The sheaf $\mathcal{O}_X$ is called the structure sheaf of the space. A very common situation is when all the stalks of the structure sheaf are local rings, in which case the pair is called a locally ringed space. Here are examples of definitions made in this way:
\begin{itemize}
  \item An $n$-dimensional $C^k$ manifold $M$ is a locally ringed space whose structure sheaf is an $\underline{\mathbf{R}}$-algebra and is locally isomorphic to the sheaf of $C^k$ real-valued functions on $R^n$.
  \item A complex analytic space is a locally ringed space whose structure sheaf is a $\underline{\mathbf{C}}$-algebra and is locally isomorphic to the vanishing locus of a finite set of holomorphic functions together with the restriction (to the vanishing locus) of the sheaf of holomorphic functions on $C^n$ for some $n$.
  \item A scheme is a locally ringed space that is locally isomorphic to the spectrum of a ring.
  \item A semialgebraic space is a locally ringed space that is locally isomorphic to a semialgebraic set in Euclidean space together with its sheaf of semialgebraic functions.

\end{itemize}

\subsection{Sheaves of modules}
Let $(X, \mathcal{O}_X)$ be a ringed space. A sheaf of modules is a sheaf $\mathcal{M}$ such that on every open set $U$ of $X, \mathcal{M}(U)$ is an $\mathcal{O}_X(U)$-module and for every inclusion of open sets $V \subseteq U$, the restriction map $\mathcal{M}(U) \to \mathcal{M}(V)$ is a homomorphism of $\mathcal{O}_X(U)$-modules.

Most important geometric objects are sheaves of modules. For example, there is a one-to-one correspondence between vector bundles and locally free sheaves of $\mathcal{O}_X$-modules. Sheaves of solutions to differential equations are D-modules, that is, modules over the sheaf of differential operators.

A particularly important case are abelian sheaves, which are modules over the constant sheaf $\underline{\mathbf{Z}}$. Every sheaf of modules is an abelian sheaf.
The condition that a module is finitely generated or finitely presented can also be formulated for a sheaf of modules. $\mathcal{M}$ is finitely generated if, for every point $x$ of $X$, there exists an open neighborhood $U$ of $x$, a natural number $n$ (possibly depending on $U$), and a surjective morphism of sheaves $\mathcal{O}_X^n|_U \to \mathcal{M}|_U$. Similarly, $\mathcal{M}$ is finitely presented if in addition there exists a natural number $m$ (again possibly depending on $U$) and a morphism of sheaves $\mathcal{O}_X^m|_U \to \mathcal{O}_X^n|_U$ such that the sequence of morphisms $\mathcal{O}_X^m|_U \to \mathcal{O}_X^n|_U \to \mathcal{M}$ is exact. Equivalently, the kernel of the morphism $\mathcal{O}_X^n|_U \to \mathcal{M}$ is itself a finitely generated sheaf.

These, however, are not the only possible finiteness conditions on a sheaf. The most important finiteness condition for a sheaf is coherence. $\mathcal{M}$ is coherent if it is of finite type and if, for every open set $U$ and every morphism of sheaves $\phi : \mathcal{O}_X^n \to \mathcal{M}$ (not necessarily surjective), the kernel of $f$ is of finite type. $\mathcal{O}_X$ is coherent if it is coherent as a module over itself. Note that coherence is a strictly stronger condition than finite presentation: $\mathcal{O}_X$ is always finitely presented as a module over itself, but it is not always coherent. For example, let $X$ be a point, let $\mathcal{O}_X$ be the ring $R = C[x_1, x_2, ...]$ of complex polynomials in countably many indeterminates. Choose $n = 1$, and for the morphism $f$, take the map that sends every variable to zero. The kernel of this map is not finitely generated, so $\mathcal{O}_X$ is not coherent.

\section{Scheme}

In mathematics, schemes connect the fields of algebraic geometry, commutative algebra and number theory. Schemes were introduced by Alexander Grothendieck in 1960 in his treatise \'{E}l\'{e}ments de g\'{e}om\'{e}trie alg\'{e}brique, with the aim of developing the formalism needed to solve deep problems of algebraic geometry, such as the Weil conjectures (the last of which was proved by Pierre Deligne). Schemes enlarge the notion of algebraic variety to include nilpotent elements (the equations $x = 0$ and $x^2 = 0$ define the same points, but different schemes), and "varieties" defined over any commutative ring.

To be technically precise, a scheme is a topological space together with commutative rings for all of its open sets, which arises from gluing together spectra (spaces of prime ideals) of commutative rings along their open subsets. In other words, it is a locally ringed space which is locally a spectrum of a commutative ring. There are many ways one can qualify a scheme. According to a basic idea of Grothendieck, conditions should be applied to a morphism of schemes. Any scheme S has a unique morphism to $Spec(Z)$, so this attitude, part of the relative point of view, doesn't lose anything. For details on the development of scheme theory, which quickly becomes technically demanding, see first glossary of scheme theory.

\subsection{History and motivation}
The algebraic geometers of the Italian school had often used the somewhat foggy concept of "generic point" when proving statements about algebraic varieties. What is true for the generic point is true for all points of the variety except a small number of special points. In the 1920s, Emmy Noether had first suggested a way to clarify the concept: start with the coordinate ring of the variety (the ring of all polynomial functions defined on the variety); the maximal ideals of this ring will correspond to ordinary points of the variety (under suitable conditions), and the non-maximal prime ideals will correspond to the various generic points, one for each subvariety. By taking all prime ideals, one thus gets the whole collection of ordinary and generic points. Noether did not pursue this approach.

In the 1930s, Wolfgang Krull turned things around and took a radical step: start with any commutative ring, consider the set of its prime ideals, turn it into a topological space by introducing the Zariski topology, and study the algebraic geometry of these quite general objects. Others did not see the point of this generality and Krull abandoned it.

Andr� Weil was especially interested in algebraic geometry over finite fields and other rings. In the 1940s he returned to the prime ideal approach; he needed an abstract variety (outside projective space) for foundational reasons, particularly for the existence in an algebraic setting of the Jacobian variety. In Weil's main foundational book (1946), generic points are constructed by taking points in a very large algebraically closed field, called a universal domain.

In 1944 Oscar Zariski defined an abstract Zariski�Riemann space from the function field of an algebraic variety, for the needs of birational geometry: this is like a direct limit of ordinary varieties (under 'blowing up'), and the construction, reminiscent of locale theory, used valuation rings as points.

In the 1950s, Jean-Pierre Serre, Claude Chevalley and Masayoshi Nagata, motivated largely by the Weil conjectures relating number theory and algebraic geometry, pursued similar approaches with prime ideals as points. According to Pierre Cartier, the word scheme was first used in the 1956 Chevalley Seminar, in which Chevalley was pursuing Zariski's ideas; and it was Andr� Martineau who suggested to Serre the move to the current spectrum of a ring in general.

\subsection{Modern definitions of the objects of algebraic geometry}
Alexander Grothendieck then gave the decisive definition, bringing to a conclusion a generation of experimental suggestions and partial developments.[citation needed] He defined the spectrum of a commutative ring as the space of prime ideals with Zariski topology, but augments it with a sheaf of rings: to every Zariski-open set he assigns a commutative ring, thought of as the ring of "polynomial functions" defined on that set. These objects are the "affine schemes"; a general scheme is then obtained by "gluing together" several such affine schemes, in analogy to the fact that general varieties can be obtained by gluing together affine varieties.

The generality of the scheme concept was initially criticized: some schemes are removed from having straightforward geometrical interpretation, which made the concept difficult to grasp. However, admitting arbitrary schemes makes the whole category of schemes better-behaved. Moreover, natural considerations regarding, for example, moduli spaces, lead to schemes that are "non-classical". The occurrence of these schemes that are not varieties (nor built up simply from varieties) in problems that could be posed in classical terms made for the gradual acceptance of the new foundations of the subject.

Subsequent work on algebraic spaces and algebraic stacks by Deligne, Mumford, and Michael Artin, originally in the context of moduli problems, has further enhanced the geometric flexibility of modern algebraic geometry. Grothendieck advocated certain types of ringed toposes as generalisations of schemes, and following his proposals relative schemes over ringed toposes were developed by M. Hakim. Recent ideas about higher algebraic stacks and homotopical or derived algebraic geometry have regard to further expanding the algebraic reach of geometric intuition, bringing algebraic geometry closer in spirit to homotopy theory.

\subsection{Definitions}

An affine scheme is a locally ringed space isomorphic to the spectrum of a commutative ring. We denote the spectrum of a commutative ring A by $Spec(A)$. A scheme is a locally ringed space $X$ admitting a covering by open sets $U_i$, such that the restriction of the structure sheaf $O_X$ to each $U_i$ is an affine scheme. Therefore one may think of a scheme as being covered by "coordinate charts" of affine schemes. The above formal definition means exactly that schemes are obtained by glueing together affine schemes for the Zariski topology.

In the early days, this was called a prescheme, and a scheme was defined to be a separated prescheme. The term prescheme has fallen out of use, but can still be found in older books, such as Grothendieck's �l�ments de g�om�trie alg�brique and Mumford's Red Book .

\subsection{The category of schemes}

Schemes form a category if we take as morphisms the morphisms of locally ringed spaces. (See also: morphism of schemes.)

Morphisms from schemes to affine schemes are completely understood in terms of ring homomorphisms by the following contravariant adjoint pair: For every scheme $X$ and every commutative ring $A$ we have a natural equivalence

$$\operatorname{Hom}_{\rm Schemes}(X, \operatorname{Spec}(A)) \cong \operatorname{Hom}_{\rm CRing}(A, {\mathcal O}_X(X))$$
Since $Z$ is an initial object in the category of rings, the category of schemes has $Spec(Z)$ as a final object.

The category of schemes has finite products, but one has to be careful: the underlying topological space of the product scheme of $(X, O_X)$ and $(Y, O_Y)$ is normally not equal to the product of the topological spaces $X$ and $Y$. In fact, the underlying topological space of the product scheme often has more points than the product of the underlying topological spaces. For example, if $K$ is the field with nine elements, then $Spec K \times Spec K \approx Spec (K \otimes_Z K) \approx Spec (K \otimes_{Z/3Z} K) \approx Spec (K � K)$, a set with two elements, though $Spec K$ has only a single element.

For a scheme $S$, the category of schemes over $S$ has also fibre products, and since it has a final object $S$, it follows that it has finite limits.

\subsection{$O_X$ modules}
Just as the $R$-modules are central in commutative algebra when studying the commutative ring $R$, so are the $O_X$-modules central in the study of the scheme X with structure sheaf $O_X$. (See locally ringed space for a definition of $O_X$-modules.) The category of $O_X$-modules is abelian. Of particular importance are the coherent sheaves on $X$, which arise from finitely generated (ordinary) modules on the affine parts of $X$. The category of coherent sheaves on $X$ is also abelian.

The sections of the structure sheaf $O_X$ of $X$ are called regular functions, which are defined on each open subsets $U \in X$. The invertible subsheaf of $O_X$, denoted $O^*_X$, consists only of the invertible germs of regular functions under the multiplication. In most situations, the sheaf $K_X$ is defined on an open affine subset Spec $A$ of $X$ as the total quotient rings $Q(A)$ (though there are cases where the definition is more complicated). The sections of $K_X$ are called rational functions on $X$. The invertible subsheaf of $K_X$ is denoted by $K^*_X$. The equivalent class of this invertible sheaf turns to be an abelian group with tensor products and isomorphic to $H_1(X, O^*_X)$, which is called Picard group. On projective varieties the sections of the structure sheaf $O_X$ defined on each open subsets $U$ of $X$ are also called regular functions though there are no global sections except for constants.

\section{Noncommutative torus}
In mathematics, and more specifically in the theory of $C^*$-algebras, the \textbf{noncommutative tori} $A_\theta$ (also known as irrational rotation algebras for irrational values of $\theta$) are a family of noncommutative $C^*$-algebras which generalize the algebra of continuous functions on the 2-torus. Many topological and geometric properties of the classical 2-torus have algebraic analogues for the noncommutative tori, and as such they are fundamental examples of a noncommutative space in the sense of Alain Connes.
For any irrational number $\theta$, the noncommutative torus $A_\theta$ is the $C^*$-subalgebra of $B(L^2(\mathbb{T}))$, the algebra of bounded linear operators of square-integrable functions on the unit circle of $C$ generated by unitary elements $U$ and $V$, where $U(f)(z)=zf(z)$ and $V(f)(z)=f(e^{-2\pi i\theta z})$. A quick calculation shows that $VU = e^{-2\pi i\theta} UV$
\subsection{Alternative characterizations}
\begin{itemize}
  \item \textbf{Universal property}: $A_\theta$ can be defined (up to isomorphism) as the universal $C^*$-algebra generated by two unitary elements $U$ and $V$ satisfying the relation $VU = e^{2\pi i \theta}UV$. This definition extends to the case when $\theta$ is rational. In particular when $\theta = 0$, $A_\theta$  is isomorphic to continuous functions on the 2-torus by the Gelfand transform.
  \item \textbf{Irrational rotation algebra}: Let the infinite cyclic group $Z$ act on the circle $S^1$ by the rotation action by angle $2\pi i \theta$. This induces an action of $Z$ by automorphisms on the algebra of continuous functions $C(S^1)$. The resulting $C^*$-crossed product $C(S^1) \otimes Z$ is isomorphic to $A_\theta$. The generating unitaries are the generator of the group $Z$ and the identity function on the circle $z : S^1 \to C$
  \item \textbf{Twisted group algebra}: The function $s : Z^2 \times Z^2 \to C$; $s((m,n), (p,q)) = e^{2\pi inp\theta}$ is a group 2-cocycle on $Z^2$, and the corresponding twisted group algebra $C^*(Z^2; s)$ is isomorphic to $A_\theta$.
\end{itemize}

\subsection{Classification and $K$-theory}

The $K$-theory of $A_\theta$ is $Z^2$ in both even dimension and odd dimension, and so does not distinguish the irrational rotation algebras. But as ordered groups, $K_0 \approx Z + \theta Z$. Therefore, two noncommutative tori $A_\theta$ and $A_\eta$ are isomorphic if and only if either $\theta + \eta$ or $\theta - \eta$ is an integer.

Two irrational rotation algebras $A_\theta$ and $A_\eta$ are strongly Morita equivalent if and only if $\theta$ and $ \eta$ are in the same orbit of the action of $SL(2, Z)$ on $R$ by fractional linear transformations. In particular, the noncommutative tori with $\theta$ rational are Morita equivalent to the classical torus. On the other hand, the noncommutative tori with $\theta$ irrational are simple $C^*$-algebras

\section{Grothendieck topology}

In category theory, a branch of mathematics, a Grothendieck topology is a structure on a category C which makes the objects of $C$ act like the open sets of a topological space. A category together with a choice of Grothendieck topology is called a site.

Grothendieck topologies axiomatize the notion of an open cover. Using the notion of covering provided by a Grothendieck topology, it becomes possible to define sheaves on a category and their cohomology. This was first done in algebraic geometry and algebraic number theory by Alexander Grothendieck to define the �tale cohomology of a scheme. It has been used to define other cohomology theories since then, such as $l$-adic cohomology, flat cohomology, and crystalline cohomology. While Grothendieck topologies are most often used to define cohomology theories, they have found other applications as well, such as to John Tate's theory of rigid analytic geometry.

There is a natural way to associate a site to an ordinary topological space, and Grothendieck's theory is loosely regarded as a generalization of classical topology. Under meager point-set hypotheses, namely sobriety, this is completely accurate�it is possible to recover a sober space from its associated site. However simple examples such as the indiscrete topological space show that not all topological spaces can be expressed using Grothendieck topologies. Conversely, there are Grothendieck topologies which do not come from topological spaces.
Andr\'{e} Weil's famous Weil conjectures proposed that certain properties of equations with integral coefficients should be understood as geometric properties of the algebraic variety that they define. His conjectures postulated that there should be a cohomology theory of algebraic varieties which gave number-theoretic information about their defining equations. This cohomology theory was known as the "Weil cohomology", but using the tools he had available, Weil was unable to construct it.
In the early 1960s, Alexander Grothendieck introduced �tale maps into algebraic geometry as algebraic analogues of local analytic isomorphisms in analytic geometry. He used �tale coverings to define an algebraic analogue of the fundamental group of a topological space. Soon Jean-Pierre Serre noticed that some properties of \'{e}tale coverings mimicked those of open immersions, and that consequently it was possible to make constructions which imitated the cohomology functor H1. Grothendieck saw that it would be possible to use Serre's idea to define a cohomology theory which he suspected would be the Weil cohomology. To define this cohomology theory, Grothendieck needed to replace the usual, topological notion of an open covering with one that would use \'{e}tale coverings instead. Grothendieck also saw how to phrase the definition of covering abstractly; this is where the definition of a Grothendieck topology comes from.

\subsection{Definition}

\begin{itemize}
  \item \textbf{Motivation} The classical definition of a sheaf begins with a topological space $X$. A sheaf associates information to the open sets of $X$. This information can be phrased abstractly by letting $O(X)$ be the category whose objects are the open subsets $U$ of $X$ and whose morphisms are the inclusion maps $V \to U$ of open sets $U$ and $V$ of $X$. We will call such maps open immersions, just as in the context of schemes. Then a presheaf on $X$ is a contravariant functor from $O(X)$ to the category of sets, and a sheaf is a presheaf which satisfies the gluing axiom. The gluing axiom is phrased in terms of pointwise covering, i.e., $\{U_i\}$ covers $U$ if and only if  $\cup_i U_i = U$. In this definition, $U_i$ is an open subset of $X$. Grothendieck topologies replace each $U_i$ with an entire family of open subsets; in this example, $U_i$ is replaced by the family of all open immersions $V_{ij} \to U_i$. Such a collection is called a sieve. Pointwise covering is replaced by the notion of a covering family; in the above example, the set of all $\{V_{ij} \to U_i\}_j$ as $i$ varies is a covering family of $U$. Sieves and covering families can be axiomatized, and once this is done open sets and pointwise covering can be replaced by other notions which describe other properties of the space $X$.
  \item \textbf{Sieves} In a Grothendieck topology, the notion of a collection of open subsets of $U$ stable under inclusion is replaced by the notion of a sieve. If $c$ is any given object in $C$, a sieve on $c$ is a subfunctor of the functor $Hom(-, c)$; (this is the Yoneda embedding applied to $c$). In the case of $O(X)$, a sieve $S$ on an open set $U$ selects a collection of open subsets of $U$ which is stable under inclusion. More precisely, consider that for any open subset $V$ of $U$, $S(V)$ will be a subset of $Hom(V, U)$, which has only one element, the open immersion $V \to U$. Then $V$ will be considered "selected" by $S$ if and only if $S(V)$ is nonempty. If $W$ is a subset of $V$ then there is a morphism $S(V) \to S(W)$ given by composition with the inclusion $W \to V$. If $S(V)$ is non-empty, it follows that $S(W)$ is also non-empty.

If $S$ is a sieve on $X$, and $f: Y \to X$ is a morphism, then left composition by $f$ gives a sieve on $Y$ called the pullback of $S$ along $f$, denoted by $f^\ast S$. It is defined as the fibered product $S \times_{Hom(-, X)} Hom(-, Y)$ together with its natural embedding in $Hom(-, Y)$. More concretely, for each object $Z$ of $C$, $f^\ast S(Z) = \{ g: Z \to Y | fg \in S(Z) \}$, and $f^\ast S$ inherits its action on morphisms by being a subfunctor of $Hom(-, Y)$. In the classical example, the pullback of a collection $\{V_i\}$ of subsets of $U$ along an inclusion $W \to U$ is the collection $\{V_i \cap W\}$.
  \item \textbf{Grothendieck topology} A Grothendieck topology $J$ on a category $C$ is a collection, for each object $c$ of $C$, of distinguished sieves on $c$, denoted by $J(c)$ and called covering sieves of $c$. This selection will be subject to certain axioms, stated below. Continuing the previous example, a sieve $S$ on an open set $U$ in $O(X)$ will be a covering sieve if and only if the union of all the open sets $V$ for which $S(V)$ is nonempty equals $U$; in other words, if and only if $S$ gives us a collection of open sets which cover $U$ in the classical sense

     \item \textbf{Axioms} The conditions we impose on a Grothendieck topology are
              \begin{description}
                \item [T1] (Base change) If $S$ is a covering sieve on $X$, and $f: Y \to X$ is a morphism, then the pullback $f^\ast S$ is a covering sieve on $Y$.
                \item[T2] (Local character) Let $S$ be a covering sieve on $X$, and let $T$ be any sieve on $X$. Suppose that for each object $Y$ of $C$ and each arrow $f: Y \to X$ in $S(Y)$, the pullback sieve $f^\ast T$ is a covering sieve on $Y$. Then $T$ is a covering sieve on $X$.
                \item[T3]  (Identity) $Hom(-, X)$ is a covering sieve on $X$ for any object $X$ in $C$.
              \end{description}
        The base change axiom corresponds to the idea that if $\{U_i\}$ covers $U$, then $\{U_i \cap V\}$ should cover $U \cap V$. The local character axiom corresponds to the idea that if $\{U_i\}$ covers $U$ and $\{V_{ij}\}_j \in J_i$ covers $U_i$ for each $i$, then the collection $\{V_{ij}\}$ for all $i$ and $j$ should cover $U$. Lastly, the identity axiom corresponds to the idea that any set is covered by all its possible subsets.
 \item \textbf{Grothendieck pretopologies} In fact, it is possible to put these axioms in another form where their geometric character is more apparent, assuming that the underlying category $C$ contains certain fibered products. In this case, instead of specifying sieves, we can specify that certain collections of maps with a common codomain should cover their codomain. These collections are called covering families. If the collection of all covering families satisfies certain axioms, then we say that they form a Grothendieck pretopology. These axioms are:   \begin{description}
  \item[(PT 0)] (Existence of fibered products) For all objects $X$ of $C$, and for all morphisms $X_0 \to X$ which appear in some covering family of $X$, and for all morphisms $Y \to X$, the fibered product $X_0 �_X Y$ exists.
  \item[(PT 1)] (Stability under base change) For all objects $X$ of $C$, all morphisms $Y \to X$, and all covering families $\{X_\alpha \to X\}$, the family $\{X_\alpha �_X Y \to Y\}$ is a covering family.
  \item[(PT 2)] (Local character) If $\{X_\alpha \to X\}$ is a covering family, and if for all $\alpha, \{X_{\beta \alpha} \to X\alpha \}$ is a covering family, then the family of composites $\{X_{\beta\alpha} \to X_\alpha \to X\}$ is a covering family.
  \item[(PT 3)] (Isomorphisms) If $f: Y \to X$ is an isomorphism, then $\{f\}$ is a covering family.
  For any pretopology, the collection of all sieves that contain a covering family from the pretopology is always a Grothendieck topology.

For categories with fibered products, there is a converse. Given a collection of arrows $\{X_\alpha \to X\}$, we construct a sieve $S$ by letting $S(Y)$ be the set of all morphisms $Y \to X$ that factor through some arrow $X_\alpha \to X$. This is called the sieve generated by $\{X_\alpha \to X\}$. Now choose a topology. Say that $\{X_\alpha \to X\}$ is a covering family if and only if the sieve that it generates is a covering sieve for the given topology. It is easy to check that this defines a pretopology.
(PT 3) is sometimes replaced by a weaker axiom:
  \item[(PT 3')]  (Identity) If $1_X : X \to X$ is the identity arrow, then $\{1_X\}$ is a covering family.
  (PT 3) implies (PT 3'), but not conversely. However, suppose that we have a collection of covering families that satisfies (PT 0) through (PT 2) and (PT 3'), but not (PT 3). These families generate a pretopology. The topology generated by the original collection of covering families is then the same as the topology generated by the pretopology, because the sieve generated by an isomorphism $Y \to X$ is $Hom(-, X)$. Consequently, if we restrict our attention to topologies, (PT 3) and (PT 3') are equivalent.
\end{description}

\end{itemize}

\subsection{Sites and sheaves}

Let $C$ be a category and let $J$ be a Grothendieck topology on $C$. The pair $(C, J)$ is called a site.

A presheaf on a category is a contravariant functor from $C$ to the category of all sets. Note that for this definition $C$ is not required to have a topology. A sheaf on a site, however, should allow gluing, just like sheaves in classical topology. Consequently, we define a sheaf on a site to be a presheaf $F$ such that for all objects $X$ and all covering sieves $S$ on $X$, the natural map $Hom(Hom(-, X), F) \to Hom(S, F)$, induced by the inclusion of $S$ into $Hom(-, X)$, is a bijection. Halfway in between a presheaf and a sheaf is the notion of a separated presheaf, where the natural map above is required to be only an injection, not a bijection, for all sieves $S$. A morphism of presheaves or of sheaves is a natural transformation of functors. The category of all sheaves on $C$ is the topos defined by the site $(C, J)$.

Using the Yoneda lemma, it is possible to show that a presheaf on the category $O(X)$ is a sheaf on the topology defined above if and only if it is a sheaf in the classical sense.

Sheaves on a pretopology have a particularly simple description: For each covering family $\{X_\alpha \to X\}$, the diagram

$$F(X) \rightarrow \prod_{\alpha\in A} F(X_\alpha) {{{} \atop \longrightarrow}\atop{\longrightarrow \atop {}}} \prod_{\alpha,\beta \in A} F(X_\alpha\times_X X_\beta)$$
must be an equalizer. For a separated presheaf, the first arrow need only be injective.

Similarly, one can define presheaves and sheaves of abelian groups, rings, modules, and so on. One can require either that a presheaf $F$ is a contravariant functor to the category of abelian groups (or rings, or modules, etc.), or that F be an abelian group (ring, module, etc.) object in the category of all contravariant functors from $C$ to the category of sets. These two definitions are equivalent.

\subsection{Examples of sites}
\begin{itemize}
  \item \textbf{The discrete and indiscrete topologies} Let $C$ be any category. To define the discrete topology, we declare all sieves to be covering sieves. If $C$ has all fibered products, this is equivalent to declaring all families to be covering families. To define the indiscrete topology, we declare only the sieves of the form $Hom(-, X)$ to be covering sieves. The indiscrete topology is also known as the biggest or chaotic topology, and it is generated by the pretopology which has only isomorphisms for covering families. A sheaf on the indiscrete site is the same thing as a presheaf.
  \item \textbf{The canonical topology} Let $C$ be any category. The Yoneda embedding gives a functor $Hom(-, X)$ for each object $X$ of $C$. The canonical topology is the biggest (finest) topology such that every representable presheaf, i.e. presheaf of the form $Hom(-, X)$, is a sheaf. A covering sieve or covering family for this site is said to be strictly universally epimorphic. A topology which is less fine than the canonical topology, that is, for which every covering sieve is strictly universally epimorphic, is called subcanonical. Subcanonical sites are exactly the sites for which every presheaf of the form $Hom(-, X)$ is a sheaf. Most sites encountered in practice are subcanonical.
  \item \textbf{Small site associated to a topological space} We repeat the example which we began with above. Let $X$ be a topological space. We defined $O(X)$ to be the category whose objects are the open sets of $X$ and whose morphisms are inclusions of open sets. Note that for an open set $U$ and a sieve $S$ on U, the set $S(V)$ contains either zero or one element for every open set $V$. The covering sieves on an object $U$ of $O(X)$ are those sieves $S$ satisfying the following condition:

\textbf{If $W$ is the union of all the sets $V$ such that $S(V)$ is non-empty, then $W = U$.}
This notion of cover matches the usual notion in point-set topology.

This topology can also naturally be expressed as a pretopology. We say that a family of inclusions $\{V_\alpha \subset U\}$ is a covering family if and only if the union $\cup V_\alpha$ equals $U$. This site is called the small site associated to a topological space $X$.
  \item \textbf{Big site associated to a topological space} Let $Spc$ be the category of all topological spaces. Given any family of functions $\{u_\alpha : V_\alpha \to X\}$, we say that it is a surjective family or that the morphisms ua are jointly surjective if $\cup u_\alpha(V_\alpha)$ equals $X$. We define a pretopology on $Spc$ by taking the covering families to be surjective families all of whose members are open immersions. Let $S$ be a sieve on $Spc.$ $S$ is a covering sieve for this topology if and only if:
\begin{itemize}
  \item For all $Y$ and every morphism $f : Y \to X$ in $S(Y)$, there exists a $V$ and a $g : V \to X$ such that $g$ is an open immersion, $g$ is in $S(V)$, and f factors through $g$.
  \item If $W$ is the union of all the sets $f(Y)$, where $f : Y \to X$ is in $S(Y)$, then $W = X$.
 \end{itemize}
Fix a topological space $X$. Consider the comma category $Spc/X$ of topological spaces with a fixed continuous map to $X$. The topology on $Spc$ induces a topology on $Spc/X$. The covering sieves and covering families are almost exactly the same; the only difference is that now all the maps involved commute with the fixed maps to $X$. This is the big site associated to a topological space $X$ . Notice that $Spc$ is the big site associated to the one point space. This site was first considered by Jean Giraud.
  \item \textbf{The big and small sites of a manifold} Let $M$ be a manifold. $M$ has a category of open sets $O(M)$ because it is a topological space, and it gets a topology as in the above example. For two open sets $U$ and $V$ of $M$, the fiber product $U \times_M V$ is the open set $U \cap V$, which is still in $O(M)$. This means that the topology on $O(M)$ is defined by a pretopology, the same pretopology as before.

Let $Mfd$ be the category of all manifolds and continuous maps. (Or smooth manifolds and smooth maps, or real analytic manifolds and analytic maps, etc.) $Mfd$ is a subcategory of $Spc$, and open immersions are continuous (or smooth, or analytic, etc.), so $Mfd$ inherits a topology from $Spc$. This lets us construct the big site of the manifold $M$ as the site $Mfd/M$. We can also define this topology using the same pretopology we used above. Notice that to satisfy $(PT 0)$, we need to check that for any continuous map of manifolds $X \to Y$ and any open subset $U$ of $Y$, the fibered product $U \times_Y X$ is in $Mfd/M$. This is just the statement that the preimage of an open set is open. Notice, however, that not all fibered products exist in $Mfd$ because the preimage of a smooth map at a critical value need not be a manifold.
  \item \textbf{Topologies on the category of schemes}
  The category of schemes, denoted Sch, has a tremendous number of useful topologies. A complete understanding of some questions may require examining a scheme using several different topologies. All of these topologies have associated small and big sites. The big site is formed by taking the entire category of schemes and their morphisms, together with the covering sieves specified by the topology. The small site over a given scheme is formed by only taking the objects and morphisms which are part of a cover of the given scheme.

The most elementary of these is the Zariski topology. Let $X$ be a scheme. $X$ has an underlying topological space, and this topological space determines a Grothendieck topology. The Zariski topology on $Sch$ is generated by the pretopology whose covering families are jointly surjective families of scheme-theoretic open immersions. The covering sieves $S$ for $Zar$ are characterized by the following two properties:
\begin{itemize}
  \item For all $Y$ and every morphism $f : Y \to X$ in $S(Y)$, there exists a $V$ and a $g : V \to X$ such that $g$ is an open immersion, $g$ is in $S(V)$, and f factors through $g$.

  \item If $W$ is the union of all the sets $f(Y)$, where $f : Y \to X$ is in $S(Y)$, then $W = X$.

\end{itemize}
Despite their outward similarities, the topology on $Zar$ is not the restriction of the topology on $Spc!$ This is because there are morphisms of schemes which are topologically open immersions but which are not scheme-theoretic open immersions. For example, let $A$ be a non-reduced ring and let $N$ be its ideal of nilpotents. The quotient map $A \to A/N$ induces a map $Spec(A/N) \to Spec(A)$ which is the identity on underlying topological spaces. To be a scheme-theoretic open immersion it must also induce an isomorphism on structure sheaves, which this map does not do. In fact, this map is a closed immersion.

The \'{e}tale topology is finer than the Zariski topology. It was the first Grothendieck topology to be closely studied. Its covering families are jointly surjective families of \'{e}tale morphisms. It is finer than the Nisnevich topology, but neither finer nor coarser than the $cdh$ and l' topologies.

There are two flat topologies, the $fppf$ topology and the $fpqc$ topology. $fppf$ stands for fid\`{e}lement plate de pr\'{e}sentation finie, and in this topology, a morphism of affine schemes is a covering morphism if it is faithfully flat, of finite presentation, and is quasi-finite. fpqc stands for fid\`{e}lement plate et quasi-compacte, and in this topology, a morphism of affine schemes is a covering morphism if it is faithfully flat. In both categories, a covering family is defined be a family which is a cover on Zariski open subsets. In the fpqc topology, any faithfully flat and quasi-compact morphism is a cover. These topologies are closely related to descent. The fpqc topology is finer than all the topologies mentioned above, and it is very close to the canonical topology.

Grothendieck introduced crystalline cohomology to study the p-torsion part of the cohomology of characteristic p varieties. In the crystalline topology which is the basis of this theory, covering maps are given by infinitesimal thickenings together with divided power structures. The crystalline covers of a fixed scheme form a category with no final object.
\end{itemize}

\subsection{Continuous and cocontinuous functors}
There are two natural types of functors between sites. They are given by functors which are compatible with the topology in a certain sense.

\begin{itemize}
  \item \textbf{Continuous functors} If $(C, J)$ and $(D, K)$ are sites and $u : C \to D$ is a functor, then $u$ is continuous if for every sheaf $F$ on $D$ with respect to the topology $K$, the presheaf $F_u$ is a sheaf with respect to the topology $J$. Continuous functors induce functors between the corresponding topoi by sending a sheaf $F$ to $F_u$. These functors are called pushforwards. If $\tilde{C}$ and $\tilde{D}$ denote the topoi associated to $C$ and $D$, then the pushforward functor is $u_s : \tilde{D} \to \tilde{C}$.

$u_s$ admits a left adjoint $u^s$ called the pullback. $u^s$ need not preserve limits, even finite limits.

In the same way, $u$ sends a sieve on an object $X$ of $C$ to a sieve on the object $uX$ of $D$. A continuous functor sends covering sieves to covering sieves. If $J$ is the topology defined by a pretopology, and if $u$ commutes with fibered products, then $u$ is continuous if and only if it sends covering sieves to covering sieves and if and only if it sends covering families to covering families. In general, it is not sufficient for $u$ to send covering sieves to covering sieves
  \item \textbf{Cocontinuous functors} Again, let $(C, J)$ and $(D, K)$ be sites and $v : C \to D$ be a functor. If $X$ is an object of $C$ and $R$ is a sieve on $vX$, then $R$ can be pulled back to a sieve S as follows: A morphism $f : Z \to X$ is in $S$ if and only if $v(f) : vZ \to vX$ is in $R$. This defines a sieve. $v$ is cocontinuous if and only if for every object $X$ of $C$ and every covering sieve $R$ of $vX$, the pullback $S$ of $R$ is a covering sieve on $X$.

Composition with $v$ sends a presheaf $F$ on $D$ to a presheaf $Fv$ on $C$, but if $v$ is cocontinuous, this need not send sheaves to sheaves. However, this functor on presheaf categories, usually denoted $\hat{v^*}$, admits a right adjoint $\hat{v_*}$. Then $v$ is cocontinuous if and only if $\hat{v_*}$ sends sheaves to sheaves, that is, if and only if it restricts to a functor $v_* : \tilde{C} \to \tilde{D}$. In this case, the composite of $\hat{v^*}$ with the associated sheaf functor is a left adjoint of $v_*$ denoted $v^*$. Furthermore, $v*$ preserves finite limits, so the adjoint functors $v_*$ and $v^*$ determine a geometric morphism of topoi $\tilde{C} \to \tilde{D}$.
  \item \textbf{Morphisms of sites} A continuous functor $u : C \to D$ is a morphism of sites $D \to C$ (not $C \to D$) if $u^s$ preserves finite limits. In this case, $u^s$ and $u_s$ determine a geometric morphism of topoi $\tilde{C} \to \tilde{D}$. The reasoning behind the convention that a continuous functor $C \to D$ is said to determine a morphism of sites in the opposite direction is that this agrees with the intuition coming from the case of topological spaces. A continuous map of topological spaces $X \to Y$ determines a continuous functor $O(Y) \to O(X)$. Since the original map on topological spaces is said to send $X$ to $Y$, the morphism of sites is said to as well.

A particular case of this happens when a continuous functor admits a left adjoint. Suppose that $u : C \to D$ and $v : D \to C$ are functors with $u$ right adjoint to $v$. Then $u$ is continuous if and only if $v$ is cocontinuous, and when this happens, $u_s$ is naturally isomorphic to v* and us is naturally isomorphic to $v^*$. In particular, $u$ is a morphism of sites.
\end{itemize}

\section{Commutative ring}
In ring theory, a branch of abstract algebra, a commutative ring is a ring in which the multiplication operation is commutative. The study of commutative rings is called commutative algebra.

Some specific kinds of commutative rings are given with the following chain of class inclusions:

$$Commutative rings \supset integral domains \supset integrally closed domains \supset unique factorization domains$$
 $$\supset principal ideal domains \supset Euclidean domains \supset fields \supset finite fields$$

\subsection{Definition and first examples}

A ring is a set $R$ equipped with two binary operations, i.e. operations combining any two elements of the ring to a third. They are called addition and multiplication and commonly denoted by $"+"$ and $"�"$; e.g. $a + b$ and $a � b$. To form a ring these two operations have to satisfy a number of properties: the ring has to be an \textbf{abelian group} under addition as well as a \textbf{monoid} under multiplication, where multiplication distributes over addition; i.e., $a � (b + c) = (a � b) + (a � c)$. The identity elements for addition and multiplication are denoted $0$ and $1$, respectively.

If the multiplication is commutative, i.e.

$$a � b = b � a$$,
then the ring $R$ is called commutative. In the remainder of this article, all rings will be commutative, unless explicitly stated otherwise.

\subsection{Examples}
\begin{itemize}
  \item An important example, and in some sense crucial, is the ring of integers $Z$ with the two operations of addition and multiplication. As the multiplication of integers is a commutative operation, this is a commutative ring. It is usually denoted $Z$ as an abbreviation of the German word Zahlen (numbers).

  \item A field is a commutative ring where every non-zero element a is invertible; i.e., has a multiplicative inverse $b$ such that $a � b = 1$. Therefore, by definition, any field is a commutative ring. The rational, real and complex numbers form fields.

\item If $R$ is a given commutative ring, then the set of all polynomials in the variable $X$ whose coefficients are in $R$ forms the polynomial ring, denoted $R[X]$. The same holds true for several variables.

\item If $V$ is some topological space, for example a subset of some $R^n$, real- or complex-valued continuous functions on $V$ form a commutative ring. The same is true for differentiable or holomorphic functions, when the two concepts are defined, such as for $V$ a complex manifold.

\end{itemize}

\subsection{Ideals and the spectrum}
In contrast to fields, where every nonzero element is multiplicatively invertible, the theory of rings is more complicated. There are several notions to cope with that situation. First, an element a of ring $R$ is called a unit if it possesses a multiplicative inverse. Another particular type of element is the zero divisors, i.e. a non-zero element a such that there exists a non-zero element b of the ring such that $ab = 0$. If $R$ possesses no zero divisors, it is called an integral domain since it closely resembles the integers in some ways.

Many of the following notions also exist for not necessarily commutative rings, but the definitions and properties are usually more complicated. For example, all ideals in a commutative ring are automatically two-sided, which simplifies the situation considerably.
\subsection{Ideals and factor rings}

The inner structure of a commutative ring is determined by considering its ideals, i.e. nonempty subsets that are closed under multiplication with arbitrary ring elements and addition: for all $r \in R$, $i$ and $j$ in $I$, both $r_i$ and $i + j$ are required to be in $I$. Given any subset $F = \{f_j\}_{j \in J}$ of $R$ (where $J$ is some index set), the ideal generated by $F$ is the smallest ideal that contains $F$. Equivalently, it is given by finite linear combinations

$$r_1f_1 + r_2f_2 + ... + r_nf_n$$.
An ideal generated by one element is called a principal ideal. A ring all of whose ideals are principal is called a principal ideal ring; two important cases are $Z$ and $k[X]$, the polynomial ring over a field $k$. Any ring has two ideals, namely the zero ideal $\{0\}$ and $R$, the whole ring. An ideal is proper if it is strictly smaller than the whole ring. An ideal that is not strictly contained in any proper ideal is called maximal. An ideal $m$ is maximal if and only if $R / m$ is a field. Except for the zero ring, any ring (with identity) possesses at least one maximal ideal; this follows from Zorn's lemma.

The definition of ideals is such that "dividing" $I$ "out" gives another ring, the factor ring $R / I$: it is the set of cosets of $I$ together with the operations

$$(a + I) + (b + I) = (a + b) + I$$ and $$(a + I)(b + I) = ab + I$$
For example, the ring $Z/nZ$ (also denoted $Z_n$), where $n$ is an integer, is the ring of integers modulo $n$. It is the basis of modular arithmetic.

\subsection{Localizations}

The localization of a ring is the counterpart to factor rings insofar as in a factor ring $R / I$ certain elements (namely the elements of $I$) become zero, whereas in the localization certain elements are rendered invertible, i.e. multiplicative inverses are added to the ring. Concretely, if $S$ is a multiplicatively closed subset of $R$ (i.e. whenever $s, t \in S$ then so is $st$) then the localization of $R$ at $S$, or ring of fractions with denominators in $S$, usually denoted $S^{-1}R$ consists of symbols

$$\frac{r}{s} \mbox{with} r \in R, s \in S $$
subject to certain rules that mimick the cancellation familiar from rational numbers. Indeed, in this language $Q$ is the localization of $Z$ at all nonzero integers. This construction works for any integral domain $R$ instead of $Z$. The localization $(R \setminus \{0\})^{-1}R$ is called the quotient field of $R$. If $S$ consists of the powers of one fixed element $f$, the localisation is written $R_f$.

 \subsection{Prime ideals and the spectrum}

A particularly important type of ideals is prime ideals, often denoted p. This notion arose when algebraists (in the $19^{th}$ century) realized that, unlike in $Z$, in many rings there is no unique factorization into prime numbers. (Rings where it does hold are called unique factorization domains.) By definition, a prime ideal is a proper ideal such that, whenever the product $ab$ of any two ring elements $a$ and $b$ is in $p$, at least one of the two elements is already in $p$. (The opposite conclusion holds for any ideal, by definition). Equivalently, the factor ring $R / p$ is an integral domain. Yet another way of expressing the same is to say that the complement $R \setminus p$ is multiplicatively closed. The localisation $(R \setminus p)^{-1}R$ is important enough to have its own notation: $R_p$. This ring has only one maximal ideal, namely $pR_p$. Such rings are called local.

By the above, any maximal ideal is prime. Proving that an ideal is prime, or equivalently that a ring has no zero-divisors can be very difficult.
Prime ideals are the key step in interpreting a ring geometrically, via the spectrum of a ring Spec $R$: it is the set of all prime ideals of $R$ As noted above, there is at least one prime ideal, therefore the spectrum is nonempty. If $R$ is a field, the only prime ideal is the zero ideal, therefore the spectrum is just one point. The spectrum of $Z$, however, contains one point for the zero ideal, and a point for any prime number $p$ (which generates the prime ideal $pZ$). The spectrum is endowed with a topology called the Zariski topology, which is determined by specifying that subsets $D(f) = \{p \in Spec(R), f \notin p\}$, where $f$ is any ring element, be open. This topology tends to be different from those encountered in analysis or differential geometry; for example, there will generally be points which are not closed. The closure of the point corresponding to the zero ideal $0 \subset Z$, for example, is the whole spectrum of $Z$.

The notion of a spectrum is the common basis of commutative algebra and algebraic geometry. Algebraic geometry proceeds by endowing $Spec(R)$ with a sheaf $\mathcal{O}$ (an entity that collects functions defined locally, i.e. on varying open subsets). The datum of the space and the sheaf is called an affine scheme. Given an affine scheme, the underlying ring $R$ can be recovered as the global sections of $\mathcal{O}$. Moreover, the established one-to-one correspondence between rings and affine schemes is also compatible with ring homomorphisms: any $f: R \to S$ gives rise to a continuous map in the opposite direction

$Spec(S) \to Spec(R), q \to f^{-1}(q)$, i.e. any prime ideal of $S$ is mapped to its preimage under $f$, which is a prime ideal of $R$.
The spectrum also makes precise the intuition that localisation and factor rings are complementary: the natural maps $R \to R_f$ and $R \to R / fR$ correspond, after endowing the spectra of the rings in question with their Zariski topology, to complementary open and closed immersions respectively.

Altogether the equivalence of the two said categories is very apt to reflect algebraic properties of rings in a geometrical manner. Affine schemes are�much the same way as manifolds are locally given by open subsets of Rn�local models for schemes, which are the object of study in algebraic geometry. Therefore, many notions that apply to rings and homomorphisms stem from geometric intuition.

\subsection{Ring homomorphisms}
As usual in algebra, a function f between two objects that respects the structures of the objects in question is called homomorphism. In the case of rings, a ring homomorphism is a map $f: R \to S$ such that

$$f(a + b) = f(a) + f(b), f(ab) = f(a)f(b)$$ and $$f(1) = 1$$
These conditions ensure $f(0) = 0$, but the requirement that the multiplicative identity element $1$ is preserved under $f$ would not follow from the two remaining properties. In such a situation $S$ is also called an $R$-algebra, by understanding that $s \in S$ may be multiplied by some $r$ of $R$, by setting

$$r � s := f(r) � s$$
The kernel and image of $f$ are defined by $ker (f) = \{r \in R, f(r) = 0\}$ and $im (f) = f(R) = \{f(r), r \in R\}$. The kernel is an ideal of $R$, and the image is a subring of $S$.

\subsection{Modules}
The outer structure of a commutative ring is determined by considering linear algebra over that ring, i.e., by investigating the theory of its modules, which are similar to vector spaces, except that the base is not necessarily a field, but can be any ring $R$. The theory of $R$-modules is significantly more difficult than linear algebra of vector spaces. Module theory has to grapple with difficulties such as modules not having bases, that the rank of a free module (i.e. the analog of the dimension of vector spaces) may not be well-defined and that submodules of finitely generated modules need not be finitely generated (unless $R$ is Noetherian, see below).

Ideals within a ring $R$ can be characterized as $R$-modules which are submodules of $R$. On the one hand, a good understanding of $R$-modules necessitates enough information about $R$. Vice versa, however, many techniques in commutative algebra that study the structure of $R$, by examining its ideals, proceed by studying modules in general.
\subsection{Noetherian rings}
A ring is called Noetherian (in honor of Emmy Noether, who developed this concept) if every ascending chain of ideals

$$0 \subseteq I_0 \subseteq I_1 ... \subseteq In \subseteq I_{n + 1} \subseteq ...$$
becomes stationary, i.e. becomes constant beyond some index $n$. Equivalently, any ideal is generated by finitely many elements, or, yet equivalent, submodules of finitely generated modules are finitely generated. A ring is called Artinian (after Emil Artin), if every descending chain of ideals

$$R \supseteq I_0 \supseteq I_1 ... \supseteq I_n \supseteq I_{n + 1} \supseteq ...$$
becomes stationary eventually. Despite the two conditions appearing symmetric, Noetherian rings are much more general than Artinian rings. For example, $Z$ is Noetherian, since every ideal can be generated by one element, but is not Artinian, as the chain

$$Z \supseteq 2Z \supseteq 4Z \supseteq 8Z \supseteq ...$$
shows. In fact, by the Hopkins�Levitzki theorem, every Artinian ring is Noetherian.

Being Noetherian is an extremely important finiteness condition. The condition is preserved under many operations that occur frequently in geometry: if $R$ is Noetherian, then so is the polynomial ring $R[X_1, X_2, ..., X_n]$ (by Hilbert's basis theorem), any localization $S^{-1}R$, factor rings $R / I$.

\section{Spectrum of a ring}
In abstract algebra and algebraic geometry, the spectrum of a commutative ring $R$, denoted by $Spec(R)$, is the set of all prime ideals of $R$. It is commonly augmented with the Zariski topology and with a structure sheaf, turning it into a locally ringed space.
\subsection{Zariski topology}

For any ideal $I$ of $R$, define $V_I$ to be the set of prime ideals containing $I$. We can put a topology on $Spec(R)$ by defining the collection of closed sets to be

$$\{ V_I \colon I \text{ is an ideal of } R \}$$
This topology is called the Zariski topology.

A basis for the Zariski topology can be constructed as follows. For $f\in R$, define $D_f$ to be the set of prime ideals of $R$ not containing $f$. Then each $D_f$ is an open subset of $Spec(R)$, and $\{D_f:f\in R\}$ is a basis for the Zariski topology.

$Spec(R)$ is a compact space, but almost never Hausdorff: in fact, the maximal ideals in $R$ are precisely the closed points in this topology. However, $Spec(R)$ is always a Kolmogorov space. It is also a spectral space
\subsection{Sheaves and schemes}

Given the space $X=Spec(R)$ with the Zariski topology, the structure sheaf $O_X$ is defined on the $D_f$ by setting $\Gamma(D_f, O_X) = R_f$, the localization of $R$ at the multiplicative system $\{1,f,f^2,f^3,...\}$. It can be shown that this satisfies the necessary axioms to be a $B$-Sheaf. Next, if $U$ is the union of $\{D_{f_i}\}_{i\in I}$, we let $\Gamma(U,O_X) = \lim_{i\in I} R_{f_i}$, and this produces a sheaf; see the Gluing axiom article for more detail.

If $R$ is an integral domain, with field of fractions $K$, then we can describe the ring $\Gamma(U,O_X)$ more concretely as follows. We say that an element $f$ in $K$ is regular at a point $P$ in $X$ if it can be represented as a fraction $f = a/b$ with $b$ not in $P$. Note that this agrees with the notion of a regular function in algebraic geometry. Using this definition, we can describe $\Gamma(U,O_X)$ as precisely the set of elements of $K$ which are regular at every point $P$ in $U$.

If $P$ is a point in $Spec(R)$, that is, a prime ideal, then the stalk at $P$ equals the localization of $R$ at $P$, and this is a local ring. Consequently, $Spec(R)$ is a locally ringed space.

Every locally ringed space isomorphic to one of this form is called an affine scheme. General schemes are obtained by "gluing together" several affine schemes.

\subsection{Functoriality}

It is useful to use the language of category theory and observe that Spec is a functor. Every ring homomorphism $f : R \to S$ induces a continuous map $Spec(f) : Spec(S) \to Spec(R)$ (since the preimage of any prime ideal in $S$ is a prime ideal in $R$). In this way, Spec can be seen as a contravariant functor from the category of commutative rings to the category of topological spaces. Moreover for every prime $P$ the homomorphism $f$ descends to homomorphisms

$$O_{f^{-1}}(P) \to O_P$$
of local rings. Thus Spec even defines a contravariant functor from the category of commutative rings to the category of locally ringed spaces. In fact it is the universal such functor and this can be used to define the functor Spec up to natural isomorphism.

The functor Spec yields a contravariant equivalence between the category of commutative rings and the category of affine schemes; each of these categories is often thought of as the opposite category of the other.

\subsection{Motivation from algebraic geometry}
Following on from the example, in algebraic geometry one studies algebraic sets, i.e. subsets of $K^n$ (where $K$ is an algebraically closed field) that are defined as the common zeros of a set of polynomials in $n$ variables. If $A$ is such an algebraic set, one considers the commutative ring $R$ of all polynomial functions $A \to K$. The maximal ideals of $R$ correspond to the points of $A$ (because $K$ is algebraically closed), and the prime ideals of $R$ correspond to the subvarieties of $A$ (an algebraic set is called irreducible or a variety if it cannot be written as the union of two proper algebraic subsets).

The spectrum of $R$ therefore consists of the points of $A$ together with elements for all subvarieties of $A$. The points of $A$ are closed in the spectrum, while the elements corresponding to subvarieties have a closure consisting of all their points and subvarieties. If one only considers the points of $A$, i.e. the maximal ideals in $R$, then the Zariski topology defined above coincides with the Zariski topology defined on algebraic sets (which has precisely the algebraic subsets as closed sets).

One can thus view the topological space $Spec(R)$ as an "enrichment" of the topological space $A$ (with Zariski topology): for every subvariety of $A$, one additional non-closed point has been introduced, and this point "keeps track" of the corresponding subvariety. One thinks of this point as the generic point for the subvariety. Furthermore, the sheaf on $Spec(R)$ and the sheaf of polynomial functions on $A$ are essentially identical. By studying spectra of polynomial rings instead of algebraic sets with Zariski topology, one can generalize the concepts of algebraic geometry to non-algebraically closed fields and beyond, eventually arriving at the language of schemes.

\subsection{Global Spec}
There is a relative version of the functor Spec called global Spec, or relative Spec, and denoted by Spec. For a scheme $Y$, and a quasi-coherent sheaf of $O_Y$-algebras $A$, there is a unique scheme $Spec(A)$, and a morphism $f \colon \bold{Spec} \ A \to Y$ such that for every open affine $U \subseteq Y$, there is an isomorphism induced by $f: f^{-1}(U) \cong \mathrm{Spec} \  A(U)$, and such that for open affines $U \subseteq V$, the inclusion $f^{-1}(U) \to f^{-1}(V)$ induces the restriction map $A(V) \to A(U)$. That is, as ring homomorphisms induce opposite maps of spectra, the restriction maps of a sheaf of algebras induce the inclusion maps of the spectra that make up the Spec of the sheaf.

\subsection{Representation theory perspective}
From the perspective of representation theory, a prime ideal $I$ corresponds to a module $ R/I$, and the spectrum of a ring corresponds to irreducible cyclic representations of $R$, while more general subvarieties correspond to possibly reducible representations that need not be cyclic. Recall that abstractly, the representation theory of a group is the study of modules over its group algebra.

The connection to representation theory is clearer if one considers the polynomial ring $R=K[x_1,\dots,x_n]$ or, without a basis, $R=K[V]$. As the latter formulation makes clear, a polynomial ring is the group algebra over a vector space, and writing in terms of $x_i$ corresponds to choosing a basis for the vector space. Then an ideal $I$, or equivalently a module $R/I$, is a cyclic representation of $R$ (cyclic meaning generated by 1 element as an $R$-module; this generalizes 1-dimensional representations).

In the case that the field is algebraically closed (say, the complex numbers), every maximal ideal corresponds to a point in $n$-space, by the nullstellensatz (the maximal ideal generated by $(x_1-a_1), (x_2-a_2),\ldots,(x_n-a_n)$ corresponds to the point $(a_1,\ldots,a_n))$. These representations of $K[V]$ are then parametrized by the dual space $V^*$, the covector being given by sending each $x_i$ to the corresponding $a_i$. Thus a representation of $K^n$ ($K$-linear maps $K^n \to K$) is given by a set of $n$ numbers, or equivalently a covector $K^n \to K$.

Thus, points in n-space, thought of as the max spec of $R=K[x_1,\dots,x_n]$, correspond precisely to 1-dimensional representations of $R$, while finite sets of points correspond to finite-dimensional representations (which are reducible, corresponding geometrically to being a union, and algebraically to not being a prime ideal). The non-maximal ideals then correspond to infinite-dimensional representations.

\subsection{Functional analysis perspective}
The term "spectrum" comes from the use in operator theory. Given a linear operator $T$ on a finite-dimensional vector space $V$, one can consider the vector space with operator as a module over the polynomial ring in one variable $R=K[T]$, as in the structure theorem for finitely generated modules over a principal ideal domain. Then the spectrum of $K[T]$ (as a ring) equals the spectrum of $T$ (as an operator).

Further, the geometric structure of the spectrum of the ring (equivalently, the algebraic structure of the module) captures the behavior of the spectrum of the operator, such as algebraic multiplicity and geometric multiplicity. For instance, for the 2�2 identity matrix has corresponding module:

$$K[T]/(T-1) \oplus K[T]/(T-1)$$
the 2�2 zero matrix has module

$$K[T]/(T-0) \oplus K[T]/(T-0)$$
showing geometric multiplicity 2 for the zero eigenvalue, while a non-trivial 2�2 nilpotent matrix has module

$$K[T]/T^2$$
showing algebraic multiplicity 2 but geometric multiplicity 1.

In more detail:
\begin{itemize}
  \item the eigenvalues (with geometric multiplicity) of the operator correspond to the (reduced) points of the variety, with multiplicity;

  \item the primary decomposition of the module corresponds to the unreduced points of the variety;

  \item a diagonalizable (semisimple) operator corresponds to a reduced variety;

  \item a cyclic module (one generator) corresponds to the operator having a cyclic vector (a vector whose orbit under $T$ spans the space);

  \item the last invariant factor of the module equals the minimal polynomial of the operator, and the product of the invariant factors equals the characteristic polynomial.

\end{itemize}

\section{Topological Quantum Field Theory}\label{tqft}
A topological quantum field theory (or topological field theory or TQFT) is a
quantum field theory which computes topological invariants.
In physics, topological quantum field theories are the low energy
effective theories of topologically ordered states, such as fractional quantum Hall states,
string-net condensed states, and other strongly correlated quantum liquid states.
In 1988 Atiyah gave a description of topological QFT with axioms.
The basic idea is that a
TQFT is a functor from a certain category of cobordisms to the category of vector spaces.
\begin{defn}
In dimension $d$, TQFT is a monoidal functor $ Z:Cob(d+1) \longrightarrow Vect$, where $Cob(d+1)$ is the
category whose objects are closed, oriented $d$-manifolds $M$ without boundary. The cobordism morphism
$W:M \longrightarrow M'$ is a smooth, oriented, compact $d+1$-dimensional manifold with boundary $\partial W=M \sqcup -M'$.
Two cobordisms $W_1$ and $W_2$ are equivalent if there is an orientation-preserving diffeomorphism $f:W_1 \longrightarrow W_2$.
$Vect$ is a symmetric monoidal category of finite dimensional complex vector space where morphisms are linear maps $L:V_1 \longrightarrow V_2$
with dual $L^*:V_2 \longrightarrow V_1$.
\end{defn}
 This functor satisfies the following axioms:
\begin{enumerate}
\item For a cobordism $W$ with boundary $\partial W=M_1 \sqcup -M_2$, then $Z(W)=Z(M_1)\longrightarrow Z(M_2)$ is a homomorphism, \ie a linear map of vector spaces.
\item  Z is {\em involutory}, that is, $Z(-M)=Z(M)^*$, where $-M$ denotes $M$ with the opposite
orientation and $Z(M)^*=\Hom(Z(M),\C)$ is the dual vector space.
\item Z is {\em multiplicative}, that is, $Z(M_1\sqcup M_2)=Z(M_1) \otimes Z(M_2)$.
\item Z is {\em Associative}: for composite cobordisms (gluing) $W=W_1 \cup_{M_2} W_2$ with $\partial W_1=M_1 \sqcup -M_2$ and
$\partial W_2=M_2 \sqcup -M_3$ then $ Z(W)=Z(W_1)\circ Z(W_2) \in Hom(Z(M_1),Z(M_3))$.
\item $Z(\emptyset)=\C$.
\item $Z(M \times I)$ is the identity endomorphism of $Z(M)$.
\end{enumerate}
\begin{rems}
\begin{itemize}
\item The identity endomorphism of $Z(M)$ in $(6)$ and the functoriality  of $Z$ imply
homotopy invariance.
\item Let $W$ be a closed (d+1)-dimensional manifold (with empty boundary). Then by $(5)$ the vector $Z(W)$ is just a complex number. This means that  a TQFT assigns a numerical invariant to each closed (d+1)-manifolds.
\item Let $ W=W_1 \bigsqcup_{M_2} W_2$ as in $(3)$ with $M_1=M_3=\emptyset$. Then one
can cut $W$ along a $d$-manifold $M_2$ and one obtains
$$Z(W)=\langle Z(W_1),Z(W_2)\rangle$$
where $\langle\cdot,\cdot\rangle$ denotes the pairing of the vector space $Z(M_2)$ with its
dual $Z(M_2)^*=Z(-M_2)$.
\end{itemize}
\end{rems}
Here one can give a physical explanation of the meaning to these axioms. In dimension $3$ we can suppose that $M$ is a physical space with an imaginary time $M \times I$ and a Hilbert space $Z(M)$
of the quantum theory associated to the Hamiltonian $H$
with evolution operator $e^{itH}$ (where $t$ is the coordinate on the interval $I$).
In axiom $(6)$ the Hamiltonian $H$ vanishes. Thus having a topological QFT implies that there is no real dynamics taking place along the cylinder $M \times I$. Notice that, for a manifold $W$ with $\partial W=\bar{M_1} \cup M_2$, there can still be an interesting propagation
from $M_1$ to $M_2$ and this reflects the nontrivial topology of $W$.

Topological quantum field theories had many important applications in modern geometry, among
these the work of Gromov  on pseudo-holomorphic curves in symplectic geometry.
In TQFT
a vector $Z(W)$ in the Hilbert space $Z(M)$ is called a {\em vacuum state} if
$\partial W=M$ and for a closed manifold $W$ the number $Z(W)$ is the {\em vacuum-vacuum} expectation value. In analogy with the statical mechanics it is also called the partition function.

\section{2-Category}
In category theory, a 2-category is a small category $ \mathcal{C}_2$ with ``morphisms between morphisms". 2-categories are the first case of higher order categories and they are constructed
as follows:
\begin{itemize}
\item $ \mathcal{C}_2$ is defined as a small category enriched over $Cat$ which is  defined as a category of small categories and functors.
Here we mean by enriched category a category whose $Hom-Sets$ are replaced by objects
from some other category. More precisely, a 2-category consists of the following data.
\item A class of objects $(A, B, ....) \in Cat $ called $0$-cells.
\item For all $0$-cells $A$ and $B$, we can define a set $\mathcal{C}_2(A,B)$ which is defined as
a $Hom_{\mathcal{C}_2}(A,B)$ of objects $f: A \longrightarrow B$ which are called
$1$-cells.
\item  A morphism $\alpha: f_1 \longrightarrow f_2$ for any two morphisms $f_1$ and $f_2$
of $\mathcal{C}_2$. These 2-morphisms are called 2-cells.
\item The 2-categorical compositions of 2-morphisms is denoted as  $\bullet$ and is called {\em vertical composition}.
\item For all objects $A$, $ B$ and $ C$, there is a functor
$$\circ:\mathcal{C}_2(A,B) \times \mathcal{C}_2(B,C)\longrightarrow \mathcal{C}_2(A,C) $$
called {\em horizontal composition},
which is associative and admits the identity 2-cells $I_A$ as identities.
\item For any object $A$ there is a functor from the terminal category (with one object and
one arrow) to $\mathcal{C}_2(A,A)$.
\end{itemize}
The notion of 2-category differs from the more general notion of a bicategory in that composition of 1-morphisms is required to be strictly associative, whereas in a bicategory it needs only be associative
up to a $2$-isomorphism.

There are three different ways to obtain a category from a 2-category, all of which we
use in Chapter 1.  They are summarized as follows.

\begin{itemize}
\item {\em Forgetting 2-morphisms}: one is left with the category consisting of the objects
and 1-morphisms of the 2-category.
\item {\em Forgetting objects}: one obtains a category whose objects are the 1-morphisms
of the 2-category and whose morphisms are the 2-morphisms of the 2-category.
\item {\em Equivalence relation}: one uses the 2-morphisms to define an equivalence
relation on the set of 1-morphisms and obtains in this way a category whose objects
are the same as the objects of the 2-category and whose morphisms are the
equivalence classes of 1-morphisms of the two category modulo the equivalence
relation generated by the 2-morphisms.
\end{itemize}

\section{Group Rings}\label{gring}
A group ring is a ring $R[G]$ constructed from a ring $R$ and a group $G$.
As an $R$-module, the ring $R[G]$ is the free module over $R$ generated by the elements
of $G$, that is, the elements of the group ring are finite linear combinations of elements of $G$
with coefficients in $R$,
$$R[G]=\{\sum_{g\in G} \alpha_g g|\alpha_g \in R \}  $$
with all but finitely many of the $\alpha_g$ being $0$.

The $R$-module $R[G]$ is a ring with addition of formal linear combinations
\begin{equation}\label{RGadd}
(\sum_{g \in G, a_g\in R} a_gg)+(\sum_{g \in G,b_g \in R} b_gg)=\sum_{g \in G} (a_g+b_g)g
\end{equation}
and multiplication defined by the group operation in $G$ extended by
linearity and distributivity, and the requirement that elements of $R$ commute with elements of $G$,
\begin{equation}\label{RGmultiply}
(\sum_{g \in G, a_g \in R} a_gg) (\sum_{h \in G, b_h \in R} b_hh)=\sum_{g,h \in G} (a_g b_h)gh.
\end{equation}
If $R$ has a unit element, then this is the unique bilinear multiplication for which
$(1 g)(1 h) = (1 gh)$. In this case, $G$ can be identified with the elements $1g$ of $R[G]$.
The identity element of $G$ is the multiplicative unit in the ring $R[G]$.
If $R$ is commutative, then $R[G]$ is an associative algebra over $R$.
If $R=F$ is a field, then $F[G]$ is an algebra, called the group algebra.

We have the following equivalent descriptions of the group ring
\begin{defn}\cite{pass}
Let $G$ be a group and $R$ a ring. Define the set $R[G]$ to be one of the following equivalent statements:
\begin{itemize}
\item The set of all formal $R$-linear combinations of elements of $G$.
\item The set of all functions $f: G \longrightarrow R$ with $f(g) = 0$ for all but finitely many $g \in G$.
\item The free $R$-module with basis $G$.
\end{itemize}
The ring structure is given as above by \eqref{RGadd} and \eqref{RGmultiply}.
\end{defn}

If $R$ and $G$ are both commutative,\ie $R$ is a commutative ring and $G$ is an abelian group,
then $R[G]$ is commutative.
If $H$ is a subgroup of $G$, then $R[H]$ is a subring of $R[G]$. Similarly, if $S$ is a subring of $R$,
then $S[G]$ is a subring of $R[G]$.

\subsection{Group algebra over a finite group}

We recall briefly the example of group algebras for finite groups. These occur naturally
in the theory of group representations of finite groups. As we have seen above,
when $R$ is a field $F$ the group algebra $F[G]$ is a vector space over $F$,
with a canonical basis $e_g$ given by the elements $g\in G$ and with elements given
by formal sums
$$v=\sum_{g\in G}x_g e_g$$
As we saw in general for group rings, the algebra structure is defined by the multiplication
in the group,
$$e_g.e_h=e_{gh} $$
Thinking of the free vector space as $F$-valued functions on $G$, the algebra multiplication
can be written equivalently as convolution of functions.

The group algebra is an algebra over itself; under the correspondence of representations
over $R$ and $R[G]$ modules, it is the regular representation of the group.
Written as a representation, it is the representation $g\mapsto \rho_g$ with the
action given by $\rho(g).e_h=e{gh}$, or
$$\rho(g).r=\sum_{h\in G}k_h\rho(g).e_h= \sum_{h\in G}k_h.e_{gh} $$

For a finite group, the dimension of the vector space $F[G]$ is equal to the number of elements in the group. The field $F$ is commonly taken to be the complex numbers $\C$ or the reals $\R$.
The group algebra $\C[G]$ of a finite group over the complex numbers is a semisimple ring. This result, Maschke's theorem, allows us to understand $C[G]$ as a finite product of matrix rings with
entries in $\C$.

\subsection{Groupoids, semigroups, semigroupoids and their rings}

In the first chapter of our work we introduced algebras that are generalizations of group rings.
They are generalizations in two different senses. First of all one can pass from groups to
groupoids and define the groupoid ring $R[\cG]$. In a different direction one has generalizations
where one replaces the group by a semigroup and have the corresponding semigroup ring $R[S]$.
In our case, we work with a generalization of both of these concepts which is a semigroupoid $\cS$
and the corresponding ring $R[\cS]$. We recall here these different notions and stress the
way in which they differ from one another and from the original notion of group ring recalled
above.

\subsection{Groupoid Ring}
A groupoid $\mathfrak{G}$ is a small category in which each morphism is an isomorphism.
Thus $\mathfrak{G}$ has a set of morphisms, which we call elements of $\mathfrak{G}$, a set
$Y$ of objects together with range (target) and source functions $r, s: \mathfrak{G}
\longrightarrow Y$ such that, for $g_1,g_2 \in \mathfrak{G}$ with $r (g_1)=s (g_2)$, then
the product or composite $g_2g_1=g_2\circ g_1$ exists, with
$s (g_2g_1)=s(g_1)$ and $r(g_2g_1)=r(g_2)$.
The composition is associative. For $\gamma:Y \longrightarrow \mathfrak{G}$ and
for an element $ x\in Y$ the element $\gamma(x)$ is denoted by $1_x$ and it acts as the identity, and
each element $g$ has an inverse $g^{-1}$ such that $s(g^{-1})=r (g)$ and $r(g^{-1})=s(g)$, with
$g^{-1}g =\gamma(s( g))$ and $gg^{-1} =\gamma(r(g))$.

In a groupoid $\mathfrak{G}$, for $y_1,y_2 \in Y$ we define the set $\mathfrak{G}(y_1,y_2)$
of all morphisms with initial point $y_1$ and final point $y_2$. We say that $\mathfrak{G}$
is transitive if, for all $y_1,y_2 \in Y$  the set $\mathfrak{G}(y_1,y_2)$ is non-empty.
For $y \in Y$ we denote the set $\{g \in \mathfrak{G}: s(g)=y\} $ by $\mathfrak{G}_y$.
Let $\mathfrak{G}$ be a groupoid. The transitive component of $x\in Y$, denoted by
$C(\mathfrak{G})_y$, is the full subgroupoid of $\mathfrak{G}$ on those objects $x \in Y $ such that $\mathfrak{G}(y,x)$ is non-empty.

\begin{defn}
Let $\mathfrak{G}$ be a groupoid and $R$ a ring or a field. The groupoid ring (or groupoid
algebra in the field case) $R[\mathfrak{G}]$
consists of all finite formal sums of the form $\sum_i^n r_ig_i$ where $r_i \in R$ and
$g_i \in \mathfrak{G}$, which satisfy the following conditions.
\begin{enumerate}
\item If $ \sum_{i=1}^n r_ig_i=\sum_{i=1}^n s_ig_i$ then $r_i=s_i$, for $i=1,2,...,n$.\\
\item $(\sum_{i=1}^n r_ig_i)+(\sum_{i=1}^n s_ig_i)=\sum_{i=1}^n (r_i+s_i)g_i$.\\
\item $(\sum_{i=1}^n r_ig_i)(\sum_{i=1}^n s_ig_i)=(\sum_{i=1}^n k_it_i) $ where
$ k_i=\sum_{i,j}^n r_is_j$ and  $t_i=g_ig_j$.\\
\item $r_ig_i=g_ir_i$ for all $r_i \in R$ and $g_i \in \mathfrak{G}$.\\
\item $ r\sum_{i=1}^n r_ig_i=\sum_{i=1}^n rr_ig_i$, for $r,r_i \in R$.\\
\end{enumerate}
\end{defn}

Notice that since $1 \in R$ and $g_i \in \mathfrak{G}$, we have $\mathfrak{G}=1.\mathfrak{G }\subset R[\mathfrak{G}]$
and $R \subset \mathfrak{G}$ if and only if $\mathfrak{G}$ has identity, otherwise $R$ not in $\mathfrak{G}$.

\subsection{Semigroup Ring}

We now similarly recall the notion of semigroup ring, which is another
generalization of the concept of group ring recalled in \S \ref{gring} above.

The construction of the semigroup ring is not far from what we said before for
the group ring. We try to illustrate the concept of semigroup from another
perspective. Let $\cS$ be a semigroup and let  $R$ be a commutative ring.
We define the semigroup ring $R[\cS]$
to be the set of functions $f:\cS \longrightarrow R$ that send all but
finitely elements of $\cS$ to zero,
$$ f(s)=\sum_{m \in \cS} a_m \delta_m(s), $$
where $a_m\in R$ and $\delta_m(s)=\delta_{m,s}$ is the Kronecker delta
function, and all but finitely many of the coefficients are $a_m=0$.
Clearly the set of such functions has the structure of an $R$-module if
is a ring, or a vector space if $R$ is a field. From the product in the
semigroup $\cS$ we can also define a product on the semigroup
ring $R[\cS]$ as follows. Let  $(x,y)$ be a pair of elements of $\cS$ with $x y=s \in
\cS$ then we set
\begin{equation}\label{semgrringprod}
(f*g)(s)=\sum_{x y=s}f(x)g(y)  .
\end{equation}
This is analogous to the way one defines the product in the group ring.
In fact, it takes the product of all non-zero components of $f$ and $g$ and
collects the resulting terms whose indices multiply to the same element of the
semigroup. With this additive and multiplicative structure, one can check that,
as in the case of groups, the set $R[\cS]$ is in fact a ring (or an algebra if $R$
is a field).

In this text we have assumed the convention that semigroups have
a unit. However, the definition above makes sense also for the case
where one does not require $\cS$ to have a unit. In some text the
semigroup ring of a semigroup with unit is called a {\em monoid
ring}. It is then a unital ring with a unit given by the identity (unit) of
the semigroup. If $\cS$ is a group we recover the same definition of
{\em group ring} discussed in \S \ref{gring} above.
If in any of these cases we start with a commutative
semigroup we get a commutative ring.

Notice that if $\cS$ is a group, for the group ring $R[\cS]$, since $xy=s$
only if $y=x^{-1}s$, we can rewrite the product formula
\eqref{semgrringprod} in the equivalent form
$$ (f*g)(s)=\sum_{x \in \cS}f(x)g(x^{-1}s).  $$
This way of multiplying two functions on a group is called
{\em convolution} product.

The notion of semigroupoid and semigroupoid ring is described
in detail in Chapter 1. It is similar to the groupoid case, in as the
compositions are only defined when the range of the first element
agrees with the source of the second, and it is also similar to
the semigroup case, in the sense that not all elements have an
inverse. The notion of semigroupoid ring or algebra that we
consider there is still a natural generalization of the notion
of groupoid ring, as the ones we recalled in this appendix.
%----------------------------------------------------------------------------------------------------------------------

\section{Creation and annihilation operators}
The unitaries $U_k: f \mapsto (U_k f)(n)=f(n+k)$, for $k\in \Z$,
acting on $\ell^2(\Z)$, induce isometries
\begin{equation}\label{Squot}
(S_n f)(k)=\left\{ \begin{array}{ll} f(k+n) & k+n\geq 0 \\[2mm]
0 & k+n <0, \end{array}\right.
\end{equation}
acting on the Hilbert space $\ell^2(\N\cup\{0\})=\ell^2(\Z/V)$ with
$V=\{\pm 1\}$.
\begin{lem}\label{lemSstarS}
The operators $S_n$ of \eqref{Squot}, for $n\in \Z$,
satisfy the relations $S_n^*=S_{-n}$ and
\begin{equation}\label{SstarS}
S^*_n S_n =P_n \ \ \ \text{ and } \ \ \  S_n S_n^*= P_{-n}
\end{equation}
with $P_n$ the projection $P_n f(k)=f(k)\chi_{[n,\infty)}(k)$, which
is the identity for $n<0$. The operators $S_n$ also satisfy the relations
\begin{equation}\label{Snm}
S_n S_m = P_{-n} S_{n+m} .
\end{equation}
\end{lem}
\proof First notice that the $S_k$ satisfy
\begin{equation}\label{Squotstar}
(S_n^* f)(k)=\left\{\begin{array}{ll} f(k-n) & k-n \geq 0 \\[2mm]
0 & k-n< 0. \end{array}\right.
\end{equation}
In fact, we have
$$ \langle S_n f,\psi \rangle=\sum_{k\in\N\cup\{ 0 \}} f(k+n)
\chi_{[0,\infty)}(k+n) \psi(k) $$
$$ =\sum_{u\in\Z} f(u) \chi_{[0,\infty)}(u) \chi_{[0,\infty)}(u-n)
\psi(u-n) =\sum_{u\in\N\cup\{ 0 \}} f(u) \chi_{[0,\infty)}(u-n)
\psi(u-n). $$
Thus, we have $S_n^* =S_{-n}$. we then have
$$ S_n^* S_n f (u)=\chi_{[n,\infty)}(k) f(k) =P_{n} f (k) $$
and
$$ S_n S_n^* f(k)= \chi_{[-n,\infty)}(k) f(k) =P_{-n} f (k). $$
This is in fact a particular case of the following relations.
The relation $U_n U_m=U_{n+m}$ satisfied by the unitaries acting on
$\ell^2(\Z)$ descends to the relation \eqref{Snm} between the
isometries $S_n$ acting on $\ell^2(\N\cup\{0\})$. In fact, we have
$$
S_n P_m f(k) = \chi_{[-n,\infty)}(k)\, \chi_{[m-n,\infty)}(k) f(k+n) =
P_{m-n} S_n f(k)
$$
where
$$ P_m S_n f(k) = \chi_{[-n,\infty)}(k) \, \chi_{[m,\infty)}(k)\, f(k+n). $$
Thus, we see that $S_n S_m= P_{-n} S_{n+m}$, since
$$ (S_n S_m f)(k)= \chi_{[-n,\infty)}(k) \, \chi_{[-(n+m),\infty)}(k)
f(k+m+n) = (P_{-n}S_{n+m} \, f)(u) . $$
\endproof
Thus, we see that, even in the case of a commutative group like $\Z$,
where the algebra of the $U_m$ is commutative, we obtain a
noncommutative algebra of isometries $S_m$,
$$ S_m S_n = P_{-m} S_{n+m} \neq  P_{-n} S_{n+m} = S_n S_m. $$
Notice however that, if $n$ and $m$ are both positive, then
$P_{-m}=1=P_{-n}$ so that $S_n S_m=S_m S_n=S_{n+m}$.
Notice also that the fact that the algebra generated by the isometries
$S_n$ is associative follows from the fact that the projections $P_n$
commute among themselves, as they are given by multiplication
operators by the characteristic functions $\chi_{[n,\infty)}$. In
fact, we have
$$ S_n (S_m S_k)= S_n P_{-m} S_{m+k} = P_{-m-n} S_n S_{m+k} =
P_{-m-n} P_{-n} S_{n+m+k} $$
$$
(S_n S_m) S_k = P_{-n} S_{n+m} S_k = P_{-n} P_{-(n+m)} S_{n+m+k},
$$
with $P_{-(n+m)} P_{-n}=P_{-n} P_{-(n+m)}$.
For $n>0$, we also have $S_{-1}^n = P_1 S_{-2} S_{-1}^{n-1}=P_1\cdots
P_{n-1} S_{-n}=P_{n-1} S_{-n}=S_{-n}$, since $P_{n-1} S_{-n} f(k)=
\chi_{[n-1,\infty)}(k) \chi_{[n,\infty)}(k) f(k-n)=\chi_{[n,\infty)}(k)
f(k-n)$. Similarly, $S_1^n=P_{-1} S_2
S_1^{n-1}=P_{-1}\cdots P_{-n+1} S_n=S_n$ since $P_{-n+k}=1$, for
$k=0,\ldots,n-1$.

Clearly, the algebra of the $S_n$ we described here is generated by
a single isometry $S_{-1}$, which is the isometry that describes the
``phase'' part of the creation operator in quantum mechanics,
see \cite{Feyn}.
In fact, recall that the creation and annihilation operators $a^\dag$
and $a$, with $a^*=a^\dag$, act on $\ell^2(\N\cup\{ 0\})$ by
\begin{equation}\label{aadag}
a^\dag e_n = \sqrt{n+1}\, e_{n+1} \ \ \ \text{ and } \ \ \ a\,e_n
=\sqrt{n}\, e_{n-1},
\end{equation}
with the commutation relation $[a^\dag,a]=1$. It is well known that,
while the operators $a^\dag$ and $a$ do not have a polar decomposition
in terms of a unitary and a self-adjoint operator, they have a
decomposition in terms of an isometry and a self-adjoint operator in
the form
\begin{equation}\label{polaradaga}
a^\dag = N^{1/2} S_{-1} \ \ \ \text{ and } \ \ \ a= S_1 N^{1/2},
\end{equation}
where $N\,e_n=n\, e_n$ is the grading operator on
$\ell^2(\N\cup\{0\})$ and $S_1$ and $S_{-1}$ are the isometries
described above, $S_1 e_n =e_{n-1}$ and $S_{-1} e_n=e_{n+1}$.
Notice that the grading operator $N$ acting on $\ell^2(\N\cup\{0\})$
defines a $1^+$-summable self-adjoint operator with
compact resolvent and with the
property that the commutators with the operators $S_n$ are bounded.
Namely, one has the commutation relation
\begin{equation}\label{Nboundcomm}
[N,S_n] = -n S_n .
\end{equation}
The Hamiltonian associated to the creation and annihilation operators
is of the form
\begin{equation}\label{Hamilt}
H = a^\dag a , \ \ \ \text{ with } \ \ \ \Sp(H)=\N\cup \{ 0 \}.
\end{equation}
The corresponding partition function at inverse temperature $\beta>0$
is of the form
\begin{equation}\label{Zbeta}
Z(\beta)=\Tr(e^{-\beta H})=\sum_{n\geq 0} \exp(-\beta n)=
(1-\exp(-\beta))^{-1}.
\end{equation}
%----------------------------------------------------------------------------------------------------------------------------------------------------------------------
\section{A quick introduction to Dirac operators}

\subsection{ Clifford Algebra}
Let $V \simeq \R^n$ be a vector space with non degenerate symmetric bilinear
form $g$. Over a field of characteristic different than $2$, such a bilinear form
can always be determined by the corresponding {\em quadric form} $q$ defined
as $q(v)=g(v,v)$, by setting $2g(u,v)=q(u,v)-q(u)-q(v)$.

\begin{defn}
The Clifford Algebra $CL(V,g) $ is an algebra over $\R$ generated  by the vectors
$v \in V$, subject to the relation $uv+vu=2g(u,v) $ for all $u,v \in V$.
\end{defn}

\subsubsection{ concepts in Riemannian geometry}
Let $M$ be a compact smooth $n$-dimensional manifold without boundary.
Define a Riemannian metric on $M$ to be a symmetric bilinear form
$$ g:\mathfrak{F}(M)\times \mathfrak{F}(M) \longrightarrow C(M),$$
where $\mathfrak{F}(M)=\Gamma (M,T_{\mathbb{C}}M)$ is the space of
continuous vector fields on $M$, and $C(M)$ is the
commutative $C^*$-algebra of continuous functions on $M$. Then $g$ satisfies
the following properties.
\begin{enumerate}
\item $g(X,Y)$ is a real function if $X,Y$ are real vector fields.
\item $g$ is $C(M)$-bilinear \ie $ g(fX,Y)=g(X,fY)=fg(X,Y)$ for all $f \in C(M)$.
In this condition $g$ is given by a continuous family of symmetric
bilinear map $g_x:T_xM \times T_xM \longrightarrow R$ , where $g_x$ is positive definite.
\item $g(X,X) \geq 0$ for $X$ real, with $g(X,X)=0$ only if $X=0$ in $\mathfrak{F}(M)$.
\end{enumerate}

\subsubsection{ Dirac Operator}
Let $(M,g)$ be a smooth compact Riemannian $m$-manifold without boundary
with a Clifford algebra bundle $CL(M)$. A Clifford module is a module over $CL(M)$.
Any Clifford module $\Lambda$ that is finitely generated and projective is of the form
$\Lambda = \Gamma (M,E)$ for $E \longrightarrow M$ a complex vector bundle.
For $ E \longrightarrow M$ a smooth complex vector bundle of Clifford modules
$\Lambda =\Gamma (M,E)$, we can define the Clifford multiplication which is
a bundle map $c:CL(M) \longrightarrow Hom(E,E)$ which is given fiberwise by maps
$c: CL(T_xM,g_x) \longrightarrow Hom_{\mathbb{C}} (E_x,E_x)$.

Any choice of a smooth connection
$$ \nabla :C^\infty (M,E) \longrightarrow C^\infty (M,T^*M \otimes E) $$
defines an operator of {\em Dirac type} by setting $\mathcal{D}=c \circ \nabla$.
We use here the identification of tangent and cotangent bundle $TM\cong T^*M$
induced by the Riemannian metric.

Consider a small open chart domain $U \subset M$, where the cotangent
bundle is trivial, \ie $ T^*M \mid_U \simeq U\times \mathbb{R}^n$. For any local coordinate $(x_1,x_2,...,x_n)$
over a chart domain $U$, the local coordinate of the cotangent bundle  $ T^*M \mid_U $ are
$(x,\xi)=(x_1,x_2,...,x_n,\xi_1,\xi_2,...,\xi_n)$ where $\xi \in T^*_x M$.
A differential operator acting on smooth local sections $f \in \Gamma(U,E)$ is an operator $P$ of the form
$$ P=\sum_{\mid \alpha\mid \leq d} a_\alpha(x) D^\alpha$$
with $a_\alpha \in \Gamma (U,End E)$, and where $D^\alpha = D^{\alpha_1}_1.D^{\alpha_2}_2.D^{\alpha_3}_3...D^{\alpha_n}_n. $
and $D_j=-i\frac{\partial}{\partial x_j}$, with a positive integer $d$ representing the order of $P$.
Let $E \longrightarrow M$ be a vector bundle of rank $r$. By Fourier transform, we can write
for $f\in C^\infty (U,\mathbb{R}^n)$:
\begin{eqnarray}\label{PP1}
Pf(x)&=&(2\pi)^n \int_{\mathbb{R}^n}\exp^{ix\xi}p(x,\xi)\hat{f}(\xi)d^n\xi \nonumber\\
&=&(2\pi)^n  \int\int_{\mathbb{R}^{2n}}\exp^{i(x-y)\xi}p(x,\xi)f(\xi)d^ny d^n\xi ,
\end{eqnarray}
where $p(x,\xi)$ is a polynomial of order $d$ in the $\xi$-variable, called the {\em complete symbol}
of the operator $P$. Then we can isolate the homogeneous part
$$p(x,\xi)=\sum_{j=0}^d p_{d-j}(x,\xi)$$
where $p_{d-j}(x,t\xi)=t^{d-j}p_{d-j}(x,\xi) $ for $ t>0$.
\begin{defn}
An element $p$ is called a {\em classical symbol} if we can find a sequence of terms
$p_{d}(x,\xi), p_{d-1}(x,\xi),....$ with
$$p(x,\xi)\sim \sum_{j>0}^\infty p_{d-j}(x,\xi)$$
such that $p_{d-j}(x,t\xi)=t^{d-j}p_{d-j}(x,\xi)$ for $ t>0$.
\end{defn}

\begin{defn}
A classical {\em pseudo-differential operator} of order $d$
over the chart domain $U \subset \mathbb{R}^n$
is an operator $P$ defined by \ref{PP1}, for which $p(x,\xi)$ is a classical symbol, whose leading term
$p_d(x,\xi)$ does not vanish. This leading term is called the {\em principal symbol} of $P$, and we also
denote it by $\sigma(P)(x,\xi)=p_d(x,\xi)$.
\end{defn}

For an operator of Dirac type, which is a first-order differential operator on $\Gamma(M,E)$, we get
$\sigma(\mathcal{D}) \in \Gamma(T^*M,\pi^*(End E))$ and for the property of
$p(x,\xi)=c(dx^j)(\xi_j-i\omega_j(x))$ we get
$$\sigma (\mathcal{D})(x,\xi)=c(\xi_jdx^j)=c(\xi)$$
and $$\sigma (\mathcal{D}^2)(x,\xi)=(\sigma (\mathcal{D})(x,\xi))^2=(c(\xi))^2=g(\xi,\xi)$$
Notice that $\sigma (\mathcal{D}^2)$ only vanishes when $\xi=0$, that is, on the zero section
of $T^*M$.

\begin{defn}
Let $P$ a classical pseudo-differential operator, then $P$ is called {\em elliptic}
if $\sigma (P)(x,\xi)$ is invertible when $\xi \neq 0$.
\end{defn}

An operators $\cD$ of Dirac type is elliptic and so is its
square $\cD^2$. On a compact manifold without boundary this
implies that it is Fredholm (has finite dimensional kernel
and cokernel), hence its index $Ind(\cD)=\dim Ker(\cD)-
\dim Coker(\cD)$ is well defined. The Atiyah-Singer index
theorem gives a local formula, in terms of integration of
a differential form, for the index. On a manifold with
boundary, the Fredholm property depends on the choice of
boundary conditions. With the Atiyah-Patodi-Singer boundary
conditions one still has a Fredholm operator and an
index formula, now with an additional term that is an
eta invariant for the operator restricted to the boundary
manifold.

%-----------------------------------------------------------------------------------------------------------------------------------------------
\section{Concepts of Cyclic Cohomology}
Cyclic cohomology of non-commutative algebras is playing in non-commutative geometry a similar
role to that of de Rham cohomology in differential topology. The first appearance of the Cyclic cohomology
was in the cohomology theory for algebras. The cyclic cohomology $HC^*(\mathcal{A})$ of an algebra $\mathcal{A}$
over $\mathbb{R} $ or $ \mathbb{C}$ is the cochain complex $\{ C^*_\lambda(\mathcal{A}),b \}$,
where $C^*_\lambda(\mathcal{A}), n\geq 0$ consists of the $(n+1)$-linear forms $\vartheta$ on $\mathcal{A}$
satisfying the cyclicity condition
\begin{eqnarray}\label{Cyc}
\vartheta (a^0,a^1,...a^n)=(-1)^n \vartheta(a^1,a^2,...,a^0)
\end{eqnarray}
where $ a^0,a^1,...,a^n \in \mathcal{A}$ and the coboundary operator is given by
$$(b\vartheta)(a^0,a^1,...,a^n)=\sum^n_{j=0}(-1)^{j}(a^0,...,a^ja^{j+1},...,a^{n+1})+(-1)^{n+1}\vartheta (a^{n+1}a^0,...,a^n) $$
$C^*_\lambda(\mathcal{A})$ then consists of all continuous $(n+1)$-linear forms on $\mathcal{A}$ satisfying \ref{Cyc}.
Cyclic cohomology provides numerical invariants of K-theory classes as follows.
For an even integer $n$, given an $n$-dimensional cyclic cocycle $\vartheta$ on $\mathcal{A}$, then
the scalar
\begin{eqnarray}\label{Tr}
\vartheta \otimes Tr(E,E,...,E)
\end{eqnarray}
is invariant under homotopy, for an idempotent $$E^2=E\in M_N(\mathcal{A})=\mathcal{A}\otimes M_N(\mathbb{C})$$
This gives the pairing $\langle[\vartheta],[E]\rangle$ between cyclic homology and K-theory.
For a manifold $M$ let $\mathcal{A}=C^\infty(M)$ with
$$\vartheta(f^0,f^1,...,f^n)=\langle\Omega,f^0df^1\wedge df^2\wedge...\wedge df^n\rangle$$
where $f^1,f^2,...,f^n \in \mathcal{A}$
and $\Omega$ is a closed $n$-dimensional de Rham form on $M$. Then the invariant \ref{Tr} up to
normalization is equal to $\langle\Omega,ch^*(\tau)\rangle$ where $ch^*(\tau)$ denotes the Chern character
of the rank $N$ vector bundle $\tau$ on $M$ whose fiber at $x \in M$ is the range of $E(x) \in M_N(\mathbb{C})$.
To any algebra $\mathcal{A}$ one can associate a module $\mathcal{A}^\natural$ over the cyclic category
by assigning to each integer $n \geq 0$ the vector space $C^n(\mathcal{A})$ of $(n+1)$-linear forms
$\vartheta(a^0,a^1,...,a^n)$ on $\mathcal{A}$ and to the generating morphisms the operators
$\delta_i:C^{n-1} \longrightarrow C^n$ and $\sigma_i:C^{n+1} \longrightarrow C^n$ defined above.
One thus, obtains the desired interpretation of the cyclic cohomology group of $\mathcal{K}$-algebra
$\mathcal{A}$ over a ground ring $\mathcal{K}$ in terms of derived functors over the cyclic category
$$HC^n(\mathcal{A})\simeq Ext^n_{\Lambda}(\mathcal{K}^\natural,\mathcal{A}^\natural)$$
and
$$HC_n(\mathcal{A})\simeq Tor_n^{\Lambda}(\mathcal{A}^\natural,\mathcal{K}^\natural)$$

%---------------------------------------------------------------
%\bibliographystyle{amsplain}
%\bibliography{}

\begin{thebibliography}{33333}

\bibitem{co} A.~Connes, C.~Consani, M.~Marcolli, {\em
Noncommutative geometry and motives: the thermodynamics of
endomotives}, to appear in Advances in Mathematics, math.QA/0512138.
%18

\bibitem{CC} A.~Chamseddine, A.~Connes, {\em Quantum gravity
boundary terms from spectral action of noncommutative space.}
Phys. Rev. Lett. 99, 071302 (2007).
%22
\bibitem{Co94} A.~Connes, {\em Noncommutative Geometry}. Academic Press, 1994.
%23
\bibitem{Co-ext}  A.~Connes, {\em Cohomologie cyclique et foncteurs ${\rm Ext}\sp n$}.
C. R. Acad. Sci. Paris S\'er. I Math. Vol.296 (1983), N.23, 953--958.
%24
\bibitem{Co} A.~Connes, Lecture at Oberwolfach, September 2007.
%25
\bibitem{co1} A.~Connes, M.~Marcolli, {\em Noncommutative geometry, quantum fields and motives}, Colloquium Publications, Vol.55,
American Mathematical Society, 2008.
%26
\bibitem{CoSka}  A.~Connes, G.~Skandalis, {\em The longitudinal index theorem for foliations}.  Publ. Res. Inst. Math. Sci.  20  (1984),
no. 6, 1139--1183.
%63

\bibitem{Cu} J.~Cuntz, {\em Cyclic theory and the bivariant Chern-Connes
character}, in  ``Noncommutative geometry'',  pp.73--135, Lecture Notes
in Math., 1831, Springer, Berlin, 2004.
%60

\bibitem{Nistor} V.~Nistor, {\em A bivariant Chern--Connes character},
Annals of Math. Vol.138 (1993) 555--590.
%61




\bibitem{Nistor0} Connes, Alain (1994), {\em Non-commutative geometry}, Boston, MA: Academic Press, ISBN 978-0-12-185860-5

\bibitem{Nistor1} Connes, Alain; Marcolli, Matilde (2008), {\em A walk in the noncommutative garden}, An invitation to noncommutative geometry, World Sci. Publ., Hackensack, NJ, pp. 1–128, arXiv:math/0601054

\bibitem{Nistor2} Connes, Alain; Marcolli, Matilde (2008), {\em Noncommutative geometry, quantum fields and motives}, American Mathematical Society Colloquium Publications, 55, Providence, R.I.: American Mathematical Society, ISBN 978-0-8218-4210-2, MR 2371808

\bibitem{Nistor3} Gracia-Bondia, Jose M; Figueroa, Hector; Varilly, Joseph C (2000), {\em Elements of Non-commutative geometry}, Birkhauser, ISBN 978-0-8176-4124-5

\bibitem{Nistor4} Landi, Giovanni (1997), {\em An introduction to noncommutative spaces and their geometries}, Lecture Notes in Physics. New Series m: Monographs, 51, Berlin, New York: Springer-Verlag, pp. arXiv:hep–th/9701078, arXiv:hep-th/9701078, ISBN 978-3-540-63509-3, MR 1482228

\bibitem{Nistor5} Van Oystaeyen, Fred; Verschoren, Alain (1981), {\em Non-commutative algebraic geometry}, Lecture Notes in Mathematics, 887, Springer-Verlag, ISBN 978-3-540-11153-5

\bibitem{Nistor6} Grensing, Gerhard (2013). {\em Structural Aspects of Quantum Field Theory and Noncommutative Geometry}. World Scientific. ISBN 978-981-4472-69-2.

\bibitem{Nistor7} M.R. Douglas and N. A. Nekrasov (2001) {\em Noncommutative field theory}, Rev. Mod. Phys. 73: 977–1029.

\bibitem{Nistor8} Szabo, R. (2003) {\em Quantum Field Theory on Noncommutative Spaces}, Physics Reports 378: 207-99.

\bibitem{Nistor9} V. Moretti (2003), {\em Aspects of noncommutative Lorentzian geometry for globally hyperbolic spacetimes}, " Rev. Math. Phys. 15: 1171-1218.


\bibitem{Nistor10} Arveson, W. (1976), {\em An Invitation to C*-Algebra, Springer-Verlag}, ISBN 0-387-90176-0. An excellent introduction to the subject, accessible for those with a knowledge of basic functional analysis.

\bibitem{Nistor11} Dixmier, Jacques (1969), {\em Les C*-algèbres et leurs représentations}, Gauthier-Villars, ISBN 0-7204-0762-1. This is a somewhat dated reference, but is still considered as a high-quality technical exposition. It is available in English from North Holland press.

\bibitem{Nistor12} Doran, Robert S.; Belfi, Victor A. (1986), {\em Characterizations of C*-algebras: The Gelfand-Naimark Theorems}, CRC Press, ISBN 978-0-8247-7569-8.


\bibitem{Nistor 13} Emch, G. (1972),{\em  Algebraic Methods in Statistical Mechanics and Quantum Field Theory}, Wiley-Interscience, ISBN 0-471-23900-3.

\bibitem{Nistor14} A.I. Shtern (2001) [1994], {\em C* algebra}, in Hazewinkel, Michiel, Encyclopedia of Mathematics, Springer Science Business Media B.V. / Kluwer Academic Publishers, ISBN 978-1-55608-010-4

\bibitem{Nisto15r} Sakai, S. (1971), {\em C*-algebras and W*-algebras}, Springer, ISBN 3-540-63633-1.

\bibitem{Nistor16} Segal, Irving (1947), {\em Irreducible representations of operator algebras}, Bulletin of the American Mathematical Society, 53 (2): 73–88, doi:10.1090/S0002-9904-1947-08742-5.

\bibitem{Nistor17} Isaacs, I. Martin (1993), {\em Algebra, a graduate course (1st ed.)}, Brooks/Cole, ISBN 0-534-19002-2

\bibitem{Nistor18} Jacobson, N. (1945), {\em Structure theory of simple rings without finiteness assumptions}, Transactions of the American Mathematical Society, 57: 228–245, doi:10.2307/1990204, JSTOR 1990204

\bibitem{MMAZ} Matilde Marcolli and Ahmad Zainy Al-Yasry, {\em Coverings, correspondences, and noncommutative geometry}
Journal of Geometry and Physics 58 (2008) no. 12, 1639–1661, arXiv:0807.2924
\bibitem{MMAADD} Domenic Denicola, Matilde Marcolli and Ahmad Zainy Al-Yasry, {\em Spin foams and noncommutative geometry}, Class. Quantum Grav. 27 (2010) 205025 (53pp), arXiv:1005.1057

\end{thebibliography}
\end{document}